\documentclass[11pt]{article}

\usepackage{amsfonts}
\usepackage{amscd}
\usepackage{amssymb}
\usepackage{amsthm}
\usepackage{amsmath, xspace}
\usepackage{blkarray}
\usepackage{color}
\usepackage{graphics}
\usepackage{graphicx}
\usepackage{enumitem}
\usepackage{fancyhdr}
\usepackage{mathdots}
\usepackage{mathrsfs}
\usepackage{multicol}
\usepackage{stmaryrd}
\usepackage{youngtab}
\usepackage[all]{xy}

\usepackage[plainpages,backref]{hyperref}

\usepackage{lscape}

\theoremstyle{plain}
\newtheorem{thm}{Theorem}[section]

\newtheorem{prop}[thm]{Proposition}

\theoremstyle{definition}
	
\newtheorem{remark}[thm]{Remark}
\theoremstyle{example}

\theoremstyle{remark}

\numberwithin{equation}{section}

\providecommand{\keywords}[1]{\textbf{\textit{Key words---}} #1}

\def\cB{\mathcal{B}}

\def\cL{\mathcal{L}}

\def\CC{\mathbb{C}}

\def\FF{\mathbb{F}}

\def\kk{\Bbbk}

\def\RR{\mathbb{R}}

\def\ZZ{\mathbb{Z}}

\def\fa{\mathfrak{a}}

\def\fg{\mathfrak{g}}
\def\fh{\mathfrak{h}}

\def\fo{\mathfrak{o}}

\def\fsl{\mathfrak{sl}}

\def\Card{\mathrm{Card}}
\def\diag{\mathrm{diag}}

\def\Ind{\mathrm{Ind}}

\def\Res{\mathrm{Res}}

\def\wt{\mathrm{wt}}

\def\linspan{\textnormal{-span}}

\def\al{\alpha}

\usepackage{etex} 
\usepackage{pictexwd}
\usepackage{tikz}
\usepgflibrary{shapes.geometric}
\usetikzlibrary{arrows, calc, positioning}

\tikzstyle{V}=[draw, fill =black, circle, inner sep=0pt, minimum size=1.5pt]
\tikzstyle{wV}=[draw, fill =white, circle, inner sep=0pt, minimum size=4.5pt]
\tikzstyle{bV}=[draw, fill =black, circle, inner sep=0pt, minimum size=4.5pt]
\tikzstyle{over}=[draw=white,double=black,line width=2pt, double distance=.5pt]

\def\Over[#1,#2][#3,#4]{ 
	\draw[style=over]   (#1,#2) .. controls ++(0,#4*.5-#2*.5) and ++(0,-#4*.5+#2*.5) .. (#3,#4);}
\def\Under[#1,#2][#3,#4]{ 
	\draw  (#1,#2) .. controls ++(0,#4*.5-#2*.5) and ++(0,-#4*.5+#2*.5) .. (#3,#4);}
\def\Cross[#1,#2][#3,#4]{
	\Under[#3,#2][#1,#4]\Over[#1,#2][#3,#4]}

\def\Tops[#1][#2][#3]{
	\foreach\x in {#1}{
		\draw (\x+.1,#2) -- (\x+.1,#2+.15) (\x-.1,#2) -- (\x-.1,#2+.15) ;
		\draw (\x+.1,#2+.15) arc (0:360:1mm and .5mm);}
	\foreach \x in {1,...,#3} {\draw (\x,#2)  to (\x,#2+.05) node[V]{};}
	}
\def\Bottoms[#1][#2][#3]{
	\foreach\x in {#1}{
		\draw (\x+.1,#2) -- (\x+.1,#2-.1) (\x-.1,#2) -- (\x-.1,#2-.1) ;
		\draw (\x+.1,#2-.1) arc (0:-180:1mm and .5mm);}
	\foreach \x in {1,...,#3} {\draw (\x,#2)  to (\x,#2-.05) node[V]{};}
	}
\def\Caps[#1][#2,#3][#4]{
	\Tops[#1][#3][#4]
	\Bottoms[#1][#2][#4]
	}
\def\Pole[#1][#2,#3]{
	\shade[left color=white,right color=white] (#1+.1,#2) rectangle (#1-.1,#3);
	\draw[over] (#1+.1,#2) to (#1+.1,#3) (#1-.1,#2) to (#1-.1,#3) ;}
\def\Label[#1,#2][#3][#4]{
	\node[above] at (#3,#2+.1) {#4};
	\node[below] at (#3,#1-.1) {#4};		}

\newcommand{\posleq}[1]{
	\hspace{0.1cm}
	\begin{tikzpicture}
	\draw (-0.8ex, -0.5ex) -- (0.8ex, -0.5ex);
	\draw (-0.8ex, 0.4ex) -- (0.7ex, -0.2ex);
	\draw (-0.8ex, 0.4ex) -- (0.7ex, 1ex);
	\draw (0.4ex,0.4ex) --(1.1ex, 0.4ex);
	\draw (0.75ex,0.75ex) --(0.75ex, 0.05ex);
	\end{tikzpicture}
	\hspace{0.1cm}
	}
\newcommand{\negleq}[1]{
	\hspace{0.1cm}
	\begin{tikzpicture}
	\draw (-0.8ex, -0.5ex) -- (0.8ex, -0.5ex);
	\draw (-0.8ex, 0.4ex) -- (0.7ex, -0.2ex);
	\draw (-0.8ex, 0.4ex) -- (0.7ex, 1ex);
	\draw (0.4ex,0.4ex) --(1.1ex, 0.4ex);
	\end{tikzpicture}
	\hspace{0.1cm}
	}
	
\newcommand{\zeroleq}[1]{
	\hspace{0.1cm}
	\begin{tikzpicture}
	\draw (-0.8ex, -0.5ex) -- (0.8ex, -0.5ex);
	\draw (-0.8ex, 0.4ex) -- (0.7ex, -0.2ex);
	\draw (-0.8ex, 0.4ex) -- (0.7ex, 1ex);
	\draw  (0.75ex,0.4ex) ellipse (0.2ex and 0.35ex);
	\end{tikzpicture}
	\hspace{0.1cm}
	}
	
\newcommand{\posgeq}[1]{
	\hspace{0.1cm}
	\begin{tikzpicture}
	\draw (-0.8ex, -0.5ex) -- (0.8ex, -0.5ex);
	\draw (0.8ex, 0.4ex) -- (-0.7ex, -0.2ex);
	\draw (0.8ex, 0.4ex) -- (-0.7ex, 1ex);
	\draw (-0.4ex,0.4ex) --(-1.1ex, 0.4ex);
	\draw (-0.75ex,0.75ex) --(-0.75ex, 0.05ex);
	\end{tikzpicture}
	\hspace{0.1cm}
	}
\newcommand{\neggeq}[1]{
	\hspace{0.1cm}
	\begin{tikzpicture}
	\draw (-0.8ex, -0.5ex) -- (0.8ex, -0.5ex);
	\draw (0.8ex, 0.4ex) -- (-0.7ex, -0.2ex);
	\draw (0.8ex, 0.4ex) -- (-0.7ex, 1ex);
	\draw (-0.4ex,0.4ex) --(-1.1ex, 0.4ex);
	\end{tikzpicture}
	\hspace{0.1cm}
	}
	
\newcommand{\zerogeq}[1]{
	\hspace{0.1cm}
	\begin{tikzpicture}
	\draw (-0.8ex, -0.5ex) -- (0.8ex, -0.5ex);
	\draw (0.8ex, 0.4ex) -- (-0.7ex, -0.2ex);
	\draw (0.8ex, 0.4ex) -- (-0.7ex, 1ex);
	\draw  (-0.75ex,0.4ex) ellipse (0.2ex and 0.35ex);
	\end{tikzpicture}
	\hspace{0.1cm}
	}

\newcommand{\posl}[1]{
	\hspace{0.1cm}
	\begin{tikzpicture}
	\draw (-0.8ex, 0.4ex) -- (0.7ex, -0.2ex);
	\draw (-0.8ex, 0.4ex) -- (0.7ex, 1ex);
	\draw (0.4ex,0.4ex) --(1.1ex, 0.4ex);
	\draw (0.75ex,0.75ex) --(0.75ex, 0.05ex);
	\end{tikzpicture}
	\hspace{0.1cm}
	}
\newcommand{\negl}[1]{
	\hspace{0.1cm}
	\begin{tikzpicture}
	\draw (-0.8ex, 0.4ex) -- (0.7ex, -0.2ex);
	\draw (-0.8ex, 0.4ex) -- (0.7ex, 1ex);
	\draw (0.4ex,0.4ex) --(1.1ex, 0.4ex);
	\end{tikzpicture}
	\hspace{0.1cm}
	}
	
\newcommand{\zerol}[1]{
	\hspace{0.1cm}
	\begin{tikzpicture}
	\draw (-0.8ex, 0.4ex) -- (0.7ex, -0.2ex);
	\draw (-0.8ex, 0.4ex) -- (0.7ex, 1ex);
	\draw  (0.75ex,0.4ex) ellipse (0.2ex and 0.35ex);
	\end{tikzpicture}
	\hspace{0.1cm}
	}
	
\newcommand{\posg}[1]{
	\hspace{0.1cm}
	\begin{tikzpicture}
	\draw (0.8ex, 0.4ex) -- (-0.7ex, 1ex);
	\draw (0.8ex, 0.4ex) -- (-0.7ex, -0.2ex);
	\draw (-0.4ex,0.4ex) --(-1.1ex, 0.4ex);
	\draw (-0.75ex,0.75ex) --(-0.75ex, 0.05ex);
	\end{tikzpicture}
	\hspace{0.1cm}
	}
\newcommand{\negg}[1]{
	\hspace{0.1cm}
	\begin{tikzpicture}
	\draw (0.8ex, 0.4ex) -- (-0.7ex, -0.2ex);
	\draw (0.8ex, 0.4ex) -- (-0.7ex, 1ex);
	\draw (-0.4ex,0.4ex) --(-1.1ex, 0.4ex);
	\end{tikzpicture}
	\hspace{0.1cm}
	}
	
\newcommand{\zerog}[1]{
	\hspace{0.1cm}
	\begin{tikzpicture}
	\draw (0.8ex, 0.4ex) -- (-0.7ex, -0.2ex);
	\draw (0.8ex, 0.4ex) -- (-0.7ex, 1ex);
	\draw  (-0.75ex,0.4ex) ellipse (0.2ex and 0.35ex);
	\end{tikzpicture}
	\hspace{0.1cm}
	}

\makeatletter
\renewcommand{\@makefnmark}{\mbox{\textsuperscript{}}}
\makeatother

\title{Positive level, negative level and level zero}
\author{
Finn McGlade\quad\ \ email:\ fmcglade@ucsd.edu \\
Arun Ram\quad\ \ email:\ aram@unimelb.edu.au \\
Yaping Yang\quad\ \ email:\ yaping.yang1@unimelb.edu.au 
\\
}
\date{\vspace{-5ex}}

\usetikzlibrary{arrows.meta}

\begin{document}

\maketitle

\begin{center}
{\sl Dedicated to I.G.\ Macdonald and A.O.\ Morris}
\end{center}

\begin{abstract}
\noindent
This is a survey on the combinatorics and
geometry of integrable representations of quantum affine algebras 
with a particular focus on level 0. Pictures and examples are included
to illustrate the affine Weyl group orbits, crystal graphs and Macdonald polynomials
that provide detailed understanding of the structure of the extremal weight modules 
and their characters.
The final section surveys the alcove walk method of working with the
positive level, negative level and level zero affine flag varieties and describes
the corresponding actions of the affine Hecke algebra.
\end{abstract}

\footnote{\today \ \ AMS Subject Classifications: Primary 17B37; Secondary 17B67.}

\setcounter{section}{-1}

\section{Introduction}

This paper is about positive level, negative level and level 0.
It was motivated by the striking result of \cite{KNS17}, which establishes a 
Pieri-Chevalley formula for the K-theory of the semi-infinite (level 0) affine flag variety.
This made us want to learn more about the level 0 integrable modules of quantum affine
algebras.  Our trek brought us face to face with a huge literature, 
including important contributions from
Drinfeld, Kashiwara, Beck, Chari, Nakajima, 
Lenart-Schilling-Shimozono, Cherednik-Orr, Naito-Sagaki, Feigin-Makedonskyi, Kato, their coauthors 
and many others.  It is a beautiful theory and we count ourselves lucky to have been drawn into it.

The main point is that the integrable modules for quantum affine algebras $\mathbf{U}$
naturally partition themselves into families: positive level, negative level and level zero.
Their structure is shadowed by the orbits of the affine Weyl group on the lattice
of weights for the affine Lie algebra, which take the shape of a concave up paraboloid at
positive level, a concave down paraboloid at negative level and a tube at level 0.
These integrable modules have crystal bases which provide detailed control of their characters.
At level 0 the characters are (up to a factor similar to a Weyl denominator) 
Macdonald polynomials specialised at $t=0$.  The next amazing feature is that there are 
Borel-Weil-Bott theorems for each case: positive level, negative level and level 0,
where, respectively, the appropriate geometry is a positive level (thin) affine flag
variety, a negative level (thick) affine flag variety, and a level zero (semi-infinite) 
affine flag variety.

This paper is a survey of the general picture of positive level, negative level, and level zero,
in the context of the combinatorics of affine Weyl groups and crystals, of the representation
theory of integrable modules for quantum affine algebras, and of the geometry of
affine flag varieties.  In recent years, the picture has become more and more rich and 
taken clearer focus.
We hope that this paper will help to bring this story to a wider audience by providing
pictures and some explicit small examples for $\widehat{\fsl}_2$ and $\widehat{\fsl}_3$.

\subsection{Orbits of the affine Weyl group $W$ action on $\fh^*$}
\label{Worbits}

For $\widehat{\fsl}_2$ the vector space
$\fh^*$ 
is three dimensional with basis $\{ \delta, \omega_1, \Lambda_0\}$, 
and the orbits of the action of the affine
Weyl group $W^{\mathrm{ad}}$ on $\fh^*$ on different levels look like
$$
\begin{matrix}
\begin{tikzpicture}[scale=0.8]
\tikzstyle{every node}=[font=\tiny]    
    \draw[style=help lines,step=1cm] (-10,-5) grid (10, 5);   
        \draw[->] (0,-5) -- (0, 6) node[anchor=south]{$\Lambda_0$-axis};
    \draw[->] (-10,0) -- (10.5,0) node[anchor=west]{$\omega_1$-axis};
        \draw[-, red, very thick] (0,-5) -- (0,5.5) node[anchor=west]{$\fh^{\alpha_1^{\vee}}$};
    \draw[-, red, very thick] (-5,-5) -- (5,5) node[anchor=west]{$\fh^{\alpha_0^{\vee}}$};  
            \draw[->, blue, very thick] (0, 0) -- (1, 0);  
             \node[blue] at (1.2, 0.1) {$\omega_1$};  
            \draw[->, blue,  very thick] (0, 0) -- (-1, 0);
                \node[blue] at (-1.2, 0.1) {$-\omega_1$};  
\node at (1, 2) {$\bullet$};  
\draw[->, very thick] (0, 0) -- (1, 2) node[anchor=south]{$\Lambda$};
\node at (-1, 2) {$\bullet$};  
\draw[->, very thick] (0, 0) -- (-1, 2) node[anchor=south]{$s_1\Lambda$};
\node at (3, 2) {$\bullet$};  
 \draw[->, very thick] (0, 0) -- (3, 2) node[anchor=south]{$s_0\Lambda$};
\node at (-3, 2) {$\bullet$};  
\draw[->, very thick] (0, 0) -- (-3, 2) node[anchor=south]{$s_1s_0\Lambda$};
\node at (5, 2) {$\bullet$};  
\draw[->, very thick] (0, 0) -- (5, 2) node[anchor=south]{$s_0s_1\Lambda$};
\node at (-5, 2) {$\bullet$};  
 \draw[->, very thick] (0, 0) -- (-5, 2) node[anchor=south]{$s_1s_0s_1\Lambda$};
\node at (7, 2) {$\bullet$};  
\draw[->, very thick] (0, 0) -- (7, 2) node[anchor=south]{$s_0s_1s_0\Lambda$};
\node at (-7, 2) {$\bullet$};  
\draw[->, very thick] (0, 0) -- (-7, 2) node[anchor=south]{$s_1s_0s_1s_0\Lambda$};
\node at (9, 2) {$\bullet$};  
\draw[->, very thick] (0, 0) -- (9, 2) node[anchor=south]{$s_0s_1s_0s_1\Lambda$};
\node at (1, -2) {$\bullet$};  
\draw[->, purple, very thick] (0, 0) -- (1, -2) node[anchor=north]{$-s_1\Lambda$};
\node at (-1, -2) {$\bullet$};  
\draw[->, purple, very thick] (0, 0) -- (-1, -2) node[anchor=north]{$-\Lambda$};
\node at (3, -2) {$\bullet$};  
 \draw[->, purple, very thick] (0, 0) -- (3, -2) node[anchor=north]{$-s_1s_0\Lambda$};
\node at (-3, -2) {$\bullet$};  
\draw[->, purple, very thick] (0, 0) -- (-3, -2) node[anchor=north]{$-s_0\Lambda$};
\node at (5, -2) {$\bullet$};  
\draw[->, purple, very thick] (0, 0) -- (5, -2) node[anchor=north]{$-s_1s_0\Lambda$};
\node at (-5, -2) {$\bullet$};  
 \draw[->, purple, very thick] (0, 0) -- (-5, -2) node[anchor=north]{$-s_0s_1\Lambda$};
\node at (7, -2) {$\bullet$};  
\draw[->, purple, very thick] (0, 0) -- (7, -2) node[anchor=north]{$-s_1s_0s_1\Lambda$};
\node at (-7, -2) {$\bullet$};  
\draw[->, purple, very thick] (0, 0) -- (-7, -2) node[anchor=north]{$-s_0s_1s_0s_1\Lambda$};
\node at (9, -2) {$\bullet$};  
\draw[->, purple, very thick] (0, 0) -- (9, -2) node[anchor=north]{$-s_1s_0s_1s_0\Lambda$};
\end{tikzpicture}
\\
\hbox{
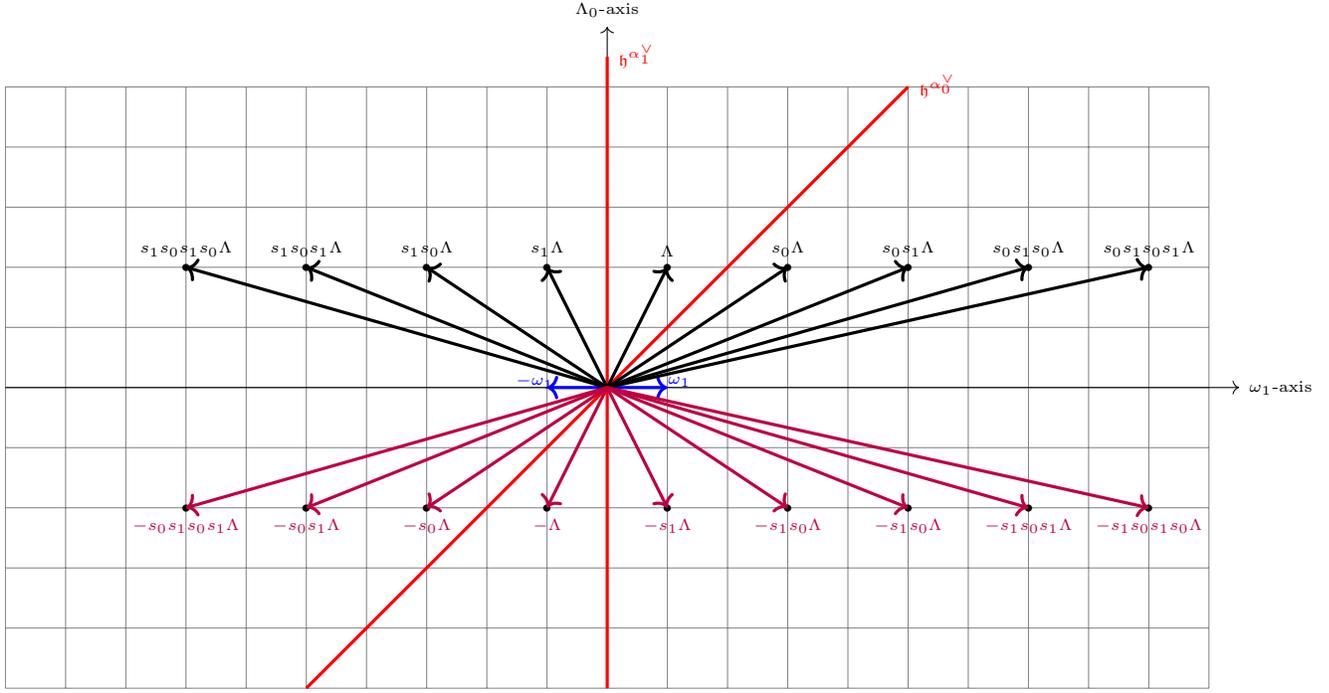
$W^{\mathrm{ad}}(\omega_1+2\Lambda_0)$ (positive level), 
$W^{\mathrm{ad}}\omega_1$ (level 0),
and $W^{\mathrm{ad}}(-\omega_1-2\Lambda_0)$ (negative level) mod $\delta$}
\end{matrix}
$$
Although informative, the picture above is misleading as it is a 
two dimensional projection (mod $\delta$) of what is actually going on.  
The $W^{\mathrm{ad}}$-action fixes the level (the coefficient of $\Lambda_0$) 
but it actually changes the $\delta$ coordinate significantly.  Let us look at the orbits
in $\delta$ and $\omega_1$ coordinates (i.e. mod $\Lambda_0$).
$$
\begin{matrix}
\begin{tikzpicture}[scale=0.7]
\tikzstyle{every node}=[font=\small]
    
    \draw[style=help lines,step=1cm] (-7.5,-6.5) grid (7.7,1.5);
    
    \draw[->] (-8,0) -- (8.5,0) node[anchor=west]{$\omega_1$-axis};
    \draw[->] (0,-6) -- (0, 1.5) node[anchor=south]{$\delta$-axis};

     \node at (1,0) {$\bullet$};
        \filldraw (1,1) node[anchor=north west,yshift=-0.1cm] {$\Lambda=\omega_1+2\Lambda_0$};
            \node at (-1,0) {$\bullet$};
 \filldraw (-1,1) node[anchor=north,yshift=-0.1cm] {$s_1\Lambda$};
     \node at (3, -1) {$\bullet$};  
    \filldraw (3, 0) node[anchor=north,yshift=-0.1cm] {$s_0\Lambda$};
         \node at (-3, -1) {$\bullet$};  
   \filldraw (-3, 0) node[anchor=north,yshift=-0.1cm] {$s_1s_0\Lambda$};
         \node at (5, -3) {$\bullet$};  
   \filldraw (5, -3) node[anchor=south west,yshift=-0.1cm] {$s_0s_1\Lambda$};
         \node at (-5, -3) {$\bullet$};  
   \filldraw (-5, -3) node[anchor=south east,yshift=-0.1cm] {$s_1s_0s_1\Lambda$};
1
         \node at (7, -6) {$\bullet$};  
   \filldraw (7, -6) node[anchor=south west,yshift=-0.1cm] {$s_0s_1s_0\Lambda$};
         \node at (-7, -6) {$\bullet$};  
   \filldraw (-7, -6) node[anchor=south east,yshift=-0.1cm] {$(s_1s_0)^2\Lambda$};

		\draw[thick] (0, 0.125)parabola (-7.54983, -7) ;
		\draw[thick] (0, 0.125)parabola (7.54983, -7) ;
\draw [blue, very thick] (1,0) to (-1,0);
\draw [red, very thick] (1,0) to (3,-1);
\draw [red, very thick] (-1,0) to (5,-3);
\draw [blue, very thick] (3,-1) to (-3,-1);
\draw [red, very thick] (-3,-1) to (7,-6);
\draw [blue, very thick] (5,-3) to (-5,-3);
\draw [blue, very thick] (7,-6) to (-7,-6);
\end{tikzpicture}
\\
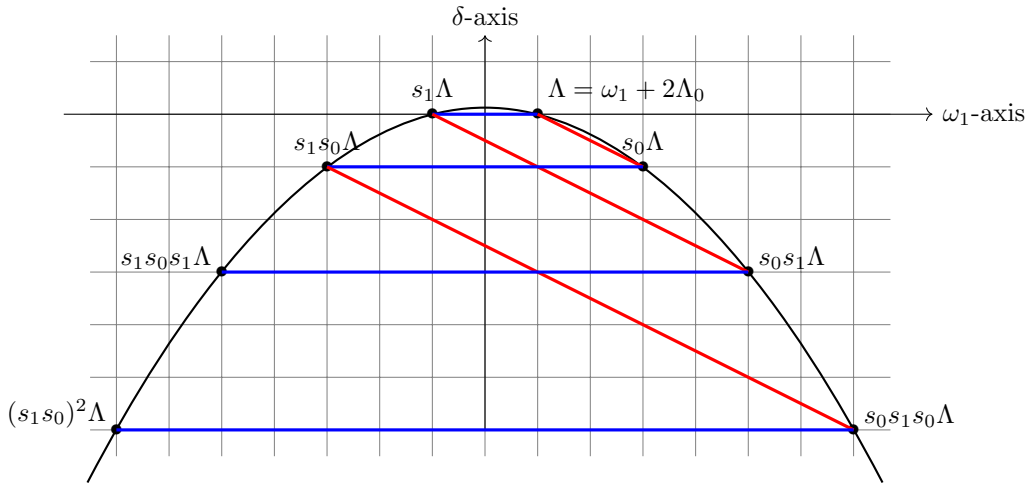
\hbox{The positive level orbit $W^{\mathrm{ad}}(\omega_1+2\Lambda_0)$ mod $\Lambda_0$}
\end{matrix}
$$
When the level (coefficient of $\Lambda_0$) is large the parabola is wide, and it gets tighter
as the level decreases.

\newpage 
$$
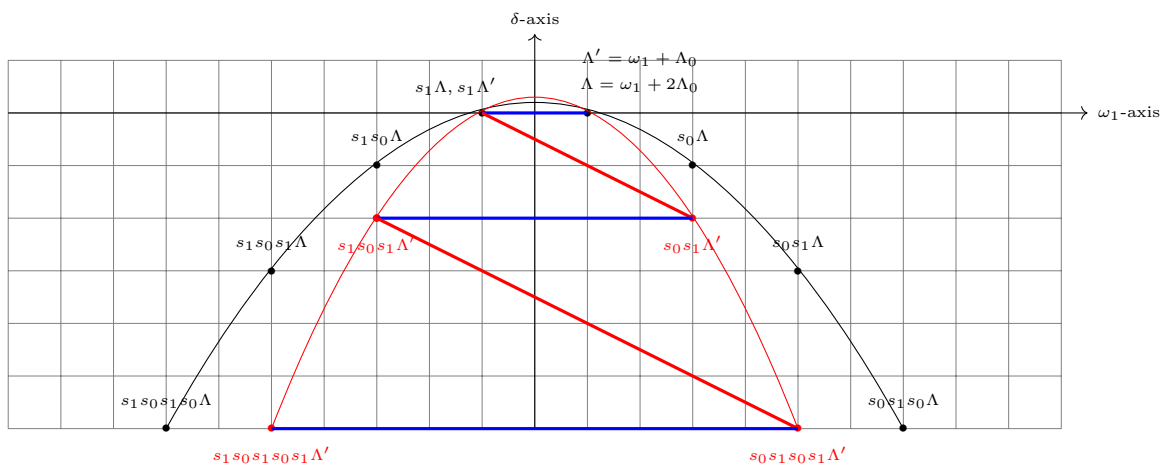
\begin{matrix}
\begin{tikzpicture}[scale=0.7]
\tikzstyle{every node}=[font=\tiny]    
    \draw[style=help lines,step=1cm] (-10,-6) grid (10,1);   
    \draw[->] (-10,0) -- (10.5,0) node[anchor=west]{$\omega_1$-axis};
    \draw[->] (0,-6) -- (0, 1.5) node[anchor=south]{$\delta$-axis};

     \node at (1,0) {$\bullet$};
        \filldraw (2,0.2) node[anchor=south] {$\Lambda=\omega_1+2\Lambda_0$};
        \filldraw (2,0.7) node[anchor=south] {$\Lambda'=\omega_1+\Lambda_0$};
            \node at (-1,0) {$\bullet$};
 \filldraw (-1.5,1) node[anchor=north,yshift=-0.1cm] {$s_1\Lambda, s_1\Lambda'$};
     \node at (3, -1) {$\bullet$};  
    \filldraw (3, 0) node[anchor=north,yshift=-0.1cm] {$s_0\Lambda$};
         \node at (-3, -1) {$\bullet$};  
   \filldraw (-3, 0) node[anchor=north,yshift=-0.1cm] {$s_1s_0\Lambda$};

         \node at (5, -3) {$\bullet$};  
   \filldraw (5, -2) node[anchor=north,yshift=-0.1cm] {$s_0s_1\Lambda$};
         \node at (-5, -3) {$\bullet$};  
   \filldraw (-5, -2) node[anchor=north,yshift=-0.1cm] {$s_1s_0s_1\Lambda$};

         \node at (7, -6) {$\bullet$};  
   \filldraw (7, -5) node[anchor=north,yshift=-0.1cm] {$s_0s_1s_0\Lambda$};
         \node at (-7, -6) {$\bullet$};  
   \filldraw (-7, -5) node[anchor=north,yshift=-0.1cm] {$s_1s_0s_1s_0\Lambda$};

		\draw (0, 0.2)parabola (-7, -6) ;
		\draw (0, 0.2)parabola (7, -6) ;
        \node[red] at (3, -2) {$\bullet$}; \node[red] at (-3, -2) {$\bullet$};
        \filldraw (3,-2) node[red, anchor=north,yshift=-0.1cm] {$s_0s_1\Lambda'$};
                \node[red] at (-3, -2) {$\bullet$};
        \filldraw (-3,-2) node[red, anchor=north,yshift=-0.1cm] {$s_1s_0s_1\Lambda'$};
        \node[red] at (5, -6) {$\bullet$};\node[red] at (-5, -6) {$\bullet$};
        \filldraw (5,-6) node[red, anchor=north,yshift=-0.1cm] {$s_0s_1s_0s_1\Lambda'$};
                \node[red] at (-3, -2) {$\bullet$};
        \filldraw (-5,-6) node[red, anchor=north,yshift=-0.1cm] {$s_1s_0s_1s_0s_1\Lambda'$};
		\draw [red] (0, 0.3)parabola (-5, -6) ;
		\draw [red] (0, 0.3)parabola (5, -6) ;
\draw [blue, very thick] (1,0) to (-1,0);
\draw [red, very thick] (-1,0) to (3,-2);
\draw [blue, very thick] (3,-2) to (-3,-2);
\draw [red, very thick] (-3,-2) to (5,-6);
\draw [blue, very thick] (5,-6) to (-5,-6);

\end{tikzpicture}
\\
\hbox{The parabolas bounding the orbits $W^{\mathrm{ad}}(\omega_1+2\Lambda_0)$ and 
$W^{\mathrm{ad}}(\omega_1+\Lambda_0)$ mod $\Lambda_0$}
\end{matrix}
$$
At level 0, the parabola pops and becomes two straight lines.
$$
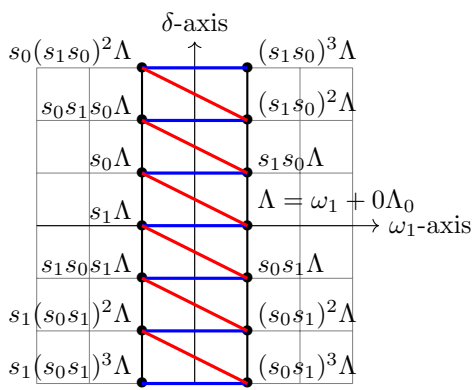
\begin{matrix}
\begin{tikzpicture}[scale=0.7]
\tikzstyle{every node}=[font=\small]
    
    \draw[style=help lines,step=1cm] (-3,-3) grid (3,3);    
    \draw[->] (-3,0) -- (3.5,0) node[anchor=west]{$\omega_1$-axis};
    \draw[->] (0,-3) -- (0, 3.5) node[anchor=south]{$\delta$-axis};
      \draw[-,thick] (-1,-3) -- (-1,3) ;
      \draw[-,thick] (1,-3) -- (1,3) ;
	\node at (1, -3) {$\bullet$};  
    \filldraw (1, -3) node[anchor=south west,yshift=-0.1cm] {$(s_0s_1)^3\Lambda$};
	\node at (1, -2) {$\bullet$};  
    \filldraw (1, -2) node[anchor=south west,yshift=-0.1cm] {$(s_0s_1)^2\Lambda$};
	\node at (1, -1) {$\bullet$};  
    \filldraw (1, -1) node[anchor=south west,yshift=-0.1cm] {$s_0s_1\Lambda$};
	\node at (1, 0) {$\bullet$};  
    \filldraw (1, 0) node[anchor=south west,yshift=0.05cm] {$\Lambda=\omega_1+0\Lambda_0$};
	\node at (1, 1) {$\bullet$};  
    \filldraw (1, 1) node[anchor=south west,yshift=-0.1cm] {$s_1s_0\Lambda$};
	\node at (1, 2) {$\bullet$};  
    \filldraw (1, 2) node[anchor=south west,yshift=-0.1cm] {$(s_1s_0)^2\Lambda$};
	\node at (1, 3) {$\bullet$};  
    \filldraw (1, 3) node[anchor=south west,yshift=-0.1cm] {$(s_1s_0)^3\Lambda$};
	\node at (-1, -3) {$\bullet$};  
    \filldraw (-1, -3) node[anchor=south east,yshift=-0.1cm] {$s_1(s_0s_1)^3\Lambda$};
	\node at (-1, -2) {$\bullet$};  
    \filldraw (-1, -2) node[anchor=south east,yshift=-0.1cm] {$s_1(s_0s_1)^2\Lambda$};
	\node at (-1, -1) {$\bullet$};  
    \filldraw (-1, -1) node[anchor=south east,yshift=-0.1cm] {$s_1s_0s_1\Lambda$};
	\node at (-1, 0) {$\bullet$};  
    \filldraw (-1, 0) node[anchor=south east,yshift=-0.1cm] {$s_1\Lambda$};
	\node at (-1, 1) {$\bullet$};  
    \filldraw (-1, 1) node[anchor=south east,yshift=-0.1cm] {$s_0\Lambda$};
	\node at (-1, 2) {$\bullet$};  
    \filldraw (-1, 2) node[anchor=south east,yshift=-0.1cm] {$s_0s_1s_0\Lambda$};
	\node at (-1, 3) {$\bullet$};  
    \filldraw (-1, 3) node[anchor=south east,yshift=-0.1cm] {$s_0(s_1s_0)^2\Lambda$};
	\foreach \y in {1,...,3}{
	\draw [blue, very thick] (-1,\y) to (1,\y);
	}	
	\foreach \y in {-3,...,0}{
	\draw [blue, very thick] (1,\y) to (-1,\y);
	}
	\foreach \x in {0,...,2}{
	\draw [red, very thick] (1,\x) to (-1,\x +1);
	}
	\foreach \x in {-3,...,-1}{
	\draw [red, very thick] (-1,\x+1) to (1,\x	);
	}
\end{tikzpicture}
\\
\hbox{The level 0 orbit $W^{\mathrm{ad}}\omega_1$ mod $\Lambda_0$}
\end{matrix}
$$
At negative level the parabola forms again, but this time facing the opposite
way, and getting wider as the level gets more and more negative.
$$
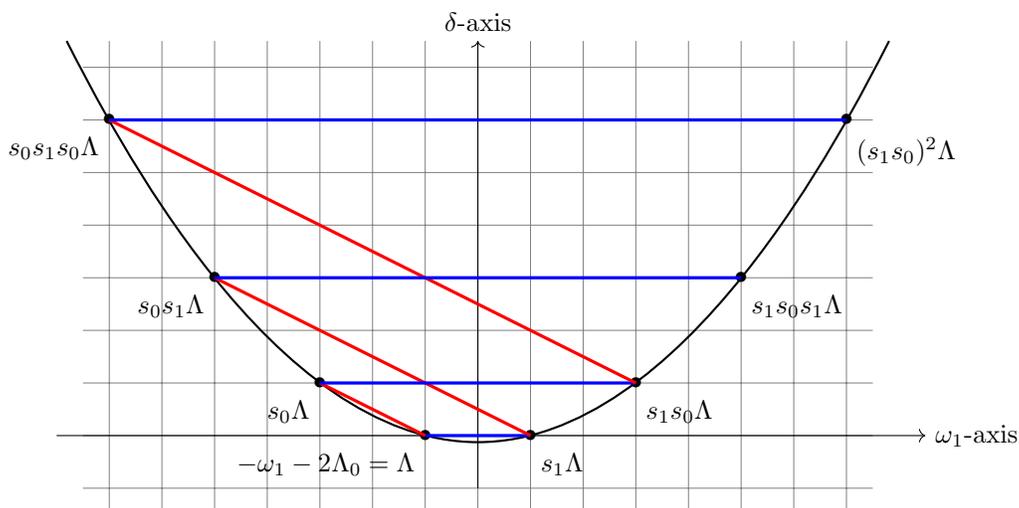
\begin{matrix}
\begin{tikzpicture}[scale=0.7]
\tikzstyle{every node}=[font=\small]
    
    \draw[style=help lines,step=1cm] (-7.5,7.5) grid (7.5,-1.5);
    
    \draw[->] (-8,0) -- (8.5,0) node[anchor=west]{$\omega_1$-axis};
    \draw[->] (0,-1) -- (0, 7.5) node[anchor=south]{$\delta$-axis};

     \node at (1,0) {$\bullet$};
        \filldraw (1,0) node[anchor=north west,yshift=-0.1cm] {$s_1\Lambda$};
            \node at (-1,0) {$\bullet$};
 \filldraw (-1,0) node[anchor=north east,yshift=-0.1cm] {$-\omega_1-2\Lambda_0=\Lambda$};
     \node at (3, 1) {$\bullet$};  
    \filldraw (3, 1) node[anchor=north west,yshift=-0.1cm] {$s_1s_0\Lambda$};
         \node at (-3, 1) {$\bullet$};  
   \filldraw (-3, 1) node[anchor=north east,yshift=-0.1cm] {$s_0\Lambda$};
         \node at (5, 3) {$\bullet$};  
   \filldraw (5, 3) node[anchor=north west,yshift=-0.1cm] {$s_1s_0s_1\Lambda$};
         \node at (-5, 3) {$\bullet$};  
   \filldraw (-5, 3) node[anchor=north east,yshift=-0.1cm] {$s_0s_1\Lambda$};
1
         \node at (7, 6) {$\bullet$};  
   \filldraw (7, 6) node[anchor=north west,yshift=-0.1cm] {$(s_1s_0)^2\Lambda$};
         \node at (-7, 6) {$\bullet$};  
   \filldraw (-7, 6) node[anchor=north east,yshift=-0.1cm] {$s_0s_1s_0\Lambda$};

		\draw[thick] (0, -0.125)parabola (-7.81024, 7.5) ;
		\draw[thick] (0, -0.125)parabola (7.81024, 7.5) ;
\draw [blue, very thick] (-1,0) to (1,0);
\draw [red, very thick] (1,0) to (-5,3);
\draw [red, very thick] (-1,0) to (-3,1);
\draw [blue, very thick] (-3,1) to (3,1);
\draw [red, very thick] (3,1) to (-7,6);
\draw [blue, very thick] (-5,3) to (5,3);
\draw [blue, very thick] (-7,6) to (7,6);
\end{tikzpicture}
\\
\hbox{The orbit $W^{\mathrm{ad}}(-\omega_1-2\Lambda_0)$ (negative level) mod $\Lambda_0$}
\end{matrix}
$$

The three different Bruhat orders on the affine Weyl group are visible on the $W^{\mathrm{ad}}$-orbits:
\begin{align*}
v \posleq\ w \qquad&\hbox{if $v(\omega_1+\Lambda_0)$ is higher than $w(\omega_1+\Lambda_0)$ in $W^{\mathrm{ad}}(\omega_1+\Lambda_0)$,} \\
v \zeroleq\ w \qquad&\hbox{if $v\omega_1$ is higher than $w(\omega_1+0\Lambda_0)$ in $W^{\mathrm{ad}}(\omega_1+0\Lambda_0)$.} \\
v \negleq\ w \qquad&\hbox{if $v(-\omega_1-\Lambda_0)$ is higher than 
$w(-\omega_1-\Lambda_0)$ in $W^{\mathrm{ad}}(-\omega_1-\Lambda_0)$.}
\end{align*}
The definitions of the Bruhat orders on $W^{\mathrm{ad}}$ and their relation to the 
closure order for Schubert cells in affine flag varieties is made precise in Section \ref{Bruhatorders}.
Indicative relations illustrating the from of the Hasse diagrams of the 
positive level, negative level, and level zero Bruhat orders for the Weyl group of
$\widehat{\fsl}_3$ are pictured in Plate A (there are additional relations which are not 
displayed in the pictures --- in an effort to make the periodicity pattern easily visible).

For the case of $\fg = \widehat{\fsl}_3$ the affine Weyl group orbits take
a similar form, with the points sitting on a downward paraboloid at positive
level, on an upward paraboloid at negative level and with the paraboloid popping
and becoming a tube at level 0 (for an example tube see the picture of $B(\omega_1+\omega_2)$ for 
$\widehat{\fsl}_3$ in Plate D).

$$
\begin{matrix}
\begin{tikzpicture}
\draw (0, 0) parabola (-3, -5) ;
\draw (0, 0)parabola (3, -5) ;
 \node at (-2,-1) {$\delta=0$};

 \filldraw (1.3, -1) circle (0.04cm) node{}; 
  \filldraw (-1.3, -1) circle (0.04cm) node{}; 
    \filldraw (0.7, -0.7) circle (0.04cm) node{}; 
        \filldraw (-0.3, -0.7) circle (0.04cm) node{}; 
            \filldraw (-0.7, -1.3) circle (0.04cm) node{}; 
        \filldraw (0.3, -1.3) circle (0.04cm) node{}; 
\tikzstyle{conefill} = [fill=blue!20,fill opacity=0.8]
 \filldraw[conefill](1.3, -1)--(0.7, -0.7)--(-0.3, -0.7)--(-1.3, -1)--(-0.7, -1.3)--(0.3, -1.3)--(1.3, -1)--cycle;  

\filldraw (1.4+0.4, -3.3-0.2) circle (0.04cm) node{};   
\filldraw (1.4+0.4*2, -3.3-0.2*2) circle (0.04cm) node{};   
\filldraw (-0.6+2/3, -3.3) circle (0.04cm) node{};   
\filldraw (-0.6+4/3, -3.3) circle (0.04cm) node{};   
\filldraw (-0.6-2/3, -3.3-0.2) circle (0.04cm) node{};   
\filldraw (-0.6-4/3, -3.3-0.4) circle (0.04cm) node{};   
\filldraw (-2.6+0.4, -3.9-0.2) circle (0.04cm) node{};   
\filldraw (-2.6+0.8, -3.9-0.4) circle (0.04cm) node{};   
\filldraw (-1.4+2/3, -4.5) circle (0.04cm) node{};   
\filldraw (-1.4+4/3, -4.5) circle (0.04cm) node{};   
\filldraw (0.6+2/3, -4.5+0.2) circle (0.04cm) node{};   
\filldraw (0.6+4/3, -4.5+0.4) circle (0.04cm) node{};   
  \node at (-3.5,-4) {$\delta=-3$};
   \tikzstyle{conefill} = [fill=blue!20,fill opacity=0.8]
  \filldraw[conefill](2.6, -3.9)--(1.4, -3.3)--(-0.6, -3.3)--(-2.6, -3.9)--(-1.4, -4.5)--(0.6, -4.5)--(2.6, -3.9) ;
\end{tikzpicture}
&
\begin{tikzpicture}
\draw (0, 0) parabola (-3, 5) ;
\draw (0, 0)parabola (3, 5) ;
 \filldraw (1.3, 1) circle (0.04cm) node{}; 
  \filldraw (-1.3, 1) circle (0.04cm) node{}; 
    \filldraw (0.7, 0.7) circle (0.04cm) node{}; 
        \filldraw (-0.3, 0.7) circle (0.04cm) node{}; 
            \filldraw (-0.7, 1.3) circle (0.04cm) node{}; 
        \filldraw (0.3, 1.3) circle (0.04cm) node{}; 
                      \tikzstyle{conefill} = [fill=blue!20,fill opacity=0.8]
  \filldraw[conefill](1.3, 1)--(0.7, 0.7)--(-0.3, 0.7)--(-1.3, 1)--(-0.7, 1.3)--(0.3, 1.3)--(1.3, 1)--cycle;  

\filldraw (1.4+0.4, 3.3+0.2) circle (0.04cm) node{};   
\filldraw (1.4+0.4*2, 3.3+0.2*2) circle (0.04cm) node{};   
\filldraw (-0.6+2/3, 3.3) circle (0.04cm) node{};   
\filldraw (-0.6+4/3, 3.3) circle (0.04cm) node{};   
\filldraw (-0.6-2/3, 3.3+0.2) circle (0.04cm) node{};   
\filldraw (-0.6-4/3, 3.3+0.4) circle (0.04cm) node{};   
\filldraw (-2.6+0.4, 3.9+0.2) circle (0.04cm) node{};   
\filldraw (-2.6+0.8, 3.9+0.4) circle (0.04cm) node{};   
\filldraw (-1.4+2/3, 4.5) circle (0.04cm) node{};   
\filldraw (-1.4+4/3, 4.5) circle (0.04cm) node{};   
\filldraw (0.6+2/3, 4.5-0.2) circle (0.04cm) node{};   
\filldraw (0.6+4/3, 4.5-0.4) circle (0.04cm) node{};   
  \node at (-3.2,4) {$\delta=3$};
 \node at (-2,1) {$\delta=0$};
   \tikzstyle{conefill} = [fill=blue!20,fill opacity=0.8]
  \filldraw[conefill](2.6, 3.9)--(1.4, 3.3)--(-0.6, 3.3)--(-2.6, 3.9)--(-1.4, 4.5)--(0.6, 4.5)--(2.6, 3.9) ;
\end{tikzpicture}
\\
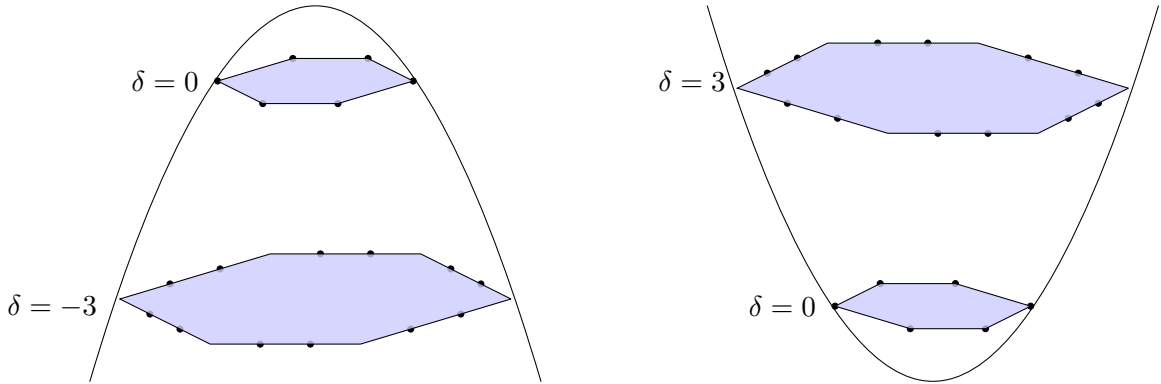
\hbox{Positive level orbit $W^{\mathrm{ad}}(\omega_1+\omega_2+2\Lambda_0)$ for $\widehat{\fsl}_3$}
&\hbox{Negative level orbit $W^{\mathrm{ad}}(-\omega_1-\omega_2-2\Lambda_0)$ for $\widehat{\fsl}_3$}
\end{matrix}
$$

\newcommand{\xslant}{-0.8} 
\newcommand{\yslant}{0.6}

\subsection{Extremal weight modules $L(\Lambda)$ and their crystals $B(\Lambda)$}

For the affine Lie algebra $\fg = \widehat{\fsl}_2$
the weights of integrable $\fg$-modules always lie in the set
$$\fh^*_\ZZ = \CC\delta + \ZZ\omega_1 + \ZZ\Lambda_0.$$
A set of representatives for the orbits of the action of $W^{\mathrm{ad}}$ on $\fh_{\ZZ}^{\ast}$ is 
$$
(\fh^{\ast})_{\mathrm{int}}=(\fh^{\ast})_{\mathrm{int}}^+\cup (\fh^{\ast})_{\mathrm{int}}^0\cup (\fh^{\ast})_{\mathrm{int}}^-,
\quad\hbox{where}\quad
\begin{array}{rl}
(\fh^{\ast})_{\mathrm{int}}^0 &= \{ a\delta + n\omega_1\in \fh_\ZZ^*\ |\  0\le n\},  \\
(\fh^{\ast})_{\mathrm{int}}^+ &= \{a\delta + m\omega_1+n\Lambda_0\in \fh_\ZZ^*\ |\ 
0\le m\le n\}, \\
(\fh^{\ast})_{\mathrm{int}}^- &= \{ a\delta  -m\omega_1 - n\Lambda_0\in \fh_\ZZ^*\ |\  
0\le m\le n\}.
\end{array}
$$
These sets are illustrated (mod $\delta$) below.

For each of the $\Lambda\in (\fh^*)_\mathrm{int}$, there is a (universal) integrable
\emph{extremal weight} module $L(\Lambda)$, which is highest weight if
$\Lambda\in (\fh^*)_{\mathrm{int}}^+$, is lowest weight (and not highest weight)
if $\Lambda\in (\fh^*)_{\mathrm{int}}^-$, and which is neither highest or lowest
weight when $\Lambda\in (\fh^*)_{\mathrm{int}}^0$.  The module $L(\Lambda)$
has a crystal $B(\Lambda)$.  

At positive level and negative levels the crystals $B(\Lambda)$ are connected,
but the crystal $B(\Lambda)$ is usually not connected in level $0$.  The connected
components and their structure are known explicitly from a combination of 
results of Kashiwara, Beck-Chari-Pressley, Nakajima, Beck-Nakajima, Fourier-Littelmann,
Ion and others.  These results are collected in Theorem \ref{levelzerostructure}
and equation \eqref{levzerochar} expresses the characters 
of the $L(\Lambda)$ in terms of Macdonald polynomials specialised at $t=0$.
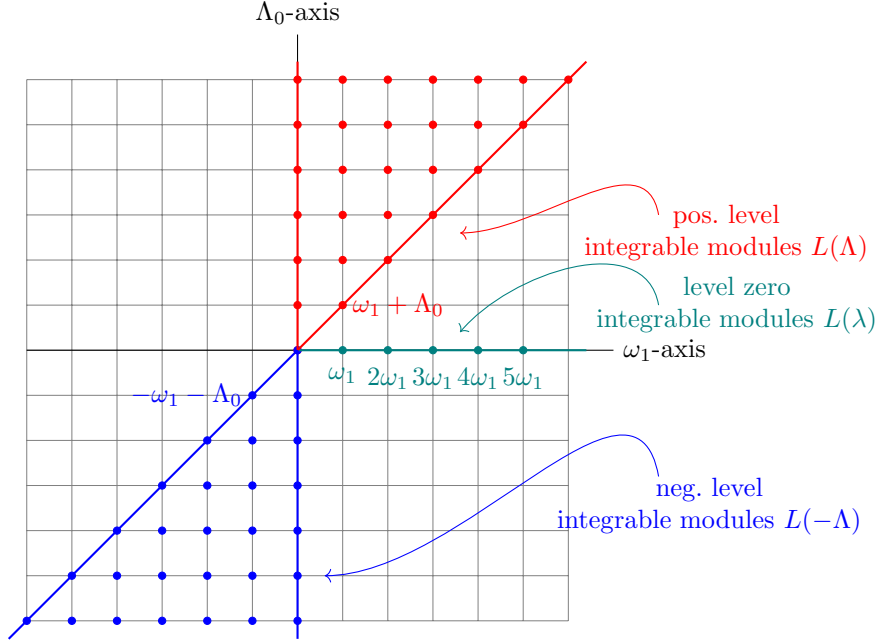
\begin{equation*}
\begin{tikzpicture}[scale=1.2]
\tikzstyle{every node}=[font=\small]
    \draw[style=help lines,step=0.5cm] (-3,-3) grid (3,3);
    \draw[-] (-3,0) -- (3.5,0) node[anchor=west]{$\omega_1$-axis};
    \draw[-] (0,-3) -- (0,3.5) node[anchor=south]{$\Lambda_0$-axis};
\draw [<-,red, in=100, out=0,  thin] (1.8,1.3) to node[anchor=west, xshift=-0.3cm, yshift=-0.5cm]{$\begin{array}{c}
\hbox{pos.\ level} \\ \hbox{integrable modules $L(\Lambda)$}\end{array}$} (4,1.5);
\draw [->,teal, in=50, out=100,  thin] (4,0.5) to node[anchor=west, yshift=-0.5cm]{$\begin{array}{c}
\hbox{level zero} \\ \hbox{integrable modules $L(\lambda)$}\end{array}$} (1.8,0.2);
\draw [<-,blue, in=100, out=0,  thin] (0.3,-2.5) to node[anchor=west, yshift=-0.4cm]{$\begin{array}{c}
\hbox{neg.\ level} \\ \hbox{integrable modules $L(-\Lambda)$}\end{array}$} (4,-1.4);

    \begin{scope}[color=teal]
    \filldraw (0.5,0) circle (0.04cm) node (A) {} node[anchor=north, yshift=-0.1cm] {$\omega_1$};
     \filldraw (1,0) circle (0.04cm) node (A) {} node[anchor=north,yshift=-0.1cm] {$2\omega_1$};
     \filldraw (1.5,0) circle (0.04cm) node (A) {} node[anchor=north,yshift=-0.1cm] {$3\omega_1$};
 \filldraw (2,0) circle (0.04cm) node (A) {} node[anchor=north,yshift=-0.1cm] {$4\omega_1$};
  \filldraw (2.5,0) circle (0.04cm) node (A) {} node[anchor=north,yshift=-0.1cm] {$5\omega_1$};
    \end{scope}
    
    \begin{scope}[color=red]
        \filldraw (0, 0) circle (0.04cm) node (B) {} node[anchor=north,yshift=-0.1cm] {};
    \filldraw (0.5,0.5) circle (0.04cm) node (B) {} node[anchor=west,yshift=0.0cm] {$\omega_1+\Lambda_0$};
     \filldraw (1,1) circle (0.04cm) node (B) {} node[anchor=north,yshift=-0.1cm] {};
     \filldraw (1.5,1.5) circle (0.04cm) node (B) {} node[anchor=north,yshift=-0.1cm] {};
 \filldraw (2,2) circle (0.04cm) node (B) {} node[anchor=north, yshift=-0.1cm] {};
  \filldraw (2.5,2.5) circle (0.04cm) node (B) {} node[anchor=north, yshift=-0.1cm] {};
   \filldraw (3,3) circle (0.04cm) node (B) {} node[anchor=north, yshift=-0.1cm] {}; 
 \filldraw (0, 0.5) circle (0.04cm) node (B) {} node[anchor=north,yshift=-0.1cm] {};
    \filldraw (0.5,1) circle (0.04cm) node (B) {} node[anchor=north,yshift=-0.1cm] {};
     \filldraw (1,1.5) circle (0.04cm) node (B) {} node[anchor=north,yshift=-0.1cm] {};
     \filldraw (1.5,2) circle (0.04cm) node (B) {} node[anchor=north,yshift=-0.1cm] {};
 \filldraw (2,2.5) circle (0.04cm) node (B) {} node[anchor=north, yshift=-0.1cm] {};
  \filldraw (2.5,3) circle (0.04cm) node (B) {} node[anchor=north, yshift=-0.1cm] {};
   \filldraw (0, 1) circle (0.04cm) node (B) {} node[anchor=north,yshift=-0.1cm] {};
    \filldraw (0.5,1.5) circle (0.04cm) node (B) {} node[anchor=north,yshift=-0.1cm] {};
     \filldraw (1,2) circle (0.04cm) node (B) {} node[anchor=north,yshift=-0.1cm] {};
     \filldraw (1.5,2.5) circle (0.04cm) node (B) {} node[anchor=north,yshift=-0.1cm] {};
 \filldraw (2,3) circle (0.04cm) node (B) {} node[anchor=north, yshift=-0.1cm] {};
   \filldraw (0, 1.5) circle (0.04cm) node (B) {} node[anchor=north,yshift=-0.1cm] {};
    \filldraw (0.5,2) circle (0.04cm) node (B) {} node[anchor=north,yshift=-0.1cm] {};
     \filldraw (1,2.5) circle (0.04cm) node (B) {} node[anchor=north,yshift=-0.1cm] {};
     \filldraw (1.5,3) circle (0.04cm) node (B) {} node[anchor=north,yshift=-0.1cm] {};
   \filldraw (0, 2) circle (0.04cm) node (B) {} node[anchor=north,yshift=-0.1cm] {};
    \filldraw (0.5,2.5) circle (0.04cm) node (B) {} node[anchor=north,yshift=-0.1cm] {};
     \filldraw (1,3) circle (0.04cm) node (B) {} node[anchor=north,yshift=-0.1cm] {};
   \filldraw (0, 2.5) circle (0.04cm) node (B) {} node[anchor=north,yshift=-0.1cm] {};
    \filldraw (0.5,3) circle (0.04cm) node (B) {} node[anchor=north,yshift=-0.1cm] {};
   \filldraw (0, 3) circle (0.04cm) node (B) {} node[anchor=north,yshift=-0.1cm] {};
    \end{scope}
        \begin{scope}[color=blue]
        \filldraw (0, 0) circle (0.04cm) node (B) {} node[anchor=north,yshift=-0.1cm] {};
    \filldraw (-0.5,-0.5) circle (0.04cm) node (B) {} node[anchor=east,yshift=-0.0cm ] {$-\omega_1-\Lambda_0$};
     \filldraw (-1,-1) circle (0.04cm) node (B) {} node[anchor=north,yshift=-0.1cm] {};
     \filldraw (-1.5,-1.5) circle (0.04cm) node (B) {} node[anchor=north,yshift=-0.1cm] {};
 \filldraw (-2,-2) circle (0.04cm) node (B) {} node[anchor=north, yshift=-0.1cm] {};
  \filldraw (-2.5,-2.5) circle (0.04cm) node (B) {} node[anchor=north, yshift=-0.1cm] {};
   \filldraw (-3,-3) circle (0.04cm) node (B) {} node[anchor=north, yshift=-0.1cm] {}; 
 \filldraw (0, -0.5) circle (0.04cm) node (B) {} node[anchor=north,yshift=-0.1cm] {};
    \filldraw (-0.5,-1) circle (0.04cm) node (B) {} node[anchor=north,yshift=-0.1cm] {};
     \filldraw (-1,-1.5) circle (0.04cm) node (B) {} node[anchor=north,yshift=-0.1cm] {};
     \filldraw (-1.5,-2) circle (0.04cm) node (B) {} node[anchor=north,yshift=-0.1cm] {};
 \filldraw (-2,-2.5) circle (0.04cm) node (B) {} node[anchor=north, yshift=-0.1cm] {};
  \filldraw (-2.5,-3) circle (0.04cm) node (B) {} node[anchor=north, yshift=-0.1cm] {};
   \filldraw (0, -1) circle (0.04cm) node (B) {} node[anchor=north,yshift=-0.1cm] {};
    \filldraw (-0.5,-1.5) circle (0.04cm) node (B) {} node[anchor=north,yshift=-0.1cm] {};
     \filldraw (-1,-2) circle (0.04cm) node (B) {} node[anchor=north,yshift=-0.1cm] {};
     \filldraw (-1.5,-2.5) circle (0.04cm) node (B) {} node[anchor=north,yshift=-0.1cm] {};
 \filldraw (-2,-3) circle (0.04cm) node (B) {} node[anchor=north, yshift=-0.1cm] {};
   \filldraw (0, -1.5) circle (0.04cm) node (B) {} node[anchor=north,yshift=-0.1cm] {};
    \filldraw (-0.5,-2) circle (0.04cm) node (B) {} node[anchor=north,yshift=-0.1cm] {};
     \filldraw (-1,-2.5) circle (0.04cm) node (B) {} node[anchor=north,yshift=-0.1cm] {};
     \filldraw (-1.5,-3) circle (0.04cm) node (B) {} node[anchor=north,yshift=-0.1cm] {};
   \filldraw (0, -2) circle (0.04cm) node (B) {} node[anchor=north,yshift=-0.1cm] {};
    \filldraw (-0.5,-2.5) circle (0.04cm) node (B) {} node[anchor=north,yshift=-0.1cm] {};
     \filldraw (-1,-3) circle (0.04cm) node (B) {} node[anchor=north,yshift=-0.1cm] {};
   \filldraw (0, -2.5) circle (0.04cm) node (B) {} node[anchor=north,yshift=-0.1cm] {};
    \filldraw (-0.5,-3) circle (0.04cm) node (B) {} node[anchor=north,yshift=-0.1cm] {};
   \filldraw (0, -3) circle (0.04cm) node (B) {} node[anchor=north,yshift=-0.1cm] {};
    \end{scope}
      \draw [teal, thick] (0,0) -- (3.2, 0);
         \draw [red, thick] (0,0) -- (3.2, 3.2);
            \draw [red, thick] (0,0) -- (0, 3.2);
            \draw [blue, thick] (0,0) -- (-3.2, -3.2);
            \draw [blue, thick] (0,0) -- (0, -3.2);
    \end{tikzpicture}
\label{integrablewtspicture}
\end{equation*}

\subsection{Affine flag varieties $G/I^+$, $G/I^0$ and $G/I^-$}

There are three kinds of affine flag varieties for the loop group
$G=G(\CC((\epsilon)))$: the positive level (thin) affine flag variety $G/I^+$, the negative level
(thick) affine flag variety $G/I^-$ and the level 0 (semi-infinite) affine flag variety $G/I^0$.
A combination of results of Kumar, Mathieu, Kashiwara, Kashiwara-Tanisaki,
Kashiwara-Shimozono, Varagnolo-Vasserot, Lusztig and Braverman-Finkelberg have made it
clear that there is a Borel-Weil-Bott theorem for each of these:
\begin{align*}
H^0(G/I^+, \cL_\Lambda) \cong L(\Lambda), \quad
&\hbox{for positive level $\Lambda\in (\fh^*)^+_{\mathrm{int}}$,} \\
H^0(G/I^0, \cL_\lambda) \cong L(\lambda), \quad
&\hbox{for level zero $\lambda\in (\fh^*)^0_{\mathrm{int}}$,} \\
H^0(G/I^-, \cL_{-\Lambda}) \cong L(-\Lambda), \quad
&\hbox{for negative $-\Lambda\in (\fh^*)^-_{\mathrm{int}}$.}
\end{align*}
These Borel-Weil-Bott theorems tightly connect the representation theory with the geometry.
In all essential aspects the combinatorics of the positive level affine flag
variety and the loop Grassmannian generalizes to the negative level and the
level 0 affine flag varieties.  

Section \ref{affineflags}
extends the results of \cite{PRS}
and displays the alcove walk combinatorics for each of the three cases (positive level,
negative level and level 0) in parallel.  In addition
it describes the method for deriving the natural affine Hecke algebra actions on the function spaces
$C(G/I^+)$, $C(I^+\backslash G/I^+)$, $C(I^+\backslash G/I^+)$ and $C(I^-\backslash G/I^+)$.

\subsection{References, technicalities and acknowledgements}

Section \ref{affLiealgebra} introduces the affine Lie algebra and the homogeneous Heisenberg subalgebra
following \cite{Kac} and Section \ref{affineWeyl} gives explicit matrices describing the actions 
of the affine Weyl group $W^{\mathrm{ad}}$ on the affine Cartan $\fh$ and its dual $\fh^*$.
Section \ref{Bruhatorders} defines the Bruhat orders on the affine Weyl group and explains their relation
to the corresponding affine flag varieties following
\cite{Kum} and \cite[\S 7 and 11]{LuICM}.  
Sections \ref{bdgps} and \ref{Macpolys} introduce the affine braid groups 
and Macdonald polynomials following \cite{RY11}.  Section \ref{Macspec} treats the specializations
of (nonsymmetric) Macdonald polynomials at $t=0$, $t=\infty$, $q=0$ and $q=\infty$ and
reviews the result of Ion \cite{Ion01} that relates Macdonald polynomials at $t=0$ to Demazure
operators.

Section 2 follows \cite{BN02}, \cite{B94}, and \cite{BCP98}, introducing the quantum affine algebra
$\mathbf{U}$, the conversion to its loop presentation, 
the PBW-type elements and the quantum homogeneous Heisenberg subalgebra.  
Section \ref{intmodules} defines integrable $\mathbf{U}$-modules
and Section \ref{Extwtmodules} introduces the extremal weight modules $L(\Lambda)$ 
following \cite[(8.2.2)]{Kas94} and \cite[\S3.1]{Kas02}.
Section \ref{Demazurechars} reviews the Demazure character formulas for extremal weight modules.
Section \ref{levelzeroloopsection} discusses the loop presentation of the level 0
extremal weight modules and the fact that these coincide with the universal standard modules of
\cite{Nak99} and the global Weyl modules of \cite{CP01}.  Letting $\mathbf{U}'$ be $\mathbf{U}$ without
the element $D$, Section \ref{findimmodules} explains how to
shrink the extremal weight module to a local Weyl module and how this provides a classification
of finite dimensional simple $\mathbf{U'}$-modules by Drinfeld polynomials.

We have made a concerted effort to make a useful survey.  In
order to simplify the exposition we have brushed under the rug a number 
of technicalities which are wisely ignored when one learns the subject (particularly 
(a) the difference between simply laced cases and the general case requires proper 
attention to the  diagonal matrix which symmetrizes the affine Cartan matrix \cite[(2.1.1)]{Kac}
causing the constants $d_i$ which pepper the quantum group literature
and (b) the machinations necessary for allowing multiple parameters $t_i^{\frac12}$ 
in Macdonald polynomials).
The reader who needs to sort out these features is advised to drink a strong double
espresso to optimise clear thinking, 
consult the references (particularly \cite{BN02} and \cite{RY11}) 
and \emph{not trust our exposition}.
Perhaps in the future a more complete (probably book length) version of this paper
will be completed which allows us to attend more carefully to these nuances
and include more detailed proofs.  Having made this point, we can say
that a careful effort has been made to provide specific references to the literature at every step and
we hope that this will be useful for the reader that wishes to go further.


We thank Martha Yip for conversations, calculations and teaching us
so much about Macdonald polynomials over the years.
We thank all the institutions which have supported our work on this paper, 
particularly the University of Melbourne and the 
Australian Research Council (grants DP1201001942, DP130100674 and DE190101231). 
Arun Ram thanks IHP (Institut Henri Poincar\'e) 
and  ICERM (Institute for Computational and Experimental Research in Mathematics)
for support, hospitality and a stimulating working environments at the thematic programs 
``Automorphic Forms, Combinatorial Representation Theory and Multiple Dirichlet Series'' and
``Combinatorics and Interactions".

Finally, it is a pleasure to dedicate this paper to Ian Macdonald and Alun Morris who forcefully
led the way to the kinds explorations of affine combinatorial representation theory that are happening 
these days.

\newpage
\noindent
\textbf{PLATE A: Bruhat orders on the affine Weyl group (partial, indicative, relations)}
$$
\begin{matrix}
\resizebox{1.2cm}{!}{
\begin{tikzpicture}[scale=1.2]
\tikzstyle{every node}=[font=\large]
	 \coordinate (omega1) at (0,1);	 
	\filldraw (0,0) node {$1$};    
	\filldraw (0,1) node {$s_0$};  
	\filldraw (0,2) node {$s_0s_1$};  
	\filldraw (0,3) node {$s_0s_1s_0$};  
	\filldraw (0,-3) node {$s_1s_0s_1$};  
	\filldraw (0,-2) node {$s_1s_0$};  
	\filldraw (0,-1) node {$s_1$};  

	\draw[-latex, black, thick]  (0,2.8) to (0,2.2);
	\draw[-latex, black, thick]  (0,1.8) to (0,1.2);
	\draw[-latex, black, thick]  (0,0.8) to (0,0.2);
	\draw[-latex, black, thick]  (0,-0.8) to (0,-0.2);
	\draw[-latex, black, thick]  (0,-1.8) to (0,-1.2);
	\draw[-latex, black, thick]  (0,-2.8) to (0,-2.2);
\end{tikzpicture}
}
&\resizebox{7.2cm}{!}{
\begin{tikzpicture}[scale=1.2]
\tikzstyle{every node}=[font=\large]
	 \coordinate (omega2) at (-0.5,0.86603);
	 \coordinate (2omega2) at  (-1,2*0.86603);
	 \coordinate (2omega2+omega1) at (-0.5, 3*0.86603);
	 \coordinate (omega2+2omega1) at (0.5,3*0.86603);
	 \coordinate (2omega1) at (1,2*0.86603);
	 \coordinate (omega1) at (0.5,0.86603);
	 \coordinate (-omega1) at (-0.5,-0.86603);
	 \coordinate (-omega1+omega2) at (-1,0);
	 \coordinate (-omega2+omega1) at (1,0);
	 \coordinate (-omega2) at  (0.5,-0.86603);
	 \coordinate (-2omega1-omega2) at (-0.5,-3*0.86603);
	 \coordinate (-2omega1) at (-1, -2*0.86603);
	 \coordinate (-2omega2) at (1,-2*0.86603);
	 \coordinate (-2omega2-omega1) at (0.5,-3*0.86603);
	 \coordinate (rho) at ($(omega1)+(omega2)$);
	 \coordinate (-rho) at ($-1*(omega1)-(omega2)$);
	 \coordinate (alpha1) at ($2*(omega1)-(omega2)$);
	 \coordinate (-alpha1) at ($-2*(omega1)+(omega2)$);
	 \coordinate (alpha2) at ($2*(omega2)-(omega1)$);
	 \coordinate (-alpha2) at ($-2*(omega2)+(omega1)$);	 
	 
	\filldraw (rho) node {$1$};    
	\filldraw (alpha1) node {$s_2$};  
	\filldraw (alpha2) node {$s_1$};  
	\filldraw (-alpha1) node {$s_1s_2$};  
	\filldraw (-alpha2) node {$s_2s_1$};  
	\filldraw (-rho) node {$s_1s_2s_1$}; 
	\filldraw ($2*(rho)+(alpha1)$) node {$s_0s_2$};    
	\filldraw ($2*(rho)$) node {$s_0$};  
	\filldraw ($2*(alpha1)$) node {$s_2s_0$};  
	\filldraw ($2*(alpha1)+(rho)$) node {$s_2s_0s_2$};  
	\filldraw ($2*(rho)+(alpha2)$) node {$s_0s_1$};    
	\filldraw ($2*(alpha2)$) node {$s_1s_0$};  
	\filldraw ($2*(alpha2)+(rho)$) node {$s_1s_0s_1$};  
	\filldraw ($2*(rho)+(alpha2)$) node {$s_0s_1$};    
	\filldraw ($-2*(alpha1)$) node {$s_1s_2s_0$};  
	\filldraw ($-2*(alpha1)+(alpha2)$) node {$s_1s_2s_0s_2$};  
	\filldraw ($-1*(alpha1)+2*(alpha2)$) node {$s_1s_0s_2$};  
	\filldraw ($2*(rho)+(alpha1)$) node {$s_0s_2$};    
	\filldraw ($-2*(alpha2)$) node {$s_2s_1s_0$};  
	\filldraw ($-2*(alpha2)+(alpha1)$) node {$s_2s_1s_0s_1$};  
	\filldraw ($-1*(alpha2)+2*(alpha1)$) node {$s_2s_0s_1$};  
	\filldraw ($-3*(rho)+(alpha1)$) node {$s_1s_2s_1s_0s_2$};    
	\filldraw ($-2*(rho)$) node {$s_1s_2s_1s_0$};  
	\filldraw ($2*(alpha1)-3*(rho)$) node {$s_2s_1s_0s_2$};  
	\filldraw ($-3*(rho)+(alpha2)$) node {$s_2s_1s_2s_0s_1$};    
	\filldraw ($2*(alpha2)-3*(rho)$) node {$s_1s_2s_0s_1$};  	

	\draw[-latex, black, thick]  ($(-alpha1)+0.15*(alpha1)$) to ($(alpha1)-0.15*(alpha1)$); 			\draw[-latex, black, thick]  ($2*(alpha1)+(rho)-0.15*(alpha1)$) to ($(rho)+0.15*(alpha1)$); 		
	\draw[-latex, black, thick]  ($2*(alpha2)+0.15*(alpha1)$) to ($2*(rho)-0.15*(alpha1)$); 		
	\draw[-latex, black, thick]  ($-2*(alpha1)+(alpha2)+0.15*(alpha1)$) to ($(alpha2)-0.15*(alpha1)$); 
	\draw[-latex, black, thick]  ($(-alpha2)+0.15*(alpha1)$) to ($-1*(alpha2)+2*(alpha1)-0.15*(alpha1)$); 
	\draw[-latex, black, thick]  ($2*(alpha2)-3*(rho)+0.15*(alpha1)$) to ($(-rho)-0.15*(alpha1)$); 
	\draw[-latex, black, thick]  ($-2*(rho)+0.15*(alpha1)$) to ($-2*(alpha2)-0.15*(alpha1)$); 
	\draw[-latex, black, thick]  ($(-alpha2)+0.15*(alpha2)$) to ($(alpha2)-0.15*(alpha2)$); 			\draw[-latex, black, thick]  ($2*(alpha2)+(rho)-0.15*(alpha2)$) to ($(rho)+0.15*(alpha2)$);
	\draw[-latex, black, thick]  ($2*(alpha1)+0.15*(alpha2)$) to ($2*(rho)-0.15*(alpha2)$); 		
	\draw[-latex, black, thick]  ($(-alpha1)+0.15*(alpha2)$) to ($-1*(alpha1)+2*(alpha2)-0.15*(alpha2)$); 		
	\draw[-latex, black, thick]  ($-2*(alpha2)+(alpha1)+0.15*(alpha2)$) to ($(alpha1)-0.15*(alpha2)$); 
	\draw[-latex, black, thick]  ($-2*(rho)+0.15*(alpha2)$) to ($-2*(alpha1)-0.15*(alpha2)$); 
	\draw[-latex, black, thick]  ($2*(alpha1)-3*(rho)+0.15*(alpha2)$) to ($(-rho)-0.15*(alpha2)$); 
		\draw[-latex, black, thick]  ($(alpha1)+0.15*(alpha2)$) to ($(rho)-0.15*(alpha2)$); 			
	\draw[-latex, black, thick]  ($(alpha2)+0.15*(alpha1)$) to ($(rho)-0.15*(alpha1)$); 			
	\draw[-latex, black, thick]  ($1.85*(rho)$) to ($1.15*(rho)$); 											
	\draw[-latex, black, thick]  ($(-alpha1)+0.15*(rho)$) to ($(alpha2)-0.15*(rho)$); 				
	\draw[-latex, black, thick]  ($(-alpha2)+0.15*(rho)$) to ($(alpha1)-0.15*(rho)$); 				
	\draw[-latex, black, thick]  ($-1*(rho)+0.15*(alpha2)$) to ($(-alpha1)-0.15*(alpha2)$); 	
	\draw[-latex, black, thick]  ($-1*(rho)+0.15*(alpha1)$) to ($(-alpha2)-0.15*(alpha1)$); 	
	\draw[-latex, black, thick]  ($-1.85*(rho)$) to ($-1.15*(rho)$); 										
	\draw[-latex, black, thick]  ($1.85*(alpha1)$) to ($1.15*(alpha1)$);									
	\draw[-latex, black, thick]  ($1.85*(alpha2)$) to ($1.15*(alpha2)$);									
	\draw[-latex, black, thick]  ($-1.85*(alpha1)$) to ($-1.15*(alpha1)$);								
	\draw[-latex, black, thick]  ($-1.85*(alpha2)$) to ($-1.15*(alpha2)$);								
	\draw[-latex, black, thick]  ($2*(alpha2)+0.85*(rho)$) to ($2*(alpha2)+0.15*(rho)$);		
	\draw[-latex, black, thick]  ($2*(alpha1)+0.85*(rho)$) to ($2*(alpha1)+0.15*(rho)$);		
	\draw[-latex, black, thick]  ($2*(-alpha1)-0.85*(rho)$) to ($2*(-alpha1)-0.15*(rho)$);		
	\draw[-latex, black, thick]  ($2*(-alpha2)-0.85*(rho)$) to ($2*(-alpha2)-0.15*(rho)$);		
	\draw[-latex, black, thick]  ($2*(-alpha1)+0.85*(alpha2)$) to ($2*(-alpha1)+0.15*(alpha2)$);	
	\draw[-latex, black, thick]  ($2*(-alpha2)+0.85*(alpha1)$) to ($2*(-alpha2)+0.15*(alpha1)$);	
	\draw[-latex, black, thick] ($2*(alpha2)-0.85*(alpha1)$) to ($2*(alpha2)-0.15*(alpha1)$); 		
	\draw[-latex, black, thick] ($2*(alpha1)-0.85*(alpha2)$) to ($2*(alpha1)-0.15*(alpha2)$); 		
	\draw[-latex, black, thick] ($-2*(alpha1)+(alpha2)+0.15*(rho)$) to ($-2*(alpha1)+(alpha2)+0.85*(rho)$); 	
	\draw[-latex, black, thick] ($-2*(alpha2)+(alpha1)+0.15*(rho)$) to ($-2*(alpha2)+(alpha1)+0.85*(rho)$); 	
	\draw[-latex, black, thick] ($(rho)+2*(alpha2)+0.15*(alpha1)$) to ($(rho)+2*(alpha2)+0.85*(alpha1)$);		 
	\draw[-latex, black, thick] ($(rho)+2*(alpha1)+0.15*(alpha2)$) to ($(rho)+2*(alpha1)+0.85*(alpha2)$);		 
	\draw[-latex, black, thick] ($2*(rho)+0.85*(alpha2)$) to ($2*(rho)+0.15*(alpha2)$);									 
	\draw[-latex, black, thick] ($2*(rho)+0.85*(alpha1)$) to ($2*(rho)+0.15*(alpha1)$);									 
	\draw[-latex, black, thick]  ($-2*(rho)-0.85*(alpha1)$) to ($-2*(rho)-0.15*(alpha1)$);									
	\draw[-latex, black, thick]  ($-2*(rho)-0.85*(alpha2)$) to ($-2*(rho)-0.15*(alpha2)$);									
	\draw[-latex, black, thick]  ($-2*(rho)-(alpha1)+0.15*(alpha2)$) to  ($-2*(rho)-(alpha1)+0.85*(alpha2)$);	
	\draw[-latex, black, thick]  ($-2*(rho)-(alpha2)+0.15*(alpha1)$) to  ($-2*(rho)-(alpha2)+0.85*(alpha1)$);	
\end{tikzpicture}
}
\\
\begin{array}{c}
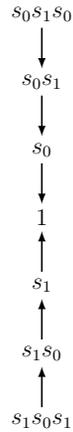
\hbox{postive level Bruhat order for $\widehat{\fsl}_2$} \\
\hbox{1 is minimal}
\end{array}
&\begin{array}{c}
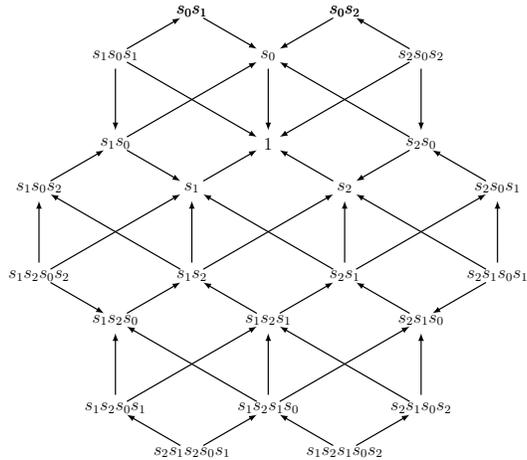
\hbox{postive level Bruhat order for $\widehat{\fsl}_3$} \\
\hbox{1 is minimal}
\end{array}
\end{matrix}
$$

$$
\begin{matrix}
\resizebox{1.2cm}{!}{
\begin{tikzpicture}[scale=1.2]
\tikzstyle{every node}=[font=\large]
	 \coordinate (omega1) at (0,1);	 
	\filldraw (0,0) node {$1$};    
	\filldraw (0,1) node {$s_0$};  
	\filldraw (0,2) node {$s_0s_1$};  
	\filldraw (0,3) node {$s_0s_1s_0$};  
	\filldraw (0,-3) node {$s_1s_0s_1$};  
	\filldraw (0,-2) node {$s_1s_0$};  
	\filldraw (0,-1) node {$s_1$};  

	\draw[-latex, red, thick]  (0,2.8) to (0,2.2);
	\draw[-latex, red, thick]  (0,1.8) to (0,1.2);
	\draw[-latex, red, thick]  (0,0.8) to (0,0.2);
	\draw[-latex, red, thick]  (0,-0.2) to (0,-0.8);
	\draw[-latex, red, thick]  (0,-1.2) to (0,-1.8);
	\draw[-latex, red, thick]  (0,-2.2) to (0,-2.8);
\end{tikzpicture}
}
&\resizebox{7.2cm}{!}{
\begin{tikzpicture}[scale=1.2]
\tikzstyle{every node}=[font=\large]
	 \coordinate (omega2) at (-0.5,0.86603);
	 \coordinate (2omega2) at  (-1,2*0.86603);
	 \coordinate (2omega2+omega1) at (-0.5, 3*0.86603);
	 \coordinate (omega2+2omega1) at (0.5,3*0.86603);
	 \coordinate (2omega1) at (1,2*0.86603);
	 \coordinate (omega1) at (0.5,0.86603);
	 \coordinate (-omega1) at (-0.5,-0.86603);
	 \coordinate (-omega1+omega2) at (-1,0);
	 \coordinate (-omega2+omega1) at (1,0);
	 \coordinate (-omega2) at  (0.5,-0.86603);
	 \coordinate (-2omega1-omega2) at (-0.5,-3*0.86603);
	 \coordinate (-2omega1) at (-1, -2*0.86603);
	 \coordinate (-2omega2) at (1,-2*0.86603);
	 \coordinate (-2omega2-omega1) at (0.5,-3*0.86603);
	 \coordinate (rho) at ($(omega1)+(omega2)$);
	 \coordinate (-rho) at ($-1*(omega1)-(omega2)$);
	 \coordinate (alpha1) at ($2*(omega1)-(omega2)$);
	 \coordinate (-alpha1) at ($-2*(omega1)+(omega2)$);
	 \coordinate (alpha2) at ($2*(omega2)-(omega1)$);
	 \coordinate (-alpha2) at ($-2*(omega2)+(omega1)$);	 
	 
	\filldraw (rho) node {$1$};    
	\filldraw (alpha1) node {$s_2$};  
	\filldraw (alpha2) node {$s_1$};  
	\filldraw (-alpha1) node {$s_1s_2$};  
	\filldraw (-alpha2) node {$s_2s_1$};  
	\filldraw (-rho) node {$s_1s_2s_1$}; 
	\filldraw ($2*(rho)+(alpha1)$) node {$s_0s_2$};    
	\filldraw ($2*(rho)$) node {$s_0$};  
	\filldraw ($2*(alpha1)$) node {$s_2s_0$};  
	\filldraw ($2*(alpha1)+(rho)$) node {$s_2s_0s_2$};  
	\filldraw ($2*(rho)+(alpha2)$) node {$s_0s_1$};    
	\filldraw ($2*(alpha2)$) node {$s_1s_0$};  
	\filldraw ($2*(alpha2)+(rho)$) node {$s_1s_0s_1$};  
	\filldraw ($2*(rho)+(alpha2)$) node {$s_0s_1$};    
	\filldraw ($-2*(alpha1)$) node {$s_1s_2s_0$};  
	\filldraw ($-2*(alpha1)+(alpha2)$) node {$s_1s_2s_0s_2$};  
	\filldraw ($-1*(alpha1)+2*(alpha2)$) node {$s_1s_0s_2$};  
	\filldraw ($2*(rho)+(alpha1)$) node {$s_0s_2$};    
	\filldraw ($-2*(alpha2)$) node {$s_2s_1s_0$};  
	\filldraw ($-2*(alpha2)+(alpha1)$) node {$s_2s_1s_0s_1$};  
	\filldraw ($-1*(alpha2)+2*(alpha1)$) node {$s_2s_0s_1$};  
	\filldraw ($-3*(rho)+(alpha1)$) node {$s_1s_2s_1s_0s_2$};    
	\filldraw ($-2*(rho)$) node {$s_1s_2s_1s_0$};  
	\filldraw ($2*(alpha1)-3*(rho)$) node {$s_2s_1s_0s_2$};  
	\filldraw ($-3*(rho)+(alpha2)$) node {$s_2s_1s_2s_0s_1$};    
	\filldraw ($2*(alpha2)-3*(rho)$) node {$s_1s_2s_0s_1$};  	
	
	\draw[-latex, red, thick]  ($(alpha1)-0.15*(alpha1)$) to ($(-alpha1)+0.15*(alpha1)$);
	\draw[-latex, red, thick]  ($2*(alpha1)+(rho)-0.15*(alpha1)$) to ($(rho)+0.15*(alpha1)$); 		
	\draw[-latex, red, thick]  ($2*(rho)-0.15*(alpha1)$) to ($2*(alpha2)+0.15*(alpha1)$); 		
	\draw[-latex, red, thick]  ($(alpha2)-0.15*(alpha1)$) to ($-2*(alpha1)+(alpha2)+0.15*(alpha1)$); 
	\draw[-latex, red, thick]  ($-1*(alpha2)+2*(alpha1)-0.15*(alpha1)$) to ($(-alpha2)+0.15*(alpha1)$); 
	\draw[-latex, red, thick]  ($(-rho)-0.15*(alpha1)$) to ($2*(alpha2)-3*(rho)+0.15*(alpha1)$); 
	\draw[-latex, red, thick]  ($-2*(alpha2)-0.15*(alpha1)$) to ($-2*(rho)+0.15*(alpha1)$); 
	\draw[-latex, red, thick]  ($(alpha2)-0.15*(alpha2)$) to ($(-alpha2)+0.15*(alpha2)$);
	\draw[-latex, red, thick]  ($2*(alpha2)+(rho)-0.15*(alpha2)$) to ($(rho)+0.15*(alpha2)$); 
	\draw[-latex, red, thick]  ($2*(rho)-0.15*(alpha2)$) to ($2*(alpha1)+0.15*(alpha2)$); 		
	\draw[-latex, red, thick]  ($-1*(alpha1)+2*(alpha2)-0.15*(alpha2)$) to ($(-alpha1)+0.15*(alpha2)$);
 	\draw[-latex, red, thick]  ($(alpha1)-0.15*(alpha2)$) to ($-2*(alpha2)+(alpha1)+0.15*(alpha2)$); 
	\draw[-latex, red, thick]  ($-2*(alpha1)-0.15*(alpha2)$) to ($-2*(rho)+0.15*(alpha2)$); 
	\draw[-latex, red, thick]  ($(-rho)-0.15*(alpha2)$) to ($2*(alpha1)-3*(rho)+0.15*(alpha2)$); 
	\draw[-latex, red, thick]    ($(rho)-0.1*(alpha2)$) to ($(alpha1)+0.1*(alpha2)$); 			
	\draw[-latex, red, thick]   ($(rho)-0.1*(alpha1)$) to ($(alpha2)+0.1*(alpha1)$); 			
	\draw[-latex, red, thick]    ($1.9*(rho)$) to ($1.1*(rho)$); 											
	\draw[-latex, red, thick]   ($(alpha2)-0.1*(rho)$) to ($(-alpha1)+0.1*(rho)$); 				
	\draw[-latex, red, thick]   ($(alpha1)-0.1*(rho)$) to ($(-alpha2)+0.1*(rho)$); 				
	\draw[-latex, red, thick]   ($(-alpha1)-0.1*(alpha2)$) to ($-1*(rho)+0.1*(alpha2)$); 	
	\draw[-latex, red, thick]   ($(-alpha2)-0.1*(alpha1)$) to ($-1*(rho)+0.1*(alpha1)$); 	
	\draw[-latex, red, thick]   ($-1.1*(rho)$) to ($-1.9*(rho)$); 										
	\draw[-latex, red, thick]   ($1.9*(alpha1)$) to ($1.1*(alpha1)$);									
	\draw[-latex, red, thick]   ($1.9*(alpha2)$) to ($1.1*(alpha2)$);									
	\draw[-latex, red, thick]   ($-1.1*(alpha1)$) to ($-1.9*(alpha1)$);								
	\draw[-latex, red, thick]   ($-1.1*(alpha2)$) to ($-1.9*(alpha2)$);								
	\draw[-latex, red, thick]   ($2*(alpha2)+0.9*(rho)$) to ($2*(alpha2)+0.1*(rho)$);		
	\draw[-latex, red, thick]   ($2*(alpha1)+0.9*(rho)$) to ($2*(alpha1)+0.1*(rho)$);		
	\draw[-latex, red, thick]  ($2*(-alpha1)-0.1*(rho)$) to ($2*(-alpha1)-0.9*(rho)$);		
	\draw[-latex, red, thick]   ($2*(-alpha2)-0.1*(rho)$) to ($2*(-alpha2)-0.9*(rho)$);								
	\draw[-latex, red, thick]  ($2*(-alpha1)+0.9*(alpha2)$) to ($2*(-alpha1)+0.1*(alpha2)$);					
	\draw[-latex, red, thick]   ($2*(-alpha2)+0.9*(alpha1)$) to ($2*(-alpha2)+0.1*(alpha1)$);					
	\draw[-latex, red, thick]  ($2*(alpha2)-0.1*(alpha1)$) to ($2*(alpha2)-0.9*(alpha1)$); 						
	\draw[-latex, red, thick]  ($2*(alpha1)-0.1*(alpha2)$) to ($2*(alpha1)-0.9*(alpha2)$); 						
	\draw[-latex, red, thick]  ($-2*(alpha1)+(alpha2)+0.9*(rho)$) to ($-2*(alpha1)+(alpha2)+0.1*(rho)$); 	
	\draw[-latex, red, thick]  ($-2*(alpha2)+(alpha1)+0.9*(rho)$) to ($-2*(alpha2)+(alpha1)+0.1*(rho)$); 	
	\draw[-latex, red, thick]  ($(rho)+2*(alpha2)+0.9*(alpha1)$) to ($(rho)+2*(alpha2)+0.1*(alpha1)$);		 
	\draw[-latex, red, thick]  ($(rho)+2*(alpha1)+0.9*(alpha2)$) to ($(rho)+2*(alpha1)+0.1*(alpha2)$);		 
	\draw[-latex, red, thick]  ($2*(rho)+0.9*(alpha2)$) to ($2*(rho)+0.1*(alpha2)$);									 
	\draw[-latex, red, thick]  ($2*(rho)+0.9*(alpha1)$) to ($2*(rho)+0.1*(alpha1)$);									 	
	\draw[-latex,red, thick]   ($-2*(rho)-0.1*(alpha1)$) to ($-2*(rho)-0.9*(alpha1)$);									
	\draw[-latex, red, thick]  ($-2*(rho)-0.1*(alpha2)$) to ($-2*(rho)-0.9*(alpha2)$);									
	\draw[-latex, red, thick]   ($-2*(rho)-(alpha1)+0.9*(alpha2)$) to ($-2*(rho)-(alpha1)+0.1*(alpha2)$);	
	\draw[-latex, red, thick]  ($-2*(rho)-(alpha2)+0.9*(alpha1)$) to ($-2*(rho)-(alpha2)+0.1*(alpha1)$);	
\end{tikzpicture}
}
\\
\begin{array}{c}
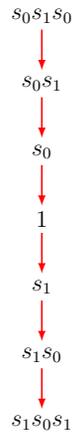
\hbox{level zero Bruhat order for $\widehat{\fsl}_2$} \\
\hbox{translation invariant}
\end{array}
&\begin{array}{c}
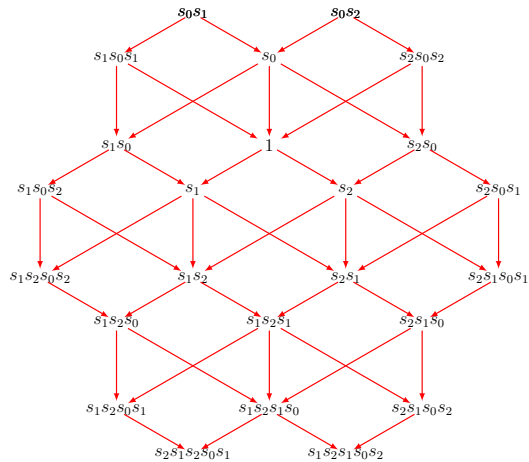
\hbox{level zero Bruhat order for $\widehat{\fsl}_3$} \\
\hbox{translation invariant}
\end{array}
\end{matrix}
$$

$$
\begin{matrix}
\resizebox{1.2cm}{!}{
\begin{tikzpicture}[scale=1.2]
\tikzstyle{every node}=[font=\large]
	 \coordinate (omega1) at (0,1);	 
	\filldraw (0,0) node {$1$};    
	\filldraw (0,1) node {$s_0$};  
	\filldraw (0,2) node {$s_0s_1$};  
	\filldraw (0,3) node {$s_0s_1s_0$};  
	\filldraw (0,-3) node {$s_1s_0s_1$};  
	\filldraw (0,-2) node {$s_1s_0$};  
	\filldraw (0,-1) node {$s_1$};  

	\draw[-latex, blue, thick]  (0,2.2) to (0,2.8);
	\draw[-latex, blue, thick]  (0,1.2) to (0,1.8);
	\draw[-latex, blue, thick]  (0,0.2) to (0,0.8);
	\draw[-latex, blue, thick]  (0,-0.2) to (0,-0.8);
	\draw[-latex, blue, thick]  (0,-1.2) to (0,-1.8);
	\draw[-latex, blue, thick]  (0,-2.2) to (0,-2.8);
\end{tikzpicture}
}
&\resizebox{7.2cm}{!}{
\begin{tikzpicture}[scale=1.2]
\tikzstyle{every node}=[font=\large]
	 \coordinate (omega2) at (-0.5,0.86603);
	 \coordinate (2omega2) at  (-1,2*0.86603);
	 \coordinate (2omega2+omega1) at (-0.5, 3*0.86603);
	 \coordinate (omega2+2omega1) at (0.5,3*0.86603);
	 \coordinate (2omega1) at (1,2*0.86603);
	 \coordinate (omega1) at (0.5,0.86603);
	 \coordinate (-omega1) at (-0.5,-0.86603);
	 \coordinate (-omega1+omega2) at (-1,0);
	 \coordinate (-omega2+omega1) at (1,0);
	 \coordinate (-omega2) at  (0.5,-0.86603);
	 \coordinate (-2omega1-omega2) at (-0.5,-3*0.86603);
	 \coordinate (-2omega1) at (-1, -2*0.86603);
	 \coordinate (-2omega2) at (1,-2*0.86603);
	 \coordinate (-2omega2-omega1) at (0.5,-3*0.86603);
	 \coordinate (rho) at ($(omega1)+(omega2)$);
	 \coordinate (-rho) at ($-1*(omega1)-(omega2)$);
	 \coordinate (alpha1) at ($2*(omega1)-(omega2)$);
	 \coordinate (-alpha1) at ($-2*(omega1)+(omega2)$);
	 \coordinate (alpha2) at ($2*(omega2)-(omega1)$);
	 \coordinate (-alpha2) at ($-2*(omega2)+(omega1)$);	 
	 
	\filldraw (rho) node {$1$};    
	\filldraw (alpha1) node {$s_2$};  
	\filldraw (alpha2) node {$s_1$};  
	\filldraw (-alpha1) node {$s_1s_2$};  
	\filldraw (-alpha2) node {$s_2s_1$};  
	\filldraw (-rho) node {$s_1s_2s_1$}; 
	\filldraw ($2*(rho)+(alpha1)$) node {$s_0s_2$};    
	\filldraw ($2*(rho)$) node {$s_0$};  
	\filldraw ($2*(alpha1)$) node {$s_2s_0$};  
	\filldraw ($2*(alpha1)+(rho)$) node {$s_2s_0s_2$};  
	\filldraw ($2*(rho)+(alpha2)$) node {$s_0s_1$};    
	\filldraw ($2*(alpha2)$) node {$s_1s_0$};  
	\filldraw ($2*(alpha2)+(rho)$) node {$s_1s_0s_1$};  
	\filldraw ($2*(rho)+(alpha2)$) node {$s_0s_1$};    
	\filldraw ($-2*(alpha1)$) node {$s_1s_2s_0$};  
	\filldraw ($-2*(alpha1)+(alpha2)$) node {$s_1s_2s_0s_2$};  
	\filldraw ($-1*(alpha1)+2*(alpha2)$) node {$s_1s_0s_2$};  
	\filldraw ($2*(rho)+(alpha1)$) node {$s_0s_2$};    
	\filldraw ($-2*(alpha2)$) node {$s_2s_1s_0$};  
	\filldraw ($-2*(alpha2)+(alpha1)$) node {$s_2s_1s_0s_1$};  
	\filldraw ($-1*(alpha2)+2*(alpha1)$) node {$s_2s_0s_1$};  
	\filldraw ($-3*(rho)+(alpha1)$) node {$s_1s_2s_1s_0s_2$};    
	\filldraw ($-2*(rho)$) node {$s_1s_2s_1s_0$};  
	\filldraw ($2*(alpha1)-3*(rho)$) node {$s_2s_1s_0s_2$};  
	\filldraw ($-3*(rho)+(alpha2)$) node {$s_2s_1s_2s_0s_1$};    
	\filldraw ($2*(alpha2)-3*(rho)$) node {$s_1s_2s_0s_1$};  	
	
	\draw[-latex, blue, thick]  ($(alpha1)-0.15*(alpha1)$) to ($(-alpha1)+0.15*(alpha1)$); 				\draw[-latex, blue, thick]  ($(rho)+0.15*(alpha1)$) to ($2*(alpha1)+(rho)-0.15*(alpha1)$); 		
	\draw[-latex, blue, thick]  ($2*(rho)-0.15*(alpha1)$) to ($2*(alpha2)+0.15*(alpha1)$); 		
	\draw[-latex, blue, thick]  ($(alpha2)-0.15*(alpha1)$) to ($-2*(alpha1)+(alpha2)+0.15*(alpha1)$); 
	\draw[-latex, blue, thick]  ($-1*(alpha2)+2*(alpha1)-0.15*(alpha1)$) to ($(-alpha2)+0.15*(alpha1)$); 
	\draw[-latex, blue, thick]  ($(-rho)-0.15*(alpha1)$) to ($2*(alpha2)-3*(rho)+0.15*(alpha1)$); 
	\draw[-latex, blue, thick]  ($-2*(alpha2)-0.15*(alpha1)$) to ($-2*(rho)+0.15*(alpha1)$); 
	\draw[-latex, blue, thick]  ($(alpha2)-0.15*(alpha2)$) to ($(-alpha2)+0.15*(alpha2)$);
	\draw[-latex, blue, thick]  ($(rho)+0.15*(alpha2)$) to ($2*(alpha2)+(rho)-0.15*(alpha2)$); 
	\draw[-latex, blue, thick]  ($2*(rho)-0.15*(alpha2)$) to ($2*(alpha1)+0.15*(alpha2)$); 		
	\draw[-latex, blue, thick]  ($-1*(alpha1)+2*(alpha2)-0.15*(alpha2)$) to ($(-alpha1)+0.15*(alpha2)$);
 	\draw[-latex, blue, thick]  ($(alpha1)-0.15*(alpha2)$) to ($-2*(alpha2)+(alpha1)+0.15*(alpha2)$); 
	\draw[-latex, blue, thick]  ($-2*(alpha1)-0.15*(alpha2)$) to ($-2*(rho)+0.15*(alpha2)$); 
	\draw[-latex, blue, thick]  ($(-rho)-0.15*(alpha2)$) to ($2*(alpha1)-3*(rho)+0.15*(alpha2)$); 
	\draw[-latex, blue, thick]    ($(rho)-0.1*(alpha2)$) to ($(alpha1)+0.1*(alpha2)$); 			
	\draw[-latex, blue, thick]	   ($(rho)-0.1*(alpha1)$) to ($(alpha2)+0.1*(alpha1)$); 			
	\draw[-latex, blue, thick]   ($1.1*(rho)$) to  ($1.9*(rho)$); 											
	\draw[-latex, blue, thick]   ($(alpha2)-0.1*(rho)$) to ($(-alpha1)+0.1*(rho)$); 				
	\draw[-latex, blue, thick]   ($(alpha1)-0.1*(rho)$) to ($(-alpha2)+0.1*(rho)$); 				
	\draw[-latex, blue, thick]   ($(-alpha1)-0.1*(alpha2)$) to ($-1*(rho)+0.1*(alpha2)$); 	
	\draw[-latex, blue, thick]   ($(-alpha2)-0.1*(alpha1)$) to ($-1*(rho)+0.1*(alpha1)$); 	
	\draw[-latex, blue, thick]   ($-1.1*(rho)$) to ($-1.9*(rho)$); 										
	\draw[-latex, blue, thick]   ($1.1*(alpha1)$) to ($1.9*(alpha1)$);									
	\draw[-latex, blue, thick]   ($1.1*(alpha2)$) to ($1.9*(alpha2)$);									
	\draw[-latex, blue, thick]   ($-1.1*(alpha1)$) to ($-1.9*(alpha1)$);								
	\draw[-latex, blue, thick]   ($-1.1*(alpha2)$) to ($-1.9*(alpha2)$);								
	\draw[-latex, blue, thick]   ($2*(alpha2)+0.1*(rho)$) to ($2*(alpha2)+0.9*(rho)$);		
	\draw[-latex, blue, thick]   ($2*(alpha1)+0.1*(rho)$) to ($2*(alpha1)+0.9*(rho)$);		
	\draw[-latex, blue, thick]   ($2*(-alpha1)-0.1*(rho)$) to ($2*(-alpha1)-0.9*(rho)$);		
	\draw[-latex, blue, thick]   ($2*(-alpha2)-0.1*(rho)$) to ($2*(-alpha2)-0.9*(rho)$);								
	\draw[-latex, blue, thick]   ($2*(-alpha1)+0.1*(alpha2)$) to ($2*(-alpha1)+0.9*(alpha2)$);					
	\draw[-latex, blue, thick]   ($2*(-alpha2)+0.1*(alpha1)$) to ($2*(-alpha2)+0.9*(alpha1)$);					
	\draw[-latex, blue, thick]  ($2*(alpha2)-0.1*(alpha1)$) to ($2*(alpha2)-0.9*(alpha1)$); 						
	\draw[-latex, blue, thick]  ($2*(alpha1)-0.1*(alpha2)$) to ($2*(alpha1)-0.9*(alpha2)$); 						
	\draw[-latex, blue, thick]  ($-2*(alpha1)+(alpha2)+0.9*(rho)$) to ($-2*(alpha1)+(alpha2)+0.1*(rho)$); 	
	\draw[-latex, blue, thick]  ($-2*(alpha2)+(alpha1)+0.9*(rho)$) to ($-2*(alpha2)+(alpha1)+0.1*(rho)$); 	
	\draw[-latex, blue, thick]  ($(rho)+2*(alpha2)+0.9*(alpha1)$) to ($(rho)+2*(alpha2)+0.1*(alpha1)$);		 
	\draw[-latex, blue, thick]  ($(rho)+2*(alpha1)+0.9*(alpha2)$) to ($(rho)+2*(alpha1)+0.1*(alpha2)$);		 
	\draw[-latex, blue, thick]  ($2*(rho)+0.1*(alpha2)$) to ($2*(rho)+0.9*(alpha2)$);									 
	\draw[-latex, blue, thick]  ($2*(rho)+0.1*(alpha1)$) to ($2*(rho)+0.9*(alpha1)$);									 
	\draw[-latex, blue, thick]   ($-2*(rho)-0.1*(alpha1)$) to ($-2*(rho)-0.9*(alpha1)$);									
	\draw[-latex, blue, thick]   ($-2*(rho)-0.1*(alpha2)$) to ($-2*(rho)-0.9*(alpha2)$);									
	\draw[-latex, blue, thick]   ($-2*(rho)-(alpha1)+0.9*(alpha2)$) to ($-2*(rho)-(alpha1)+0.1*(alpha2)$);	
	\draw[-latex, blue, thick]   ($-2*(rho)-(alpha2)+0.9*(alpha1)$) to ($-2*(rho)-(alpha2)+0.1*(alpha1)$);	
\end{tikzpicture}
}
\\
\begin{array}{c}
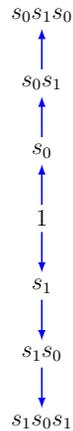
\hbox{negative level Bruhat order for $\widehat{\fsl}_2$} \\
\hbox{1 is maximal}
\end{array}
&\begin{array}{c}
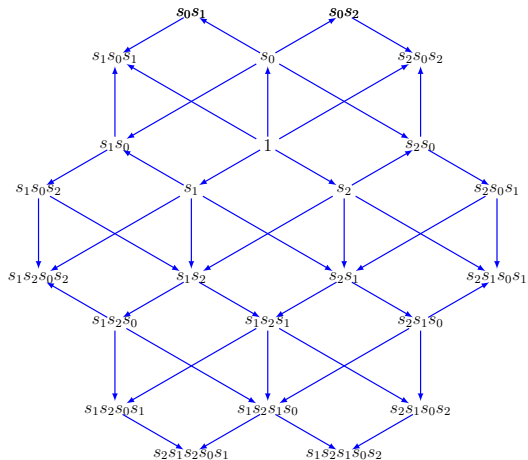
\hbox{negative level Bruhat order for $\widehat{\fsl}_3$} \\
\hbox{1 is maximal}
\end{array}
\end{matrix}
$$

\newpage
\noindent
\textbf{PLATE B: Pictures of $B(\omega_1+\Lambda_0)$, $B(\omega_1+0\Lambda_0)$
and $B(-\omega_1-\Lambda_0)$ for $\widehat{\fsl}_2$}

$$
\begin{matrix}
\begin{tikzpicture}[scale=0.7]
\tikzstyle{every node}=[font=\tiny]    
    \draw[style=help lines,step=1cm] (-10,-6) grid (10,1);   
    \draw[->] (-10,0) -- (10.5,0) node[anchor=west]{$\omega_1$-axis};
    \draw[->] (0,-6) -- (0, 1.5) node[anchor=south]{$\delta$-axis};

     \node at (1,0) {$\bullet$};
        \filldraw (1,1) node[anchor=north,yshift=-0.1cm] {$u_\Lambda$};
            \node at (-1,0) {$\bullet$};
 \filldraw (-1,1) node[anchor=north,yshift=-0.1cm] {$u_{s_1\Lambda}$};



        \node at (3, -2) {$\bullet$}; \node at (-3, -2) {$\bullet$};
        \filldraw (3,-2) node[anchor=west,yshift=0cm] {$u_{s_0s_1\Lambda}$};
                \node at (-3, -2) {$\bullet$};
        \filldraw (-3,-2) node[anchor=east,yshift=0cm] {$u_{s_1s_0s_1\Lambda}$};
        \node at (5, -6) {$\bullet$};
        \filldraw (5,-6) node[anchor=west,yshift=-0.1cm] {$u_{s_0s_1s_0s_1\Lambda}$};
                \node at (-3, -2) {$\bullet$};
		\draw (0, 0.3)parabola (-5, -6) ;
		\draw (0, 0.3)parabola (5, -6) ;
\draw [<-,blue, out=30, in=150, thick] (-1,0) to node[left, xshift=3, yshift=6]{$\tilde{f}_1$} (1,0);
\node at (1, -1) {$\bullet$}; 
\draw [<-,red, in=0, out=120,  thick] (1,-1) to node[right, xshift=5, yshift=-1]{$\tilde{f}_0$} (-1,0);
\node at (3, -2) {$\bullet$}; 
\draw [<-,red, in=0, out=120,  thick] (3,-2) to node[right, xshift=6, yshift=-1]{$\tilde{f}_0$} (1,-1);
\node at (-1, -1) {$\bullet$}; 
\draw [<-,blue, out=30, in=150, thick] (-1,-1) to node[left, yshift=4]{$\tilde{f}_1$} (1,-1);
\node at (1, -1) {$\bullet$}; 
\draw [<-,red, in=0, out=120, thick] (1,-2.2) to node[right, xshift=5, yshift=-1]{$\tilde{f}_0$} (-1,-1);
\node at (3, -3) {$\bullet$}; 
\draw [<-,red, in=0, out=120,  thick] (3,-3) to node[right, xshift=6, yshift=-1]{$\tilde{f}_0$} (1,-2.2);
\node at (1, -2.2) {$\bullet$}; 
\draw [<-,blue, out=30, in=150, thick] (1,-2) to node[left, yshift=5]{$\tilde{f}_1$} (3,-2);
\node at (1, -2) {$\bullet$}; 
\draw [<-,blue, out=30, in=150, thick] (-1,-2) to node[left, xshift=-8,yshift=2]{$\tilde{f}_1$} (1,-2);
\node at (1, -2) {$\bullet$}; 
\draw [<-,blue, out=30, in=150, thick] (-1,-2.2) to node[left, xshift=-6,yshift=-10]{$\tilde{f}_1$} (1,-2.2);
\node at (-1, -2.2) {$\bullet$}; 
\draw [<-,blue, out=30, in=150, thick] (-3,-2) to node[left, yshift=5]{$\tilde{f}_1$} (-1,-2);
\node at (-3, -2) {$\bullet$}; 
\draw [<-,red, in=20, out=120, thick] (1,-3.2) to node[right, xshift=-7, yshift=-5]{$\tilde{f}_0$} (-1,-2.2);
\node at (1, -3.2) {$\bullet$}; 
\draw [<-,red, in=0, out=120, thick] (1,-3) to node[right, xshift=6, yshift=-2.6]{$\tilde{f}_0$} (-1,-2);
\node at (1, -3) {$\bullet$}; 
\draw [<-,red, in=0, out=120, thick] (3,-4.2) to node[right, xshift=-7, yshift=-5]{$\tilde{f}_0$} (1,-3.2);
\node at (3, -4.2) {$\bullet$}; 
\node at (3, -4) {$\bullet$}; 
\node at (-1, -2) {$\bullet$}; 
\draw [<-,red, in=0, out=120, thick] (3,-4) to node[right, xshift=6, yshift=-1]{$\tilde{f}_0$} (1,-3);
\node at (1, -4) {$\bullet$}; 
\draw [<-,red, in=0, out=120, thick] (-1,-3) to node[right, xshift=-6, yshift=-5.5]{$\tilde{f}_0$} (-3,-2);
\node at (-1, -3) {$\bullet$}; 
\draw [<-,red, in=0, out=120, thick] (1,-4) to node[right, xshift=-6, yshift=-5.5]{$\tilde{f}_0$} (-1,-3);
\node at (3, -5) {$\bullet$}; 
\draw [<-,red, in=0, out=120, thick] (3,-5) to node[right, xshift=-6, yshift=-5.5]{$\tilde{f}_0$} (1,-4);
\node at (5, -6) {$\bullet$}; 
\draw [<-,red, in=0, out=120, thick] (5,-6) to node[right, xshift=-6, yshift=-5.5]{$\tilde{f}_0$} (3,-5);

\end{tikzpicture}
\\
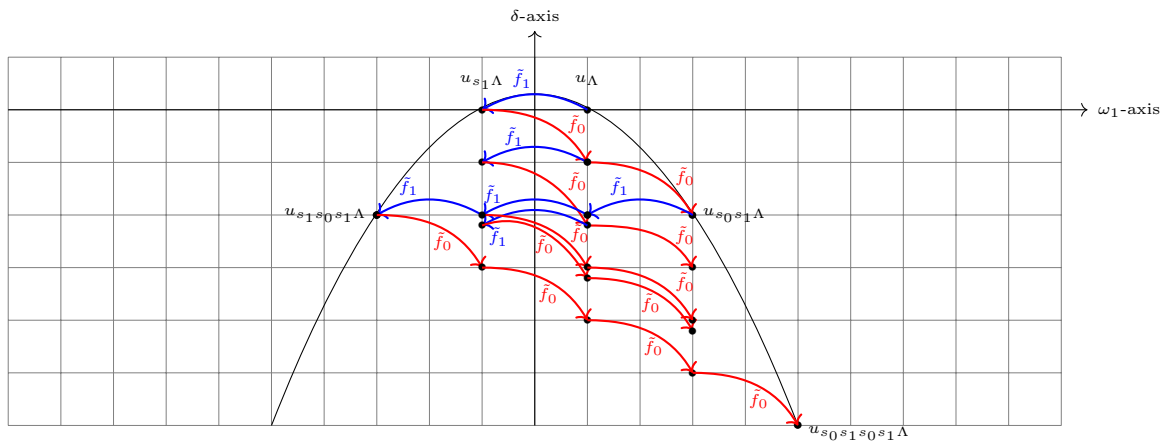
\hbox{Initial portion of the crystal graph of $B(\omega_1+\Lambda_0)$ for $\widehat{\fsl}_2$}
\end{matrix}
$$

$$
\begin{matrix}
\begin{tikzpicture}[scale=0.7]
\tikzstyle{every node}=[font=\tiny]    
    \draw[style=help lines,step=1cm] (-10,6) grid (10,-1);   
    \draw[->] (-10,0) -- (10.5,0) node[anchor=west]{$\omega_1$-axis};
    \draw[->] (0,6) -- (0, -1.5) node[anchor=north]{$\delta$-axis};

     \node at (1,0) {$\bullet$};
        \filldraw (1,0) node[anchor=north,yshift=-0.1cm] {$u_{s_1\Lambda}$};
            \node at (-1,0) {$\bullet$};
 \filldraw (-1,0) node[anchor=north,yshift=-0.1cm] {$u_{\Lambda}$};



        \node at (3, 2) {$\bullet$}; \node at (-3, 2) {$\bullet$};
        \filldraw (3,2) node[anchor=west,yshift=0cm] {$u_{s_1s_0s_1\Lambda}$};
                \node[red] at (-3, 2) {$\bullet$};
        \filldraw (-3,2) node[anchor=east,yshift=0cm] {$u_{s_0s_1\Lambda}$};
     \node at (-5, 6) {$\bullet$};
                \node at (-3, 2) {$\bullet$};
        \filldraw (-5,6) node[anchor=east,yshift=-0.1cm] {$u_{s_0s_1s_0s_1\Lambda}$};
		\draw (0, -0.3)parabola (-5, 6) ;
		\draw (0, -0.3)parabola (5, 6) ;
\draw [<-,blue, out=30, in=150, thick] (-1,0) to node[left, yshift=3, xshift=-4]{$\tilde{f}_1$} (1,0);
\draw [<-,red, in=0, out=120, thick] (1,0) to node[right, xshift=5, yshift=-2]{$\tilde{f}_0$} (-1,1);
\node at (-1, 1) {$\bullet$}; 
\draw [<-,red, in=0, out=120, thick] (-1,1) to node[right, xshift=5, yshift=-2]{$\tilde{f}_0$} (-3,2);
\node at (-3, 2) {$\bullet$}; 
\draw [<-,blue, out=30, in=150, thick] (-1,1) to node[left, yshift=2,xshift=-7]{$\tilde{f}_1$} (1,1);
\node at (1, 1) {$\bullet$}; 
\draw [<-,red, in=0, out=120, thick] (1,1) to node[right, xshift=-4, yshift=-7]{$\tilde{f}_0$} (-1,2.2);
\node at (-1, 2.2) {$\bullet$}; 
\draw [<-,red, in=0, out=120,  thick] (-1,2.2) to node[right, xshift=-1, yshift=-8]{$\tilde{f}_0$} (-3,3);
\node at (-3, 3) {$\bullet$}; 
\draw [<-,red, in=0, out=110,  thick] (1,2) to node[right, xshift=-5, yshift=-6]{$\tilde{f}_0$} (-1,3);
\draw [<-,red, in=0, out=120,  thick] (-1,3) to node[right, xshift=0, yshift=-11]{$\tilde{f}_0$} (-3,4);
\node at (-3, 4) {$\bullet$}; 
\node at (-3, 4.2) {$\bullet$}; 
\draw [<-,red, in=0, out=120, thick] (1,2.2) to node[right, xshift=6, yshift=-3.6]{$\tilde{f}_0$} (-1,3.2);
\node at (-1, 3.2) {$\bullet$}; 
\node at (-1, 3) {$\bullet$}; 
\draw [<-,red, in=0, out=120, thick] (-1,3.2) to node[right, xshift=6, yshift=-3]{$\tilde{f}_0$} (-3,4.2);
\draw [<-,blue, out=30, in=150, thick] (1,2) to node[left, yshift=3,xshift=-5]{$\tilde{f}_1$} (3,2);
\node at (1, 2) {$\bullet$}; 
\draw [<-,blue, out=10, in=170, thick] (-1,2) to node[left, yshift=-5.5,xshift=-6]{$\tilde{f}_1$} (1,2);
\node at (1, 2.2) {$\bullet$}; 
\draw [<-,blue, out=10, in=170, thick] (-1,2.2) to node[left, yshift=5,xshift=-7]{$\tilde{f}_1$} (1,2.2);
\draw [<-,blue, out=30, in=150, thick] (-3,2) to node[left, yshift=2.5,xshift=-5]{$\tilde{f}_1$} (-1,2);
\node at (-1, 2) {$\bullet$}; 
\draw [<-,red, in=0, out=120, thick] (3,2) to node[right, xshift=5, yshift=-2.5]{$\tilde{f}_0$} (1,3);
\node at (1, 3) {$\bullet$}; 
\draw [<-,red, in=0, out=120, thick] (1,3) to node[right, xshift=5, yshift=-2.5]{$\tilde{f}_0$} (-1,4);
\node at (-1, 4) {$\bullet$}; 
\draw [<-,red, in=0, out=120, thick] (-1,4) to node[right, xshift=5, yshift=-2.5]{$\tilde{f}_0$} (-3,5);
\node at (-3, 5) {$\bullet$}; 
\draw [<-,red, in=0, out=120, thick] (-3,5) to node[right, xshift=5, yshift=-2.5]{$\tilde{f}_0$} (-5,6);
\node at (-5, 6) {$\bullet$}; 

\end{tikzpicture}
\\
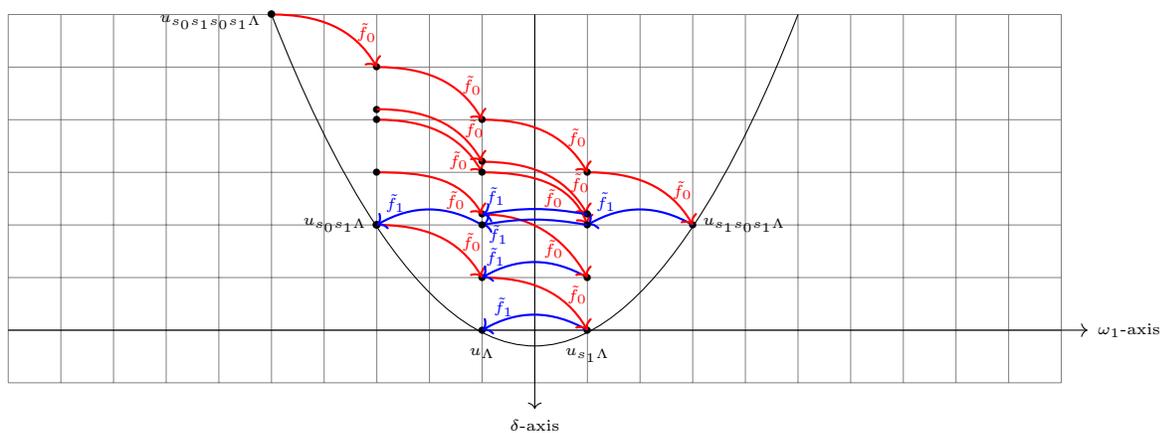
\hbox{Final portion of the crystal graph of $B(-\omega_1-\Lambda_0)$ for $\widehat{\fsl}_2$}
\end{matrix}
$$

$$\begin{matrix}
\begin{tikzpicture}[scale=0.7]
\tikzstyle{every node}=[font=\tiny]
  
    \draw[style=help lines,step=1cm] (-3,-3) grid (3,4);    
    \draw[->] (-3,0) -- (3.5,0) node[anchor=west]{$\omega_1$-axis};
    \draw[->] (0,-3.5) -- (0, 4.5) node[anchor=south]{$\delta$-axis};
      \draw[-,thick] (-1,-3) -- (-1,4) ;
      \draw[-,thick] (1,-3) -- (1,4) ;
	\node at (1, -3) {$\bullet$};  
    \filldraw (1, -3) node[anchor=south west,yshift=-0.1cm] {$u_{(s_0s_1)^3\Lambda}=u_{t_{3\alpha_1^\vee}\Lambda}$};
	\node at (1, -2) {$\bullet$};  
    \filldraw (1, -2) node[anchor=south west,yshift=-0.1cm] {$u_{(s_0s_1)^2\Lambda}=u_{t_{2\alpha_1^\vee}\Lambda}$};
	\node at (1, -1) {$\bullet$};  
    \filldraw (1, -1) node[anchor=south west,yshift=-0.1cm] {$u_{s_0s_1\Lambda} = u_{t_{\alpha_1^\vee}\Lambda}$};
	\node at (1, 0) {$\bullet$};  
    \filldraw (1, 0) node[anchor=south west,yshift=0.05cm] {$u_{\Lambda}=u_{\omega_1+0\Lambda_0}$};
	\node at (1, 1) {$\bullet$};  
    \filldraw (1, 1) node[anchor=south west,yshift=-0.1cm] {$u_{s_1s_0\Lambda} = u_{t_{-\alpha_1^\vee}\Lambda}$};
	\node at (1, 2) {$\bullet$};  
    \filldraw (1, 2) node[anchor=south west,yshift=-0.1cm] {$u_{(s_1s_0)^2\Lambda}=u_{t_{-2\alpha_1^\vee}\Lambda}$};
	\node at (1, 3) {$\bullet$};  
    \filldraw (1, 3) node[anchor=south west,yshift=-0.1cm] {$u_{(s_1s_0)^3\Lambda}=u_{t_{-3\alpha_1^\vee}\Lambda}$};
	\node at (-1, -3) {$\bullet$};  
    \filldraw (-1, -3) node[anchor=south east,yshift=-0.1cm] {$u_{s_1(s_0s_1)^3\Lambda}$};
	\node at (-1, -2) {$\bullet$};  
    \filldraw (-1, -2) node[anchor=south east,yshift=-0.1cm] {$u_{s_1(s_0s_1)^2\Lambda}$};
	\node at (-1, -1) {$\bullet$};  
    \filldraw (-1, -1) node[anchor=south east,yshift=-0.1cm] {$u_{s_1s_0s_1\Lambda}$};
	\node at (-1, 0) {$\bullet$};  
    \filldraw (-1, 0) node[anchor=south east,yshift=-0.1cm] {$u_{s_1\Lambda}$};
	\node at (-1, 1) {$\bullet$};  
    \filldraw (-1, 1) node[anchor=south east,yshift=-0.1cm] {$u_{s_0\Lambda}$};
	\node at (-1, 2) {$\bullet$};  
    \filldraw (-1, 2) node[anchor=south east,yshift=-0.1cm] {$u_{s_0s_1s_0\Lambda}$};
	\node at (-1, 3) {$\bullet$};  
    \filldraw (-1, 3) node[anchor=south east,yshift=-0.1cm] {$u_{s_0(s_1s_0)^2\Lambda}$};
	\foreach \y in {1,...,3}{
	\draw [<-,blue, out=30, in=150, thick] (-1,\y) to node[left, yshift=4,xshift=-4]{$\tilde{f}_1$} (1,\y);
	}	
	\foreach \y in {-3,...,0}{
	\draw [->,blue, out=150, in=30, thick] (1,\y) to node[left, yshift=4,xshift=-4]{$\tilde{f}_1$} (-1,\y);
	}
	\foreach \x in {0,...,2}{
	\draw [<-,red, in=0, out=140, thick] (1,\x) to node[right, xshift=5, yshift=-1]{$\tilde{f}_0$} (-1,\x +1);
	}
	\foreach \x in {-3,...,-1}{
	\draw [->,red, out=0, in=140, thick] (-1,\x+1) to node[right, xshift=5, yshift=-1]{$\tilde{f}_0$} (1,\x	);
	}
    \end{tikzpicture}
\\
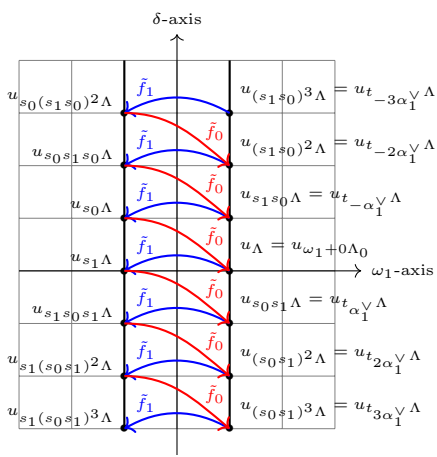
\hbox{Middle portion of the crystal graph of $B(\omega_1+0\Lambda_0)$ for $\fg = \widehat{\fsl}_2$}
\end{matrix}
$$

\newpage

\noindent
\textbf{PLATE C: Pictures for $B(2\omega_1)$.}  Representative paths from the (first five) connected components of $B(2\omega_1)$ are
$$
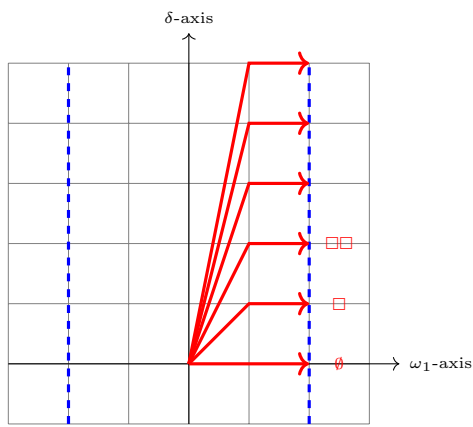
\begin{tikzpicture}[scale=0.8]
\tikzstyle{every node}=[font=\tiny]

\draw[style=help lines,step=1cm] (-3,-1) grid (3,5);    
\draw[->] (-3,0) -- (3.5,0) node[anchor=west]{$\omega_1$-axis};
\draw[->] (0,-1) -- (0, 5.5) node[anchor=south]{$\delta$-axis};
\draw[-,dashed,  blue, very thick] (-2,-1) -- (-2,5) ;
\draw[-,dashed,  blue, very thick] (2,-1) -- (2,5) ;
\filldraw (2.5, 0) node[red] {$\emptyset$};    
\filldraw (2.5, 1) node[red] {$\square$};    
\filldraw (2.5, 2) node[red] {$\square\square$};    
\draw [->,red, very thick] (0,0) to (2, 0);
\draw [->,red, very thick] (0,0) to (1,1) to (2,1);
\draw [->,red, very thick] (0,0) to (1,2) to (2,2);
\draw [->,red, very thick] (0,0) to (1,3) to (2,3);
\draw [->,red, very thick] (0,0) to (1,4) to (2,4);
\draw [->,red, very thick] (0,0) to (1,5) to (2,5);
%
\end{tikzpicture}
$$
and the paths in $B(2\omega_1)_0\subseteq B(\omega_1)\otimes B(\omega_1)$ are
$$
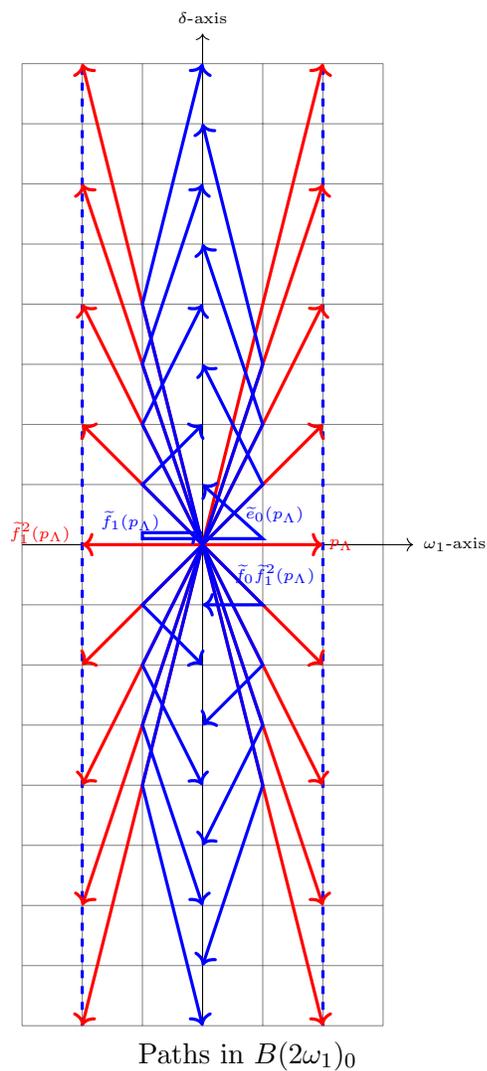
\begin{matrix}
\begin{tikzpicture}[scale=0.8]
\tikzstyle{every node}=[font=\tiny]

\draw[style=help lines,step=1cm] (-3,-8) grid (3,8);    
\draw[->] (-3,0) -- (3.5,0) node[anchor=west]{$\omega_1$-axis};
\draw[->] (0,-8) -- (0, 8.5) node[anchor=south]{$\delta$-axis};
\draw[-,dashed,  blue, very thick] (-2,-8) -- (-2,8) ;
\draw[-,dashed,  blue, very thick] (2,-8) -- (2,8) ;
\filldraw (2.3, 0) node[red] {$p_\Lambda$};    
\filldraw (-2.7, 0.2) node[red] {$\widetilde{f}_1^2(p_\Lambda)$};    
\filldraw (-1.2, 0.4) node[blue] {$\widetilde{f}_1(p_\Lambda)$};    
\filldraw (1.2, 0.5) node[blue] {$\widetilde{e}_0(p_\Lambda)$};   
\filldraw (1.2, -0.5) node[blue] {$\widetilde{f}_0\widetilde{f}_1^2(p_\Lambda)$};   
\draw [->,red, very thick] (0,0) to (2, 8);
\draw [->,red, very thick] (0,0) to (-2, 8);
\draw [->,red, very thick] (0,0) to (2, 6);
\draw [->,red, very thick] (0,0) to (-2, 6);
\draw [->,red, very thick] (0,0) to (2, 4);
\draw [->,red, very thick] (0,0) to (-2, 4);
\draw [->,red, very thick] (0,0) to (2, 2);
\draw [->,red, very thick] (0,0) to (-2, 2);
\draw [->,red, very thick] (0,0) to (2, 0);
\draw [->,red, very thick] (0,0) to (-2, 0);
\draw [->,red, very thick] (0,0) to (2, -2);
\draw [->,red, very thick] (0,0) to (-2, -2);
\draw [->,red, very thick] (0,0) to (2, -4);
\draw [->,red, very thick] (0,0) to (-2, -4);
\draw [->,red, very thick] (0,0) to (2, -6);
\draw [->,red, very thick] (0,0) to (-2, -6);
\draw [->,red, very thick] (0,0) to (2, -8);
\draw [->,red, very thick] (0,0) to (-2, -8);
\draw [->,blue, very thick] (0,0.1) to (-1, 0.1) to (-1,0.2) to (0,0.2);
\draw [->,blue, very thick] (0,0.1) to (1, 0.1) to (0,1);
\draw [->,blue, very thick] (0,0) to (-1,1) to (0,2);
\draw [->,blue, very thick] (0,0) to (1,1) to (0,3);
\draw [->,blue, very thick] (0,0) to (-1,2) to (0,4);
\draw [->,blue, very thick] (0,0) to (1,2) to (0,5);
\draw [->,blue, very thick] (0,0) to (-1,3) to (0,6);
\draw [->,blue, very thick] (0,0) to (1,3) to (0,7);
\draw [->,blue, very thick] (0,0) to (-1,4) to (0,8);
\draw [->,blue, very thick] (0,0) to (1,-1) to (0,-1);
\draw [->,blue, very thick] (0,0) to (-1,-1) to (0,-2);
\draw [->,blue, very thick] (0,0) to (1,-2) to (0,-3);
\draw [->,blue, very thick] (0,0) to (-1,-2) to (0,-4);
\draw [->,blue, very thick] (0,0) to (1,-3) to (0,-5);
\draw [->,blue, very thick] (0,0) to (-1,-3) to (0,-6);
\draw [->,blue, very thick] (0,0) to (1,-4) to (0,-7);
\draw [->,blue, very thick] (0,0) to (-1,-4) to (0,-8);
\end{tikzpicture}
&\qquad\quad
&
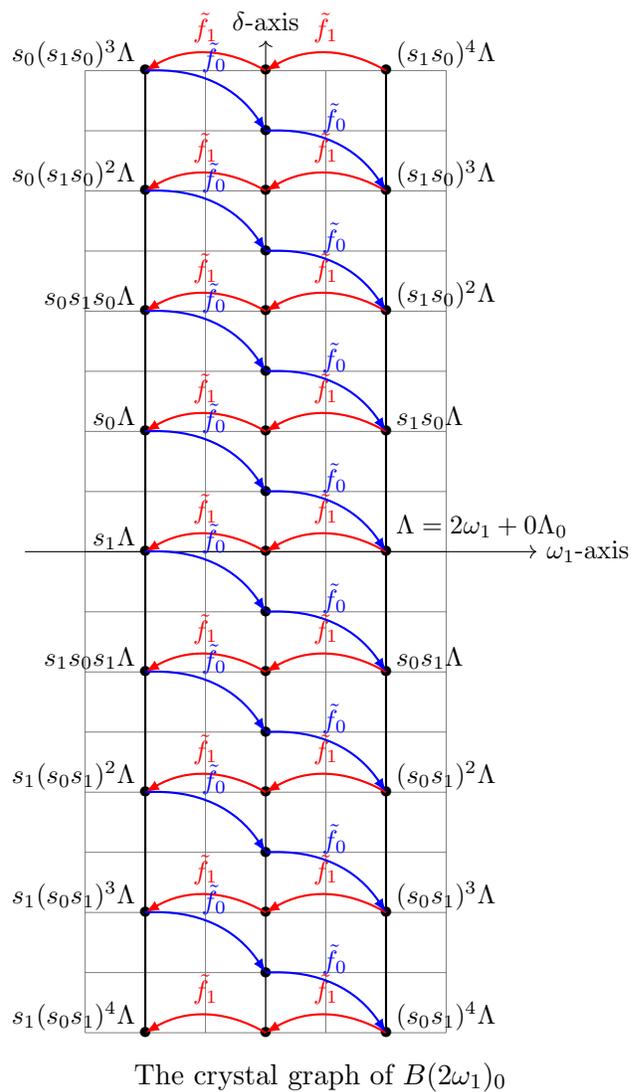
\begin{tikzpicture}[scale=0.8]
\tikzstyle{every node}=[font=\small]
    
    \draw[style=help lines,step=1cm] (-3,-8) grid (3,8);    
    \draw[->] (-4,0) -- (4.5,0) node[anchor=west]{$\omega_1$-axis};
    \draw[->] (0,-8) -- (0, 8.5) node[anchor=south]{$\delta$-axis};
      \draw[-,thick] (-2,-8) -- (-2,8) ;
      \draw[-,thick] (2,-8) -- (2,8) ;
	\node at (2, -8) {$\bullet$};  
    \filldraw (2, -8) node[anchor=south west,yshift=-0.1cm] {$(s_0s_1)^4\Lambda$};
	\node at (2, -6) {$\bullet$};  
    \filldraw (2, -6) node[anchor=south west,yshift=-0.1cm] {$(s_0s_1)^3\Lambda$};
	\node at (2, -4) {$\bullet$};  
    \filldraw (2, -4) node[anchor=south west,yshift=-0.1cm] {$(s_0s_1)^2\Lambda$};
	\node at (2, -2) {$\bullet$};  
    \filldraw (2, -2) node[anchor=south west,yshift=-0.1cm] {$s_0s_1\Lambda$};
	\node at (2, 0) {$\bullet$};  
    \filldraw (2, 0) node[anchor=south west,yshift=0.05cm] {$\Lambda=2\omega_1+0\Lambda_0$};
	\node at (2, 2) {$\bullet$};  
    \filldraw (2, 2) node[anchor=south west,yshift=-0.1cm] {$s_1s_0\Lambda$};
	\node at (2, 4) {$\bullet$};  
    \filldraw (2, 4) node[anchor=south west,yshift=-0.1cm] {$(s_1s_0)^2\Lambda$};
	\node at (2, 6) {$\bullet$};  
    \filldraw (2, 6) node[anchor=south west,yshift=-0.1cm] {$(s_1s_0)^3\Lambda$};
	\node at (2, 8) {$\bullet$};  
    \filldraw (2, 8) node[anchor=south west,yshift=-0.1cm] {$(s_1s_0)^4\Lambda$};
	\node at (-2, -8) {$\bullet$};  
    \filldraw (-2, -8) node[anchor=south east,yshift=-0.1cm] {$s_1(s_0s_1)^4\Lambda$};
	\node at (-2, -6) {$\bullet$};  
    \filldraw (-2, -6) node[anchor=south east,yshift=-0.1cm] {$s_1(s_0s_1)^3\Lambda$};
	\node at (-2, -4) {$\bullet$};  
    \filldraw (-2, -4) node[anchor=south east,yshift=-0.1cm] {$s_1(s_0s_1)^2\Lambda$};
	\node at (-2, -2) {$\bullet$};  
    \filldraw (-2, -2) node[anchor=south east,yshift=-0.1cm] {$s_1s_0s_1\Lambda$};
	\node at (-2, 0) {$\bullet$};  
    \filldraw (-2, 0) node[anchor=south east,yshift=-0.1cm] {$s_1\Lambda$};
	\node at (-2, 2) {$\bullet$};  
    \filldraw (-2, 2) node[anchor=south east,yshift=-0.1cm] {$s_0\Lambda$};
	\node at (-2, 4) {$\bullet$};  
    \filldraw (-2, 4) node[anchor=south east,yshift=-0.1cm] {$s_0s_1s_0\Lambda$};
	\node at (-2, 6) {$\bullet$};  
    \filldraw (-2, 6) node[anchor=south east,yshift=-0.1cm] {$s_0(s_1s_0)^2\Lambda$};
	\node at (-2, 8) {$\bullet$};  
    \filldraw (-2, 8) node[anchor=south east,yshift=-0.1cm] {$s_0(s_1s_0)^3\Lambda$};
	\foreach \y in {-8,...,8}{
	\node at (0, \y) {$\bullet$}; 
	}
	\foreach \y in {-8, -6, -4, -2, 0, 2, 4, 6, 8}{
	\draw[-latex,  bend right=30, thick, red] (2, \y) to node[above]{$\tilde{f}_1$} (0, \y);
	\draw[-latex,  bend right=30, thick, red] (0, \y) to node[above]{$\tilde{f}_1$} (-2, \y);
	}
		\foreach \y in {-8, -6, -4, -2, 0, 2, 4, 6}{
	\draw[-latex,  bend left=30, thick, blue] (0, \y+1) to node[above]{$\tilde{f}_0$} (2, \y);
	\draw[-latex,  bend left=30, thick, blue] (-2, \y+2) to node[above]{$\tilde{f}_0$} (0, \y+1);
	}
\end{tikzpicture}
\\
\hbox{Paths in $B(2\omega_1)_0$}
&&\hbox{The crystal graph of $B(2\omega_1)_0$}
\end{matrix}
$$

\newpage

\noindent
\textbf{PLATE D: Pictures and characters of $B(\omega_1+\omega_2)$ for $\widehat{\fsl}_3$.}
The colour red indicates change in the $\delta$-coordinate.
$$
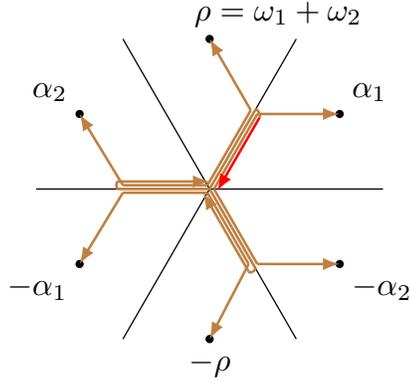
\begin{matrix}
\begin{matrix}
\resizebox{6cm}{!}{
\begin{tikzpicture}[scale=0.95,every node/.style={minimum size=1cm}]
\tikzstyle{every node}=[font=\small]
	 \coordinate (omega2) at (-0.5,0.86603);
	 \coordinate (2omega2) at  (-1,2*0.86603);
	 \coordinate (2omega2+omega1) at (-0.5, 3*0.86603);
	 \coordinate (omega2+2omega1) at (0.5,3*0.86603);
	 \coordinate (2omega1) at (1,2*0.86603);
	 \coordinate (omega1) at (0.5,0.86603);
	 \coordinate (-omega1) at (-0.5,-0.86603);
	 \coordinate (-omega1+omega2) at (-1,0);
	 \coordinate (-omega2+omega1) at (1,0);
	 \coordinate (-omega2) at  (0.5,-0.86603);
	 \coordinate (-2omega1-omega2) at (-0.5,-3*0.86603);
	 \coordinate (-2omega1) at (-1, -2*0.86603);
	 \coordinate (-2omega2) at (1,-2*0.86603);
	 \coordinate (-2omega2-omega1) at (0.5,-3*0.86603);
	 \coordinate (rho) at ($(omega1)+(omega2)$);
	 \coordinate (-rho) at ($-1*(omega1)-(omega2)$);
	 \coordinate (alpha1) at ($2*(omega1)-(omega2)$);
	 \coordinate (-alpha1) at ($-2*(omega1)+(omega2)$);
	 \coordinate (alpha2) at ($2*(omega2)-(omega1)$);
	 \coordinate (-alpha2) at ($-2*(omega2)+(omega1)$);	
	 
	\draw[thin] ($2*(-omega1)$) to  ($2*(omega1)$);
	\draw[thin] ($2*(-omega2)$) to ($2*(omega2)$);
	\draw[thin] ($2*(-omega1+omega2)$) to ($2*(-omega2+omega1)$);
	
	\draw[fill=black]
	(rho) circle (.04) node[anchor=south west, xshift=-0.3cm] {$\rho=\omega_1+\omega_2$}
	(alpha1) circle (.04) node[anchor=south west] {$\alpha_1$}
	(alpha2) circle (.04) node[anchor=south east] {$\alpha_2$}
	(-alpha1) circle (.04) node[anchor=north east] {$-\alpha_1$}
	(-alpha2) circle (.04) node[anchor=north west] {$-\alpha_2$}
	($-1*(omega1)-(omega2)$) circle (.04) node[anchor=north] {$-\rho$};

	\draw[thick, brown] ($0.05*(omega2)$) to ($(omega1)+0.05*(omega2)$);
	\draw[-latex, thick, brown] ($(omega1)+0.05*(omega2)$) to (rho);
	
	\draw[thick, brown] ($0.05*(-omega2)+0.05*(omega1)$) to ($(omega1)+0.05*(omega1)+0.05*(-omega2)$);
	\draw[-latex, thick, brown] ($(omega1)+0.05*(omega1)+0.05*(-omega2)$) to (alpha1);
	
	\draw[thick, brown] ($0.05*(omega2)$) to ($(omega2)-(omega1)+0.025*(rho)$);
	\draw[-latex, thick, brown] ($(omega2)-(omega1)+0.025*(rho)$) to (alpha2);
	
	\draw[thick, brown] ($0.05*(-omega1)$) to ($(omega2)-(omega1)-0.025*(rho)$);
	\draw[-latex, thick, brown] ($(omega2)-(omega1)-0.025*(rho)$)to (-alpha1);
	
	\draw[thick, brown] ($0.05*(-omega1)$) to ($(-omega2)-0.05*(omega1)-0.05*(-omega2)$);
	\draw[-latex, thick, brown] ($(-omega2)-0.05*(omega1)-0.05*(-omega2)$) to (-rho);
	
	\draw[thick, brown]  ($0.05*(-omega2)+0.05*(omega1)$) to ($(-omega2)+0.05*(omega1)+0.05*(-omega2)$);
	\draw[-latex, thick, brown] ($(-omega2)+0.05*(omega1)+0.05*(-omega2)$) to (-alpha2);
	
	\draw[thick, brown] (0,0) to ($1.05*(omega1)$);
	\draw[thick, out=60, in=60, brown] ($1.05*(omega1)$) to ($1.05*(omega1)-0.1*(omega2)$);
	\draw[-latex, thick, red] ($1.05*(omega1)-0.1*(omega2)$) to ($0.075*(omega1)-0.1*(omega2)$);
	
	\draw[thick, brown] (0,0) to ($1.05*(-omega2)$);
	\draw[thick, out=-60, in=-60, brown] ($1.05*(-omega2)$) to ($1.05*(-omega2)+0.05*(-alpha1)$);
	\draw[-latex, thick,brown] ($1.05*(-omega2)+0.05*(-alpha1)$) to ($0.05*(-alpha1)$);
	
	\draw[thick, brown] (0,0) to ($1.05*(omega2)-1.05*(omega1)$);
	\draw[thick, out=-180, in=-180, brown] ($1.05*(omega2)-1.05*(omega1)$) to ($1.05*(omega2)-1.05*(omega1)+0.05*(rho)$);
	\draw[-latex, thick,brown] ($1.05*(omega2)-1.05*(omega1)+0.05*(rho)$) to ($0.05*(rho)$);

\end{tikzpicture}
}
\\
B^{\mathrm{fin}}(\omega_1+\omega_2)\subseteq B(\omega_1)\otimes B(\omega_2) \\
\\
\mathrm{gchar}(B^{\mathrm{fin}}(\omega_1+\omega_2)) 
= X^\rho+X^{\alpha_1}+X^{\alpha_2}+X^{-\alpha_1} \\
\phantom{\mathrm{gchar}(B^{\mathrm{fin}}(\omega_1+\omega_2)) }\ \ 
+X^{-\alpha_2} +X^{-\rho} + 2 + q^{-1}
\end{matrix}
\,\ &
\begin{matrix}
\resizebox{6cm}{!}{
\begin{tikzpicture}[scale=1.1,every node/.style={minimum size=1cm}]
		\begin{scope}[
		yshift=0,
		every node/.append style={yslant=\yslant,xslant=\xslant},
		yslant=\yslant,xslant=\xslant
		]
		\draw [fill=black]
		(3,0) circle (.1) node[anchor=north west] {$\rho+\delta$}
		(1.5,2.5) circle (.1) node[anchor=south east] {$\alpha_2+\delta$}
		(1.5,-2.5) circle (.1) node[anchor=north west] {$\alpha_1+\delta$}
		(-1.5,2.5) circle (.1) node[anchor=south east] {$-\alpha_1+\delta$}
		(-1.5,-2.5) circle (.1) node[anchor=north west] {$-\alpha_2+\delta$}
		(-3,0) circle (.1) node[anchor=south east, xshift=-2] {$-\rho+\delta$};
		
		\draw[-latex, out=90, in=-20, thick, blue]
		(3,0) to node[right]{$\tilde{f}_1$}  (1.6,2.5);
		\draw[-latex, out=-100, in=40, thick, blue]
		(3,0) to node[right]{$\tilde{f}_2$}  (1.6,-2.5);
		\draw[-latex, out=-130, in=45, thick, blue]
		(1.5,2.5) to node[right]{$\tilde{f}_2$}  (0,0);
		\draw[-latex, out=130, in=-60, thick, blue]
		 (1.5,-2.5) to node[left]{$\tilde{f}_1$}  (0,0);
		\draw[-latex, out=-150, in=45, thick, blue]
		(0,0) to node[left, yshift=-2, xshift=4]{$\tilde{f}_2$}  (-1.5,-2.4);
		\draw[-latex, out=150, in=-60, thick, blue]
		 (0,0) to node[left]{$\tilde{f}_1$}  (-1.5,2.4);
		 \draw[-latex, out=-150, in=100, thick, blue]
		(-1.6,2.5) to node[left]{$\tilde{f}_2$}  (-3,0.1);
		\draw[-latex, out=150, in=-60, thick, blue]
		 (-1.6,-2.5) to node[left]{$\tilde{f}_1$}  (-3,-0.1);
		 
		  \coordinate (delta) at (0,0);
		 \coordinate (-pho+delta) at (-3,0);
		 \coordinate (s1s2pho+delta) at (-1.5,2.5);
		 \coordinate (s2s1pho+delta) at (-1.5,-2.5);
		 \coordinate (s1pho+delta) at (1.5,2.5);
		 \coordinate (s2pho+delta) at (1.5,-2.5);
		 \coordinate (pho+delta) at (3,0);
		\end{scope}
		
		\begin{scope}[
		yshift=-175,
		every node/.append style={yslant=\yslant,xslant=\xslant},
		yslant=\yslant,xslant=\xslant
		]
		\draw [fill=black]
		(3,0) circle (.1) node[anchor=north west] {$\rho=\omega_1+\omega_2$}
		(1.5,2.5) circle (.1) node[anchor=south east, xshift=7] {$\alpha_2$}
		(1.5,-2.5) circle (.1) node[anchor=north west] {$\alpha_1$}
		(-1.5,2.5) circle (.1) node[anchor=south east] {$-\alpha_1$}
		(-1.5,-2.5) circle (.1) node[anchor=north west, xshift=-8] {$-\alpha_2$}
		(-3,0) circle (.1) node[anchor=south east] {$-\rho$};
		\draw[-latex, out=90, in=-20, thick, blue]
		(3,0) to node[right, xshift=4, yshift=-6]{$\tilde{f}_1$}  (1.6,2.5);
		\draw[-latex, out=-100, in=40, thick, blue]
		(3,0) to node[right]{$\tilde{f}_2$}  (1.6,-2.5);
		\draw[-latex, out=-130, in=45, thick, blue]
		(1.5,2.5) to node[right, xshift=-4]{$\tilde{f}_2$}  (0,0);
		\draw[-latex, out=130, in=-60, thick, blue]
		 (1.5,-2.5) to node[left]{$\tilde{f}_1$}  (0,0);
		\draw[-latex, out=-150, in=45, thick, blue]
		(0,0) to node[left, yshift=-2, xshift=4]{$\tilde{f}_2$}  (-1.5,-2.4);
		\draw[-latex, out=150, in=-60, thick, blue]
		 (0,0) to node[left]{$\tilde{f}_1$}  (-1.5,2.4);
		\draw[-latex, out=-150, in=100, thick, blue]
		(-1.6,2.5) to node[left]{$\tilde{f}_2$}  (-3,0.1);
		\draw[-latex, out=150, in=-60, thick, blue]
		 (-1.6,-2.5) to node[left]{$\tilde{f}_1$}  (-3,-0.1);
		  \coordinate (0) at (0,0);
		 \coordinate (-pho) at (-3,0);
		 \coordinate (s1s2pho) at (-1.5,2.4);
		 \coordinate (s2s1pho) at (-1.5,-2.6);
		 \coordinate (s1pho) at (1.5,2.5);
		 \coordinate (s2pho) at (1.5,-2.5);
		 \coordinate (pho) at (3,0);
		\end{scope}

		\begin{scope}[
		yshift=-350,
		every node/.append style={yslant=\yslant,xslant=\xslant},
		yslant=\yslant,xslant=\xslant
		]
		\draw [fill=black]
		(3,0) circle (.1) node[anchor=north west] {$\rho-\delta$}
		(1.5,2.5) circle (.1) node[anchor=south east, xshift=13] {$\alpha_2-\delta$}
		(1.5,-2.5) circle (.1) node[anchor=north west] {$\alpha_1-\delta$}
		(-1.5,2.5) circle (.1) node[anchor=south east] {$-\alpha_1-\delta$}
		(-1.5,-2.5) circle (.1) node[anchor=north west] {$-\alpha_2-\delta$}
		(-3,0) circle (.1) node[anchor=south east] {$-\rho-\delta$};
		\draw[-latex, out=90, in=-20, thick, blue]
		(3,0) to node[right, xshift=4]{$\tilde{f}_1$}  (1.6,2.5);
		\draw[-latex, out=-100, in=40, thick, blue]
		(3,0) to node[right]{$\tilde{f}_2$}  (1.6,-2.5);
		\draw[-latex, out=-130, in=45, thick, blue]
		(1.5,2.5) to node[right, xshift=-4]{$\tilde{f}_2$}  (0,0);
		\draw[-latex, out=130, in=-60, thick, blue]
		 (1.5,-2.5) to node[left]{$\tilde{f}_1$}  (0,0);
		\draw[-latex, out=-150, in=45, thick, blue]
		(0,0) to node[left, yshift=-2, xshift=4]{$\tilde{f}_2$}  (-1.5,-2.4);
		\draw[-latex, out=150, in=-60, thick, blue]
		 (0,0) to node[left]{$\tilde{f}_1$}  (-1.5,2.4);
		\draw[-latex, out=-150, in=100, thick, blue]
		(-1.6,2.5) to node[left]{$\tilde{f}_2$}  (-3,0.1);
		\draw[-latex, out=150, in=-60, thick, blue]
		 (-1.6,-2.5) to node[left]{$\tilde{f}_1$}  (-3,-0.1);
		 \coordinate (-delta) at (0,0);
		 \coordinate (-pho-delta) at (-3,0);
		 \coordinate (s1s2pho-delta) at (-1.5,2.5);
		 \coordinate (s2s1pho-delta) at (-1.5,-2.5);
		 \coordinate (s1pho-delta) at (1.5,2.5);
		 \coordinate (s2pho-delta) at (1.5,-2.5);
		 \coordinate (pho-delta) at (3,0);
		\end{scope}
		
		\draw[-latex, out=-90, in=-150, thick, red] (-pho+delta) to node[above]{$\tilde{f}_0$} (0);
		\draw[-latex, out=0, in=-150, thick, red] (0) to node[right]{$\tilde{f}_0$} (pho-delta);
		\draw[-latex, out=-60, in=150, thick, red] (s2s1pho+delta) to node[left]{$\tilde{f}_0$} (s2pho);
		\draw[-latex, out=-60, in=-160, thick, red] (s2s1pho) to node[right]{$\tilde{f}_0$} (s2pho-delta);
		\draw[-latex, out=-60, in=-150, thick, red] (s1s2pho+delta) to node[below, xshift=-3]{$\tilde{f}_0$} (s1pho);
		\draw[-latex, out=-60, in=-130, thick, red] (s1s2pho) to node[below, xshift=-3]{$\tilde{f}_0$} (s1pho-delta);
		\draw[-latex, out=-60, in=-150, thick, red] (pho)(delta) to node[right, xshift=-3]{$\tilde{f}_0$} (pho);
		\draw[-latex, out=-90, in=-150, thick, red] (-pho) to node[right, xshift=-3]{$\tilde{f}_0$} (-delta);
		
		\draw[thick, black] (-pho+delta) to (-pho-delta);
		\draw[thick, black] (s2s1pho+delta) to (s2s1pho-delta);
		\draw[thick, black] (s2pho+delta) to (s2pho-delta);
		\draw[thick, black, dashed] (pho+delta) to (pho-delta); 
		\draw[thick, black] (s1s2pho+delta) to (s1s2pho-delta); 
		\draw[thick, black, dashed] (s1pho+delta) to (s1pho-delta); 
		\draw[thick, black] (-pho+delta) to (s1s2pho+delta);
		\draw[thick, black] (s1s2pho+delta) to (s1pho+delta);
		\draw[thick, black] (s1pho+delta) to (pho+delta);
		\draw[thick, black] (pho+delta) to (s2pho+delta);
		\draw[thick, black] (s2pho+delta) to (s2s1pho+delta);
		\draw[thick, black] (s2s1pho+delta) to (-pho+delta);
		\draw[thick, black] (-pho) to (s1s2pho);
		\draw[thick, black, dashed] (s1s2pho) to (s1pho);
		\draw[thick, black, dashed] (s1pho) to (pho);
		\draw[thick, black, dashed] (pho) to (s2pho);
		\draw[thick, black] (s2pho) to (s2s1pho);
		\draw[thick, black] (s2s1pho) to (-pho);
		\draw[thick, black] (-pho-delta) to (s1s2pho-delta);
		\draw[thick, black, dashed] (s1s2pho-delta) to (s1pho-delta);
		\draw[thick, black, dashed] (s1pho-delta) to (pho-delta);
		\draw[thick, black, dashed] (pho-delta) to (s2pho-delta);
		\draw[thick, black] (s2pho-delta) to (s2s1pho-delta);
		\draw[thick, black] (s2s1pho-delta) to (-pho-delta);
\end{tikzpicture}
}
\\
\hbox{$B(\omega_1+\omega_2)$ crystal graph}
\end{matrix}
\end{matrix}
$$
At $t=0$ and $t=\infty$ the 
normalized nonsymmetric Macdonald polynomials $\tilde E_{s_1s_2s_1\rho}(q,t)$ are
\begin{align*}
\tilde E_{s_1s_2s_1\rho}(q,0) 
&= X^{s_1s_2s_1\rho}+X^{s_1s_2\rho}
+X^{s_2s_1\rho}+X^{s_2\rho}+X^{s_1\rho}+X^\rho+2+q, \\
\tilde E_{s_1s_2s_1\rho}(q, \infty)
&= X^{s_1s_2s_1\rho} +q^{-1}(X^{s_1s_2\rho}+X^{s_2s_1\rho})
+q^{-2} (X^{s_1\rho} + X^{s_2\rho} + X^\rho)
+ 2q^{-2} +q^{-1}.
\end{align*}
Letting $q=e^{-\delta}$ (as in \cite[(12.1.9)]{Kac}),
the Demazure module $L(\omega_1+\omega_2)_{\le s_1s_2s_1}$ has character
\begin{align*}
\mathrm{char}(L(\omega_1+\omega_2)_{\le s_1s_2s_1})
&=\frac{1}{(1-q^{-1})^2} \tilde E_{s_1s_2s_1\rho}(q^{-1},0)
\quad\hbox{and}
\\
\mathrm{char}(L(\omega_1+\omega_2))
&= 0_q 0_q \tilde E_{s_1s_2s_1\rho}(q^{-1},0)
= 0_q 0_q\tilde E_{s_1s_2s_1\rho}(q^{-1},\infty),
\end{align*}
where 
$0_q = \cdots + q^{-3}+q^{-2}+q^{-1}+1+q+q^2+\cdots$ as in Remark \ref{RGchar}.
%

\begin{remark}
The expansion
\begin{align*}
&\frac{1}{(1-q^{-1})^2}=(1+q^{-1}+q^{-2}+\cdots) (1+q^{-1}+q^{-2}+\cdots) =1+2q^{-1}+3q^{-2}+4q^{-3}+5q^{-4}+\cdots 
\end{align*}
show that the sizes of the weight spaces of $L(\omega_1+\omega_2)_{\le s_1s_2s_1}$ 
are growing as $\delta$ increases.
Similarly, in the character formula of $L(\omega_1+\omega_2)$, the factor $0_q 0_q$
has coefficient of $q^n$ equal to
$\Card\big(\{(k_1, k_2)\in \ZZ^2\ |\  k_1+k_2=n\}\big)=\infty$. 
This shows that \emph{every weight space of the extremal weight module 
$L(\omega_1+\omega_2)$ is infinite dimensional}. 
\end{remark}

\newpage

\section{Affine Weyl groups, braid groups and Macdonald polynomials}

\subsection{The affine Lie algebra $\fg$}\label{affLiealgebra}

Let $\mathring{\fg}$ be a finite dimensional complex semisimple Lie algebra and fix a Cartan subalgebra $\fa\subseteq \mathring{\fg}$
and a symmetric, ad-invariant, nondegenerate, bilinear form 
$\langle ,\rangle\colon \mathring{\fg}\times \mathring{\fg}\to \CC$.
The \textit{affine Kac-Moody algebra} is 
$$
\fg=\left( \bigoplus_{k\in \ZZ} \mathring{\fg}\epsilon^k \right) \oplus \CC K \oplus \CC d,
\qquad\hbox{with bracket given by}\quad
[K,x\epsilon^k]=0,  \quad [K,d]=0, 
$$
\begin{equation}
[d, x\epsilon^k]=k x\epsilon^k,
\quad\hbox{and}\quad
[x\epsilon^k,y\epsilon^{\ell}] = [x,y]\epsilon^{k+\ell} + k \delta_{k,-\ell} \langle x,y\rangle K,
\end{equation}
for $x,y\in  \mathring{\fg}$ and $k,\ell \in \ZZ$ (see \cite[(7.2.2)]{Kac}).  
Let $\theta$ be the highest root of $\mathring{\fg}$ (the highest weight of the adjoint representation) and define
$$
e_0 = f_{\theta}\epsilon, 
\qquad f_0 = e_{\theta}\epsilon^{-1},
\quad\hbox{and}\quad h_0 = [e_0,f_0]
= -h_\theta + K.
$$
The miracle is that $\fg$ is a Kac-Moody Lie algebra with Chevalley generators 
\begin{equation}
e_0,\ldots ,e_n,\  h_0,\ldots h_n,\ d, \  f_0,\ldots f_n,
\qquad\hbox{which satisfy Serre relations.}
\label{KMpres}
\end{equation}
Because of \eqref{KMpres}, $\fg$ has a corresponding quantum enveloping algebra 
$\mathbf{U}= U_q\fg$.

The Cartan subalgebra of $\fg$ is 
$$
\fh=\fa \oplus \CC K\oplus \CC d, \,\ \,\ \text{where $\fa\subseteq  \mathring{\fg}$ is the Cartan subalgebra of $ \mathring{\fg}$. }
$$
Let $\mathring{R}^+$ be the set of positive roots for $\mathring{\fg}$.
For $\alpha\in\mathring{R}^+$, $k\in \ZZ$, $\ell\in \ZZ_{\ne 0}$ and $i\in \{1, \ldots, n\}$, let
$$x_{\alpha+k\delta} = e_\alpha \epsilon^k, \qquad 
x_{-\alpha+k\delta} = f_\alpha \epsilon^k, \qquad 
h_{i,\ell} = h_i\epsilon^\ell.$$
The \emph{homogeneous Heisenberg subalgebra} (see \cite[\S8.4 and \S14.8]{Kac}) is
\begin{equation}
\CC K\oplus \fa[\epsilon, \epsilon^{-1}]
\quad\hbox{with\quad}
[h_i\epsilon^k, h_j\epsilon^\ell] = k \delta_{k,-\ell} \frac{2}{\langle \alpha_i, \alpha_i\rangle}
\alpha_i(h_j) K,
\end{equation}
and $\fa[\epsilon]$ is a commutative Lie algebra with basis
$\{ h_i\epsilon^k\ |\ i\in \{1, \ldots, n\},\ k\in \ZZ_{\ge0} \}$.

\subsection{The affine Weyl group $W^{\mathrm{ad}}$ and its action on $\fh^*$ and $\fh$}\label{affineWeyl}

Let $\delta, \omega_1,\ldots, \omega_n, \Lambda_0$ be the basis in $\fh^*$ which is
the dual basis to the basis $d, h_1, \ldots, h_n, K$ of $\fh$.
The affine Weyl group $W^{\mathrm{ad}}$ 
is the subgroup of $GL(\fh^*)$ generated by the linear transformations
$s_0, s_1, \ldots, s_n$ which, in the basis $\delta, \omega_1, \ldots, \omega_n, \Lambda_0$,
are 
\begin{equation}
s_i = \begin{pmatrix}
1 &0 &\cdots &0 &0 &0 &\cdots &0 \\
0 &1 &\cdots &0 &-\alpha_i(h_1) &0 &\cdots &0 \\
0 &0 &\cdots &0 &-\alpha_i(h_2) &0 &\cdots &0 \\
\vdots &&&&\vdots &&&\vdots \\
0 &0 &\cdots &1 &-\alpha_i(h_{i-1}) &0 &\cdots &0 \\
0 &0 &\cdots &0 &-1 &0 &\cdots &0 \\
0 &0 &\cdots &0 &-\alpha_i(h_{i+1}) &1 &\cdots &0 \\
\vdots &&&&\vdots &&&\vdots \\
0 &0 &\cdots &0 &-\alpha_i(h_n) &0 &\cdots &0 \\
0 &0 &\cdots &0 &0 &0 &\cdots &1 \\
\end{pmatrix},
\qquad\hbox{for $i\in \{1, \ldots, n\}$}
\label{simatrix}
\end{equation}
and, writing $\theta = a_1\alpha_1+\cdots+a_n\alpha_n$ and
$h_\theta = [e_\theta,f_\theta] = a_1^\vee h_1+\cdots+a_n^\vee h_n$,
\begin{equation}
s_0 = 
\begin{pmatrix}
1 &a^\vee_1 &a^\vee_2 &\cdots &a^\vee_n &-1 \\
0 &1-a_1 a_1^\vee &-a_1 a_2^\vee &\cdots &-a_1 a_n^\vee &a_1 \\
0 &-a_2 a_1^\vee &1-a_2 a_2^\vee &\cdots &-a_2 a_n^\vee &a_2 \\
\vdots &&&\vdots &&\vdots \\
0 &-a_n a_1^\vee &-a_n a_2^\vee &\cdots &1-a_n a_n^\vee &a_n \\
0 &0 &0 &\cdots &0 &1
\end{pmatrix}.
\label{szeromatrix}
\end{equation}

Let $\fa_\RR^* = \hbox{$\RR$-span}\{ \alpha_1, \ldots, \alpha_n\}$.
An \emph{alcove} is a fundamental region for the action of $W^{\mathrm{ad}}$ on 
$(\RR\delta+\fa^*_\RR+\Lambda_0)/\RR\delta$.  
As explained (for example) in \cite{Ra06} and \cite{RY11}, there is a bijection 
\begin{equation}
\begin{matrix}
W^{\mathrm{ad}} &\longleftrightarrow &\{\hbox{alcoves}\} \\
1 &\longmapsto &\{ x+\Lambda_0\in \fa_\RR^*+\Lambda_0
\ |\ x(h_i) >0\ \hbox{for $i\in \{0, \ldots, n\}$}\}
\end{matrix}
\label{alcovedefn}
\end{equation}

Let $\fa^{\mathrm{ad}}_\ZZ = \hbox{$\ZZ$-span}\{h_1, \ldots, h_n\}$.
The \emph{finite Weyl group} $W_{\mathrm{fin}}$ is generated by $s_1, \ldots, s_n$.
The translation presentation of the affine Weyl group is 
\begin{equation}
W^{\mathrm{ad}}  = \fa^{\mathrm{ad}}_\ZZ \rtimes W_{\mathrm{fin}}
= \{ t_{\mu^\vee}u  \ |\  \mu^\vee\in \fa^{\mathrm{ad}}_\ZZ, u\in W_{\mathrm{fin}}\}
\quad\hbox{with}\quad
\begin{array}{l}
t_{\mu^\vee}t_{\nu^\vee} = t_{\mu^\vee+\nu^\vee}\ \hbox{and}\\
ut_{\mu^\vee}  = t_{u\mu^\vee} u,
\end{array}
\label{affWeylrelations}
\end{equation}
for $\mu^\vee, \nu^\vee\in \fa^{\mathrm{ad}}_\ZZ$ and $u\in W_{\mathrm{fin}}$.

Let $\alpha_i^\vee$ be the image of $h_i$ under the isomorphism $\fa\stackrel{\sim}{\to} \fa^*$ coming
from the nondegenerate bilinear form on $\fa$ which is the restriction of the nondegerate
bilinear form on $\mathring{\fg}$.
In matrix form with respect to the basis $\delta, \omega_1, \ldots, \omega_n, \Lambda_0$
of $\fh^*$ the action of $W^{\mathrm{ad}}$ on $\fh^*$ is given by 
\begin{equation}
t_{\mu^\vee}  =\left(
\begin{array}{c|ccc|c}
1 &k_1 &\cdots &k_n 
&- \hbox{$\frac12$} \langle \mu^\vee, \mu^\vee\rangle \\
\hline
\vdots &\phantom{\ddots} & & &\mu^\vee_1 \\
0 &&1 & &\vdots \\
\vdots &\phantom{\ddots} & & &\mu^\vee_n \\
 \hline
0 &\cdots &0 &\cdots &1
\end{array}
\right),
\quad\hbox{for}\quad
\begin{array}{l}
\mu^\vee =k_1h_1+\cdots+k_nh_n \\
\phantom{\mu^\vee}= k_1\alpha_1^\vee+\cdots+k_n\alpha_n^\vee \\
\phantom{\mu^\vee}= \frac{k_1}{d_1}\alpha_1+\cdots +\frac{k_n}{d_n}\alpha_n \\
\phantom{\mu^\vee}=\mu_1^\vee\omega_1+\cdots+\mu_n^\vee\omega_n,
\end{array}
\label{tmumatrix}
\end{equation}
so that
$-\hbox{$\frac12$}\langle \mu^\vee, \mu^\vee\rangle
=-\hbox{$\frac12$}(\mu_1^\vee k_1+\cdots \mu_n^\vee k_n)
$.  (In \eqref{tmumatrix} $d_1, \dots, d_n$ are the minimal positive integers such that
the product of the diagonal matrix
$\diag(d_1,\ldots, d_n)$ with the Cartan matrix is symmetric, see \cite[(2.1.1)]{Kac}.)

The basis $\{d, h_1, \ldots, h_n, K\}$ of $\fh$ is the dual basis to the basis
$\{\delta, \omega_1, \ldots, \omega_n, \Lambda_0\}$ of $\fh^*$.  Using the 
$W^{\mathrm{ad}}$-action on $\fh$ given by
$$s_i \mu^\vee = \mu^\vee - \alpha_i(\mu^\vee)h_i,
\qquad\hbox{for $i\in \{0, \ldots, n\}$ and $\mu^\vee\in \fh$,}
$$
the matrices for the action of $s_0, s_1, \ldots, s_n$ on $\fh$, 
in the basis $\{d, h_1, \ldots, h_n, K\}$, are the \emph{transposes of the
matrices in \eqref{simatrix} and \eqref{szeromatrix}}.

\subsection{The positive level, negative level and level 0 Bruhat orders on $W^{\mathrm{ad}}$}
\label{Bruhatorders}

In the framework of Section \ref{affineflags}, where $G= \mathring{G}(\CC((\epsilon)))$ is the loop group,
the closure orders for the Schubert cells in 
the positive level (thin) affine flag variety $G/I^+$, 
the negative level (thick) affine flag variety $G/I^-$, 
and the level 0 (semi-infinite) affine flag variety $G/I^0$
give partial orders on the affine Weyl group $W^{\mathrm{ad}}$\,:
$$\overline{I^+wI^+} = \bigsqcup_{x\posleq\  w} I^+ x I^+,
\qquad
\overline{I^+wI^0} = \bigsqcup_{x\zeroleq\  w} I^+ x I^0,
\qquad
\overline{I^+wI^-} = \bigsqcup_{x\negleq\ w} I^+ x I^-.
$$
These orders can be described combinatorially as follows.

An element $w\in W^{\mathrm{ad}}$ is \emph{dominant} if
$$w(\rho+\Lambda_0) \in \hbox{$\RR_{\ge 0}$-span}\{ \omega_1, \ldots, \omega_n\}+\Lambda_0,
\qquad\hbox{where}\quad \rho = \omega_1+\cdots+\omega_n.$$
In the identification \eqref{alcovedefn}
of elements of $W^{\mathrm{ad}}$ with alcoves,
the dominant elements of $W^{\mathrm{ad}}$ are the alcoves in the dominant
Weyl chamber.

Let $x,w\in W^{\mathrm{ad}}$ and
let $w=s_{i_1}\cdots s_{i_\ell}$ be a reduced word for $w$ in the generators $s_0, \ldots, s_n$.  
The \emph{positive level Bruhat order on $W^{\mathrm{ad}}$} is defined by
$$x\posleq\ w \quad\hbox{if $x$ has a reduced word which is a subword of $w=s_{i_1}\cdots s_{i_\ell}$}
$$
The \emph{negative level Bruhat order} on $W^{\mathrm{ad}}$ is defined by\quad
$x\negleq\  w$ if $x\posgeq\ w$.

\smallskip\noindent
The \emph{level 0} Bruhat order on $W^{\mathrm{ad}}$ is determined by
\begin{enumerate}
\item[(a)] 
$\zeroleq\ $ for dominant elements:
If $x,w$ are dominant then $x\zeroleq\ w$ if and only if $x\posleq\  w$,\quad
\item[(b)] 
$\zeroleq\ $ translation invariance:
If $\mu^\vee\in \fa^{\mathrm{ad}}_\ZZ$ and $x,w\in W$ then $x\zeroleq\ w$ if and only if $xt_{\mu^\vee} \zeroleq\  wt_{\mu^\vee}$.
\end{enumerate}

\smallskip\noindent
The \emph{positive level length} is $\ell^+\colon W^{\mathrm{ad}}\to \ZZ_{\ge 0}$ given by
$\ell^+(w) = (\hbox{length of a reduced word for $w$}).$

\smallskip\noindent
The \emph{negative level length} is $\ell^-\colon W^{\mathrm{ad}}\to \ZZ_{\le 0}$ given by
$\ell^-(w) = -\ell^+(w).$

\smallskip\noindent
The \emph{level 0 length} is $\ell^0\colon W^{\mathrm{ad}}\to \ZZ$ given by
$$\hbox{$\ell^0(w) = \ell^+(w)$ if $w$ is dominant}
\qquad\hbox{and}\qquad
\ell^0(xt_{\mu^\vee}) -\ell^0(yt_{\mu^\vee}) 
= \ell^0(x)-\ell^0(y),
$$
for $x,y\in W^{\mathrm{ad}}$ and $\mu^\vee\in \fa_\ZZ^{\mathrm{ad}}$.
Using the formula for $\ell^+$ given in \cite[(2.8)]{Mac96}, gives a formula for $\ell^0$,
\begin{equation}
\ell^0(ut_{\mu^\vee}) = \ell^+(u) + 2\langle \rho, \mu^\vee\rangle,
\qquad\hbox{for $u\in W_{\mathrm{fin}}$, $\mu^\vee\in \fa^{\mathrm{ad}}_\ZZ$.}
\end{equation}
The length functions $\ell^+$, $\ell^-$ and $\ell^0$ return, 
respectively, the dimension, the codimension
and the relative dimension of Schubert cells in the positive level, negative level and
level 0 affine flag varieties.

\subsection{The affine braid groups $\cB^{\mathrm{sc}}$ and $\cB^{\mathrm{ad}}$}\label{bdgps}

Let $\omega_1^\vee, \ldots, \omega_n^\vee$ be the basis of $\fa$ which is dual to 
the basis $\alpha_1, \ldots, \alpha_n$ of $\fa^*$.  Let
$$\fa_\ZZ^{\mathrm{ad}} = \hbox{$\ZZ$-span}\{ h_1,\ldots, h_n\}
\quad\subseteq\quad
\fa_\ZZ^{\mathrm{sc}} = \hbox{$\ZZ$-span}\{ \omega_1^\vee, \ldots, \omega_n^\vee\}.
$$
The affine braid group $\cB^{\mathrm{ad}}$ (resp.\ $\cB^{\mathrm{sc}}$) 
is generated by $T_1, \ldots, T_n$ 
and $Y^{\lambda^\vee}$, $\lambda^\vee\in \fa^{\mathrm{ad}}_\ZZ$ (resp.
$\lambda^\vee\in \fa_\ZZ^{\mathrm{sc}}$), with
relations
$$Y^{\lambda^\vee}Y^{\sigma^\vee} = Y^{\lambda^\vee+\sigma^\vee} ,
\qquad
\underbrace{T_iT_j\cdots}_{m_{ij}\ \mathrm{factors}}
=\underbrace{T_jT_i\cdots}_{m_{ij}\ \mathrm{factors}},
\qquad
\begin{array}{cl}
T_i^{-1}Y^{\lambda^\vee} = Y^{s_i\lambda^\vee}T_i^{-1} ,  
&\hbox{if $\langle \lambda^\vee, \alpha_i\rangle=0$,} \\
T_i^{-1}Y^{\lambda^\vee}T_i^{-1} = Y^{s_i\lambda^\vee},  
&\hbox{if $\langle \lambda^\vee, \alpha_i\rangle=1$,}
\end{array}
$$
for $i\in \{1, \ldots, n\}$ and $\lambda^\vee\in \fa_\ZZ^{\mathrm{ad}}$ 
(resp.\ $\fa_\ZZ^{\mathrm{sc}}$) and $m_{ij} = \alpha_i(h_j)\alpha_j(h_i)$ for $i,j\in \{1, \ldots, n\}$ with $i\ne j$.

\subsection{Macdonald polynomials}\label{Macpolys}

Let 
$$\fh^*_\ZZ = \hbox{$\ZZ$-span}\{\delta, \omega_1, \ldots, \omega_n, \Lambda_0\}
\qquad\hbox{and}\qquad
\fa_\ZZ^* = \hbox{$\ZZ$-span}\{\omega_1, \ldots, \omega_n\}.
$$
The \emph{double affine Hecke algebra} $\tilde H$ is presented by generators
$T_0, \ldots, T_n$ and $X^\mu$, $\mu \in \fh_\ZZ^*$, with relations 
\begin{equation}
X^{\lambda}X^{\mu} = X^{\lambda+\mu},
\qquad
\underbrace{T_i T_j \cdots}_{m_{ij}\ \mathrm{factors}} 
= \underbrace{T_j T_i \cdots}_{m_{ij}\ \mathrm{factors}},
\qquad
T_i^2 = (t^{\frac12}-t^{-\frac12})T_i + 1,
\label{Heckereln}
\end{equation}
\begin{align*}
T_iX^\mu
= X^{s_i\mu}T_i + (t^{\frac12}-t^{-\frac12})\frac{X^\mu -X^{s_i\mu} }{1-X^{\alpha_i} }, 
\qquad
T_i^{-1}X^{\mu}
= X^{s_i\mu}T_i^{-1} - (t^{\frac12}-t^{-\frac12})\frac{X^{\mu}-X^{s_i\mu}}{1-X^{-\alpha_i}}.
\label{TinvpastX}
\end{align*}
for $i\in \{0, \ldots, n\}$ and $\mu\in \fh_\ZZ^*$.
For $w\in W^{\mathrm{ad}}$ put 
\begin{equation}
Y^w = \begin{cases}
Y^{ws_i} T^{-1}_i, &\hbox{if $w\zerol\ ws_i$,} \\
Y^{ws_i}T_i, &\hbox{if $w\zerog\ ws_i$,} 
\end{cases}
\qquad \hbox{and let}\qquad
Y^{\lambda^\vee} = Y^{t_{\lambda^\vee}}\ \hbox{for $\lambda^\vee\in \fa^{\mathrm{ad}}_\ZZ$.}
\label{Ydefn}
\end{equation}
Putting $q=X^\delta=Y^{-K}$
then, as an algebra over $\CC[q^{\pm1}, t^{\pm\frac12}]$,
$$\hbox{$\tilde H$ has basis}\qquad
\{ X^\mu T_u Y^{\lambda^\vee}\ |\ 
\mu\in \fa_\ZZ^*+\ZZ\Lambda_0, u\in W_{\mathrm{fin}}, \lambda^\vee\in \fa^{\mathrm{ad}}_\ZZ\},
\qquad\hbox{where}\quad
T_u = T_{i_1} \cdots T_{i_\ell},
$$
for a reduced word $u = s_{i_1}\cdots s_{i_\ell}$.
The \emph{affine Hecke algebra} is the subalgebra $H$ of $\tilde H$ with basis
$\{ T_uY^{\lambda^\vee}\ |\  u\in W_{\mathrm{fin}}, \lambda^\vee\in \fa^{\mathrm{ad}}_\ZZ\}$.
The \textit{polynomial representation of $\tilde{H}$} is 
\begin{equation}
\CC\left[X\right]
=
\Ind_{H}^{\tilde{H}}(\mathbf{1})
\qquad\hbox{with basis}\qquad
\{ X^\mu\mathbf{1}\ |\ \mu\in \fa_\ZZ^*\},
\label{polyaction}
\end{equation}
and\quad $Y^K \mathbf{1} = q^{-1}\mathbf{1},$
\quad $Y^{-\alpha_i^\vee}\mathbf{1} = t \mathbf{1}$,
\quad and\quad
$T_i\mathbf{1} = t^{\frac12} \mathbf{1}$ for $i\in \{1, \ldots, n\}$.

Let
\begin{equation}
T_0^\vee = 
Y^{\alpha_0^\vee}X^{-\Lambda_0}T_0^{-1} X^{\Lambda_0}
\quad\hbox{and}\quad
T_i^\vee = T_i\ \hbox{for $i\in \{1, \ldots, n\}$.}
\label{Ticheckdefn}
\end{equation}
The automorphism of $\tilde H$ given by conjugation by
$X^{-\Lambda_0}$ is the automorphism $\tau\colon \tilde H\to \tilde H$ of \cite[(2.8)]{Ch96}.  
Extend $\tilde H$ to allow rational functions in the $Y^{\lambda^\vee}$.
For each  $i\in \{0,1,\ldots ,n\}$, \textit{the intertwiner $\tau_i^{\vee}\in \widetilde{H}$} is 
\begin{equation}
\tau_i^\vee= T_i^\vee+ \frac{t^{-\frac{1}{2}}(1-t)}{1-Y^{-\alpha_i^{\vee}}}
=
(T^{\vee}_i)^{-1}+
\frac{t^{-\frac{1}{2}}(1-t)Y^{-\alpha_i^{\vee}}}{1-Y^{-\alpha_i^{\vee}}}
\quad \hbox{so that}\quad
Y^{\lambda^\vee} \tau_i^\vee = \tau_i^\vee Y^{s_i\lambda^\vee}.
\label{intertwinerdefn} 
\end{equation}
Let, for simplicity, $\mu\in \hbox{$\ZZ$-span}\{\alpha_1, \ldots, \alpha_n\}$ 
(the general case $\mu\in \fa_\ZZ^*$ requires consideraton of 
the group $\Omega^\vee$, the quotient of the 
weight lattice by the root lattice, and is treated in detail in \cite{RY11}).
The \emph{nonsymmetric Macdonald polynomial} $E_\mu=E_\mu(q,t)$ is 
\begin{equation}
E_\mu = E_{\mu}(q,t)
=\tau^\vee_{i_1}\ldots \tau^\vee_{i_{\ell}}\mathbf{1},
\qquad\hbox{where $m_\mu=s_{i_1}\ldots s_{i_{\ell}}$ is a reduced word}
\label{nsmacddefn}
\end{equation}
for the minimal length element in the coset $t_\mu W_{\mathrm{fin}}$.
The $E_\mu$ form a basis of $\CC[X]$ consisting of eigenvectors for
the $Y^{\lambda^\vee}$ (the Cherednik-Dunkl operators).

Fix a reduced word  $m_\mu=s_{i_1}\ldots s_{i_{\ell}}$ as in \eqref{nsmacddefn}.
Identifying the elements of $W^{\mathrm{ad}}$ with alcoves as in \eqref{alcovedefn},
an alcove walk of type $\vec m_\mu = (i_1, \ldots, i_\ell)$ beginning at $1$ (the fundamental
alcove) is a sequence of steps, of types $i_1, \ldots, i_\ell$,
where a step of type $j$ is (the signs - and + indicate that $zs_j \zerog\ z$)
$$
\begin{matrix}
\beginpicture
\setcoordinatesystem units <1cm,1cm>         
\setplotarea x from -1.5 to 1.5, y from -0.5 to 0.5  
\put{$\scriptstyle{zs_j}$} at 0.6 -0.25
\put{$\scriptstyle{z}$} at -0.6 -0.25
\put{$\scriptstyle{-}$}[b] at -0.4 0.25
\put{$\scriptstyle{+}$}[b] at 0.4 0.25
\plot  0 -0.4  0 0.5 /
\arrow <5pt> [.2,.67] from -0.5 0 to 0.5 0   %
\endpicture
&\beginpicture
\setcoordinatesystem units <1cm,1cm>         
\setplotarea x from -1.5 to 1.5, y from -0.5 to 0.5  
\put{$\scriptstyle{zs_j}$} at 0.6 -0.25
\put{$\scriptstyle{z}$} at -0.6 -0.25
\put{$\scriptstyle{-}$}[b] at -0.4 0.25
\put{$\scriptstyle{+}$}[b] at 0.4 0.25
\plot  0 -0.4  0 0.5 /
\arrow <5pt> [.2,.67] from 0.5 0 to -0.5 0   %
\endpicture
&
\beginpicture
\setcoordinatesystem units <1cm,1cm>         
\setplotarea x from -1.5 to 0.5, y from -0.5 to 0.5  
\put{$\scriptstyle{zs_j}$} at 0.6 -0.25
\put{$\scriptstyle{z}$} at -0.6 -0.25
\put{$\scriptstyle{-}$}[b] at -0.4 0.35
\put{$\scriptstyle{+}$}[b] at 0.4 0.35
\plot  0 -0.4  0 0.6 /
\plot 0.5 0  0.05 0 /
\arrow <5pt> [.2,.67] from 0.05 0.1 to 0.5 0.1   %
\plot 0.05 0 0.05 0.1 /
\endpicture
&
\beginpicture
\setcoordinatesystem units <1cm,1cm>         
\setplotarea x from -1.5 to 0.5, y from -0.5 to 0.5  
\put{$\scriptstyle{zs_j}$} at 0.6 -0.25
\put{$\scriptstyle{z}$} at -0.6 -0.25
\put{$\scriptstyle{-}$}[b] at -0.4 0.35
\put{$\scriptstyle{+}$}[b] at 0.4 0.35
\plot  0 -0.4  0 0.6 /
\plot -0.5 0  -0.05 0 /
\plot -0.05 0 -0.05 0.1 /
\arrow <5pt> [.2,.67] from -0.05 0.1 to -0.5 0.1   %
\endpicture
&\hbox{with $zs_j \zerog\ z$.}
\\
\hbox{positive $j$-crossing} &\hbox{negative $j$-crossing}
&\hbox{positive $j$-fold} &\hbox{negative $j$-fold}
\end{matrix}
$$
Let $\cB(1, \vec m_\mu)$ 
be the set of alcove walks of type $\vec m_\mu = (i_1,\ldots, i_\ell)$ beginning at $1$.
For a walk $p\in \cB(1, \vec m_\mu)$ let
$$\begin{array}{l}
f^+(p) = \{ k\in \{1, \ldots, \ell\} \ |\ \hbox{the $k$th step of $p$ is a positive fold}\}, \\
f^-(p) = \{ k\in \{1, \ldots, \ell\} \ |\ \hbox{the $k$th step of $p$ is a negative fold}\}, \\
f(p) = f^+(p)\cup f^-(p) = \{ k\in \{1, \ldots, \ell\} \ |\ \hbox{the $k$th step of $p$ is a fold}\}.
\end{array}
$$
For $p\in \cB(1, \vec m_\mu)$ let $\mathrm{end}(p)$ be the endpoint of $p$ (an element of $W^{\mathrm{ad}}$)
and define the \emph{weight} $\wt(p)$ and the \emph{final direction}
$\varphi(p)$ of $p$ by 
$$X^{\mathrm{end}(p)} = X^{\mathrm{wt}(p)} T^\vee_{\varphi(p)},
\qquad\hbox{with $\wt(p)\in \fa_\ZZ^*$ and $\varphi(p)\in W_{\mathrm{fin}}$.}
$$
Using \eqref{intertwinerdefn} and doing a left to right expansion of the terms of 
$\tau_{i_1}^\vee\cdots \tau_{i_\ell}^\vee\mathbf{1}$ produces the 
monomial expansion of $E_\mu$ as sum over alcove walks as given in the following theorem.
For simplicity we state the following theorem for $\mu\in 
\hbox{$\ZZ$-span}\{\alpha_1, \ldots, \alpha_n\}$.  It holds, after a small technical adjustment to the statement, for all $\mu\in \fa_\ZZ^*$, see \cite{RY11} for details.

\begin{thm} \emph{\cite[Theorem 3.1 and Remark 3.3]{RY11}}  
Let $\mu\in \hbox{$\ZZ$-span}\{\alpha_1, \ldots, \alpha_n\}$ 
and let 
$m_\mu$ be the minimal
length element in the coset $t_\mu W_{\mathrm{fin}}$.  Fix a reduced word
$\vec m_\mu = s_{i_1}\cdots s_{i_\ell}$, let 
$$\beta_1^\vee = s_{i_{\ell}}\cdots s_{i_2}\alpha^\vee_{i_1}, \quad
\quad\beta_2^\vee = s_{i_\ell}\cdots s_{i_3}\alpha^\vee_{i_2}, \quad \ldots,
\quad \beta_\ell^\vee = \alpha_{i_\ell}^\vee,$$
and let $\mathrm{sh}(\beta_k^\vee)$ and $\mathrm{ht}(\beta_k^\vee)$ be defined by
$Y^{\beta_k^\vee}\mathbf{1} 
= q^{\mathrm{sh}(\beta_k^\vee)} t^{\mathrm{ht}(\beta_k^\vee)}\mathbf{1}$
for $k\in \{1, \ldots, \ell\}$.  Then
$$E_\mu(q,t) = \sum_{p\in \cB(1, \vec m_\mu)}
X^{wt(p)} t^{\frac12(\ell(\varphi(p))}
\prod_{k\in f^+(p)} \frac{ t^{-\frac12}(1-t) }{ 1-q^{\mathrm{sh}(\beta_k^\vee)}t^{\mathrm{ht}(\beta_k^\vee)}}
\prod_{k\in f^-(p)} \frac{ t^{-\frac12} (1-t) q^{\mathrm{sh}(\beta_k^\vee)}t^{\mathrm{ht}(\beta_k^\vee)}} 
{ 1-q^{\mathrm{sh}(\beta_k^\vee)}t^{\mathrm{ht}(\beta_k^\vee)} }.
$$
\end{thm}

\subsection{Specializations of the normalized Macdonald polynomials $\tilde E_\mu(q,t)$}\label{Macspec}

If $m_\mu = t_\mu m$ with $m\in W_{\mathrm{fin}}$ then
$E_\mu(q,t)$ has top term $t^{\frac12\ell(m)}X^\mu$
(this term is the term corresponding to the unique alcove walk in $\cB(\vec m_\mu)$ with no folds).  The \emph{normalized nonsymmetric Macdonald polynomial} is
$$\tilde E_\mu(q,t) = t^{-\frac12\ell(m)}E_\mu(q,t)
\quad\hbox{so that $\tilde E_\mu(q,t)$ has top term $X^\mu$.}
$$

\smallskip\noindent
A path $p\in \cB(1, \vec m_\mu)$ is \emph{positively folded} if there are no negative folds, i.e.\ $\# f^-(p)=0$.

\smallskip\noindent
A path $p\in \cB(1, \vec m_\mu)$ is \emph{negatively folded} if there are no positive folds, i.e.\ $\# f^+(p)=0$.

\smallskip\noindent
A path $p\in \cB(1, \vec m_\mu)$ is \emph{positive semi-infinite} if
$\displaystyle{\ell(\varphi(p)) -\ell(m) - \#f(p)+2\sum_{k\in f^-(p)} ht(\beta_k^\vee)=0.}$
A path $p\in \cB(1, \vec m_\mu)$ is \emph{negative semi-infinite} if
$\displaystyle{\ell(m) - \ell(\varphi(p)) + \#f(p)+2\sum_{k\in f^+(p)} ht(\beta_k^\vee)=0.}$

\begin{prop} \label{Macdsplzn} Let $\mu\in \fa_\ZZ^*$.  The specializations
$q=0$, $t=0$, $q^{-1}=0$ and $t^{-1}=0$ are well defined and given, respectively, by
\begin{align*}
\tilde E_\mu(0,t) 
&= \sum_{p\in \cB(1, \vec m_\mu)\atop p\ \mathrm{pos\ folded}} 
t^{\frac12 (\ell(\varphi(p)) -\ell(m) - \# f(p))}(1-t)^{\#f(p)}
X^{\wt(p)}, \\
\tilde E_\mu(q,0) 
&= \sum_{p\in \cB(1, \vec m_\mu)\atop p\ \mathrm{neg\ semi-inf}} 
q^{\sum_{k\in f^-(p)} \mathrm{sh}(\beta_k^\vee)} X^{\wt(p)}, \\
\tilde E_\mu(\infty,t) 
&= \sum_{p\in \cB(1, \vec m_\mu)\atop p\ \mathrm{neg\ folded}} 
t^{-\frac12 (\ell(\varphi(p)) - \# f(p))}(1-t^{-1})^{\#f(p)}
X^{\wt(p)}, \\
\tilde E_\mu(q,\infty) 
&= \sum_{p\in \cB(1, \vec m_\mu)\atop p\ \mathrm{pos\ semi-inf}} 
q^{-\sum_{k\in f^+(p)} \mathrm{sh}(\beta_k^\vee)} X^{\wt(p)}.
\end{align*}
\end{prop}

For $i\in \{0,1,\ldots, n\}$ and $f\in \CC[\fh_\ZZ^*]$ define
\begin{equation}
\Delta _i f = \frac{f - s_if}{1-X^{-\alpha_i}} 
\qquad\hbox{and}\qquad
D_if = (1+s_i)\frac{1}{1-X^{-\alpha_i}} f
=\frac{f - X^{-\alpha_i}s_i f}{1-X^{-\alpha_i}}.
\label{Demazureopdefn}
\end{equation}
Equations \eqref{intertwinerdefn}, \eqref{Ticheckdefn} and the last relation in \eqref{Heckereln}  give that,
as operators on the polynomial representation,
\begin{align*}
t^{\frac12} \tau_i^\vee
&=  X^{-\Lambda_0}(D_i - t\Delta_i)X^{\Lambda_0}  
+ \frac{ (1-t) Y^{-\alpha_i^\vee}} {1-Y^{-\alpha_i^\vee}}
\quad\hbox{for $i\in \{1,\ldots, n\}$, \quad and} \\
t^{\frac12}\tau_0^\vee Y^{\alpha_0^\vee}
&= Y^{-\alpha_0^\vee} t^{\frac12}\tau_0^\vee
= X^{-\Lambda_0} (D_0 - t\Delta_0) X^{\Lambda_0}
+ \frac{(1-t)Y^{-\alpha_0^\vee}}{1-Y^{-\alpha_0^\vee}}.
\end{align*}
When applied in the formula of \eqref{nsmacddefn}, these formulas for (normalized) intertwiners are
specailizable at $t=0$, giving the following result (see the examples computed in Section \ref{Ionthmexamples}).

\begin{thm} \label{MacisDem}
\emph{(\cite[\S4.1]{Ion01})} 
Let $\mu\in \fa_\ZZ^*$. There are unique
$\nu\in \fa_\ZZ^*$ and $j\in \ZZ$ such that
$$\Lambda_0+\nu\in (\fh^*)_{\mathrm{int}}^+
\quad\hbox{and}\quad
-j\delta + \mu+\Lambda_0 \in W^{\mathrm{ad}}(\Lambda_0+\nu).$$  
Let $w\in W^{\mathrm{ad}}$ be minimal length such that $-j+\mu+\Lambda_0 = w(\Lambda_0+\nu)$
and let $w= s_{i_1}\cdots s_{i_\ell}$ be a reduced word.
Letting $D_0, \ldots, D_n$ be the Demazure operators given in \eqref{Demazureopdefn},
$$\tilde E_{\mu}(q,0) = q^{-j}X^{-\Lambda_0} D_{i_1}\cdots D_{i_\ell} X^{\nu} X^{\Lambda_0}
\mathbf{1}.$$
\end{thm}

\section{Quantum affine algebras $\mathbf{U}$ and integrable modules}

\subsection{The quantum affine algebra $\mathbf{U}$}

The \emph{quantum affine algebra} $\mathbf{U}$ is the $\CC(q)$-algebra
generated by
$$E_0, \ldots, E_n, \ F_0, \ldots, F_n,, \ K^{\pm1}_0, \ldots, K^{\pm1}_n,\ 
C^{\pm\frac12}, \ D^{\pm 1}$$
with Chevalley-Serre type relations corresponding to the affine Dynkin diagram. 
Following \cite[(11)]{Dr88}, \cite[\S39]{Lu93}, \cite[end of \S 1]{B94}, there is an action of the
affine braid group $\cB^{\mathrm{sc}}$ on $\mathbf{U}$ by automorphisms.

Let $\mathbf{U}^+$ be the subalgebra of $\mathbf{U}$ generated by $E_0, E_1, \ldots, E_n$.
As explained in \cite[Lemma 1.1 (iv)]{BCP98} and \cite[(3.1)]{BN02},
there is a doubly infinite ``longest element'' for the affine Weyl group with a favourite
reduced expression $w_\infty = \cdots s_{i_{-1}}s_{i_0} s_{i_1}\cdots$.  This reduced word
is used, with the braid group action, to define \emph{root vectors} in $\mathbf{U}^+$ by
\begin{equation}
E_{\beta_0} = E_{i_0},\quad\hbox{and}\quad
E_{\beta_{-k}} = T_{i_0}^{-1}T_{i_{-1}}^{-1}\cdots T^{-1}_{i_{-(k-1)}}E_{i_{-k}}   
\quad\hbox{and}\quad
E_{\beta_k} = T_{i_1}T_{i_2}\ldots T_{i_{k-1}}E_{i_k},
\label{rootvects}
\end{equation}
for $k\in \ZZ_{>0}$. 
For $i\in \{1, \ldots, n\}$ and $r,s\in \ZZ$ define the \emph{loop generators in $\mathbf{U}^+$}
\begin{equation}
\mathbf{x}^+_{i,r} = E_{\alpha_i+r\delta} = Y^{-r\omega^\vee_i} E_i
\quad\hbox{and}\quad
\mathbf{x}^-_{i,s} = E_{-\alpha_i+s\delta} = Y^{s\omega^\vee_i} F_i,
\end{equation}
where $Y^{r\omega_i^\vee}$ and $Y^{s\omega_i^\vee}$ are elements of the braid group $\cB^{\mathrm{sc}}$
as defined in Section \ref{bdgps}.
For $r\in \ZZ_{\ge 0}$ and $s\in \ZZ_{>0}$ 
these are special cases of the root vectors in \eqref{rootvects}.
For $i\in \{1, \ldots, n\}$ and $r\in \ZZ_{>0}$ define 
$\mathbf{q}^{(i)}_r$ by
\begin{equation}
\mathbf{q}^{(i)}_r 
= \mathbf{x}^-_{i,r}\mathbf{x}^+_{i,0} - q^{-2}\mathbf{x}^+_{i,0} \mathbf{x}^-_{i,r},
\quad\hbox{let}\quad
\mathbf{q}^{(i)}_+(z) = 1+(q-q^{-1}) \sum_{s\in \ZZ_{>0}} \mathbf{q}^{(i)}_s z^{s},
\label{qplusseries}
\end{equation}
and define $\mathbf{p}^{(i)}_r$ and $\mathbf{e}^{(i)}_r$  by
\begin{equation}
\mathbf{q}^{(i)}_+(z) = \mathrm{exp}\Big( \sum_{r\in \ZZ_{>0}} (q-q^{-1}) \mathbf{p}^{(i)}_{r} z^r \Big)
\quad\hbox{and}\quad
\mathrm{exp} \Big( \sum_{r\in \ZZ_{>0}} \frac{\mathbf{p}^{(i)}_r }  {[r]} z^r \Big)
=1+\sum_{k\in \ZZ_{>0}} \mathbf{e}^{(i)}_k z^k.
\label{thepsandtheqs}
\end{equation}
For a sequence of partitions $\vec\kappa = (\kappa^{(1)}, \ldots, \kappa^{(n)})$ define
\begin{equation}
\mathbf{s}_{\vec\kappa} = \mathbf{s}_{\kappa^{(1)}}\cdots \mathbf{s}_{\kappa^{(n)}},
\qquad\hbox{where}\quad
\mathbf{s}_{\kappa^{(i)}} = \det( \mathbf{e}^{(i)}_{(\kappa^{(i)})'_r-r+s})_{1\le r,s\le m_i},
\end{equation}
where $(\kappa^{(i)})'_r$ is the length of the $r$th column of $\kappa^{(i)}$ and 
$m_i=\ell(\kappa^{(i)})$ (see \cite[Ch.\ I\ (3.5)]{Mac}).
For a sequence
$\mathbf{c} = (\cdots, c_{-3}, c_{-2}, c_{-1}, \vec \kappa, c_1, c_2, c_3, \ldots)$
with $c_i\in \ZZ_{\ge 0}$ and all but a finite number of $c_i$ equal to $0$.  
The corresponding \emph{PBW-type element of $\mathbf{U}^+$} is
\begin{equation}
E_{\mathbf{c}} 
=  (E_{\beta_0}^{(c_0)}E_{\beta_{-1}}^{(c_{-1})} E_{\beta_{-2}}^{(c_{-2})}\cdots)
\,(\mathbf{s}_{\kappa^{(1)}}\cdots \mathbf{s}_{\kappa^{(n)} })\,
( \cdots E_{\beta_{3}}^{(c_{3})}E_{\beta_{2}}^{(c_{2})}E_{\beta_{1}}^{(c_{1})} ).
\label{PBWbasis}
\end{equation}

The \emph{Cartan involution} is the $\CC$-linear 
anti-automorphism $\Omega\colon \mathbf{U}\to \mathbf{U}$ given by
$$\Omega(E_i) = F_i, \quad
\Omega(F_i) = E_i, \quad
\Omega(K_i) = K_i^{-1}, \quad
\Omega(D) = D^{-1}, \quad \Omega(q) = q^{-1},
$$
and
$\mathbf{U}^- = \Omega(\mathbf{U}^+)$.
Putting 
$\mathbf{p}^{(i)}_{-r} = \Omega(\mathbf{p}^{(i)}_{r})$, 
then (see \cite[Th.\ 2 (6)]{Dr88} and \cite[Th.\ 4.7 (2)]{B94})
\begin{equation}
[\mathbf{p}^{(i)}_{k}, \mathbf{p}^{(i)}_{l}] = \delta_{k,-l} \frac{1}{k} 
\Big(\frac{q^{k \alpha_i(h_j)} - q^{-k\alpha_i(h_j)} }{q-q^{-1} }\Big)
\frac{C^k-C^{-k}}{q-q^{-1}}.
\label{qheisenbergrelation}
\end{equation}
Define $q^{(i)}_{-s} = \Omega(q^{(i)}_s)$ and 
\begin{equation}
\mathbf{q}^{(i)}_-(z^{-1}) = 1+(q-q^{-1}) \sum_{s\in \ZZ_{>0}} \mathbf{q}^{(i)}_{-s} z^{-s}
\label{qminusseries}
\end{equation}

The \emph{Heisenberg subalgebra}
$\mathbf{H}$ is the subalgebra of $\mathbf{U}$ generated by
$\{ \mathbf{p}^{(i)}_{k}\ |\ i\in \{1, \ldots, n\}, \  k\in \ZZ_{\ne 0}\}.$
In $\mathbf{H}\cap \mathbf{U}^+$
the $\mathbf{p}^{(i)}_r$ ($r\in \ZZ_{> 0}$) are the \emph{power sums},
the $\mathbf{q}^{(i)}_r$ ($r\in \ZZ_{> 0}$) are the \emph{Hall-Littlewoods} and the
$\mathbf{e}^{(i)}_r$ ($r\in \ZZ_{> 0}$) are the \emph{elementary symmetric functions} and
the $\mathbf{s}_{\vec\kappa}$ are the \emph{Schur functions}.

\subsection{Integrable $\mathbf{U}$-modules}\label{intmodules}

As in \eqref{szeromatrix}, let $h_\theta = a^\vee_1h_1+\cdots + a^\vee_n h_n$ 
be the highest root of $\mathring{\fg}$ and let 
$$\Lambda_i = \omega_i + a^\vee_i \Lambda_0,
\quad\hbox{for $i\in \{1, \ldots, n\}$,}
$$
so that $\{\delta, \Lambda_1, \ldots, \Lambda_n, \Lambda_0\}$ is the 
dual basis in $\fh^*$ to the basis $\{d, h_1, \ldots, h_n, h_0\}$ of $\fh$.
Let 
$$\fh^*_\ZZ = \{ \Lambda\in \fh^*\ |\ \hbox{$\langle \Lambda, \alpha_i^\vee\rangle \in \ZZ$ for 
$i\in \{0,1,\ldots, n\}$}\} 
= \CC\delta + \hbox{$\ZZ$-span}\{ \Lambda_0, \ldots, \Lambda_n\}.$$
A set of representatives for the $W^{\mathrm{ad}}$-orbits on $\fh_{\ZZ}^{\ast}$ is 
\begin{equation}
(\fh^{\ast})_{\mathrm{int}}=(\fh^{\ast})_{\mathrm{int}}^+\cup (\fh^{\ast})_{\mathrm{int}}^0\cup (\fh^{\ast})_{\mathrm{int}}^-,
\quad\hbox{where}\quad
\begin{array}{l}
(\fh^{\ast})_{\mathrm{int}}^+ = \CC\delta + \ZZ_{\geq 0}\linspan\{\Lambda_0,\ldots, \Lambda_n\}, \\ \\
(\fh^{\ast})_{\mathrm{int}}^0 = \CC\delta + 0\Lambda_0
+ \ZZ_{\geq 0}\linspan\{\omega_1, \ldots, \omega_n\}, \\ \\
(\fh^{\ast})_{\mathrm{int}}^- = \CC\delta + \ZZ_{\leq 0}\linspan\{\Lambda_0,\ldots, \Lambda_n\}.
\end{array}
\label{hintdefn}
\end{equation}
For $\widehat{\fsl}_2$ these sets are pictured (mod $\delta$) in \eqref{integrablewtspicture}.

For $i\in \{0,1,\ldots, n\}$ let 
$\mathbf{U}_{(i)}$ be the subalgebra of $\mathbf{U}$ generated by $\{E_i, F_i, K^{\pm1}_i\}$.
An \textit{integrable $\mathbf{U}$-module} is a $\mathbf{U}$-module $M$ such that
if $i\in \{0,\ldots , n\}$ then 
$$
\Res_{\mathbf{U}_{(i)}}^{\mathbf{U}}(M) \qquad
\text{is a direct sum of finite dimensional $\mathbf{U}_{(i)}$-modules,}
$$
where $\Res_{\mathbf{U}_{(i)}}^{\mathbf{U}}(M)$ denotes the restriction of the $\mathbf{U}$-module
$M$ to a $\mathbf{U}_{(i)}$-module.

Let $M$ be an integrable $\mathbf{U}$-module.  Following \cite[\S5]{Lu93},
for each $w\in W^{\mathrm{ad}}$ there is a linear map 
$$T_w\colon M\to M
\quad\hbox{such that}\quad
T_w(um) = T_w(u)T_w(m),
$$
for $u\in \mathbf{U}$ and $m\in M$ (here $T_w(u)$ refers to the braid group action 
on $\mathbf{U}$).  Thus, every integrable module $M$ is a module
for the semidirect product $\cB^{\mathrm{ad}}\ltimes \mathbf{U}$ where
$\cB^{\mathrm{ad}}$ is the braid group of $W^{\mathrm{ad}}$.

\subsection{Extremal weight modules $L(\Lambda)$}\label{Extwtmodules}

Let $\Lambda\in (\fh^*)_{\mathrm{int}}$.  Following \cite[(8.2.2)]{Kas94} and \cite[\S3.1]{Kas02},
the \emph{extremal weight module $L(\Lambda)$} is the $\mathbf{U}$-module
\begin{align}
\hbox{generated by $\{u_{w\Lambda}\mid w\in W\}$}
\qquad\hbox{with relations}\quad
K_i(u_{ w\Lambda}) =q^{\langle w\Lambda, \alpha_i^{\vee}\rangle} u_{ w\Lambda},
\nonumber \\
E_i u_{w\Lambda}=0,
\quad\hbox{and}\quad
F_i^{\langle w\Lambda, \alpha_i^{\vee} \rangle} u_{w\Lambda} = u_{s_i w\Lambda},
\qquad\hbox{if $\langle w\Lambda, \alpha_i^{\vee} \rangle\in \ZZ_{\geq 0} $},
\label{KMextremalpres} \\
F_i u_{w\Lambda}=0, 
\quad\hbox{and}\quad
E_i^{-\langle w\Lambda, \alpha_i^{\vee} \rangle} u_{w\Lambda}=u_{s_i w\Lambda},
\qquad\hbox{if $\langle w\Lambda, \alpha_i^{\vee} \rangle\in \ZZ_{\leq 0} $},
\nonumber
\end{align}
for $i\in \{0, \cdots, n\}$.
The module $L(\Lambda)$ has a crystal $B(\Lambda)$ (\cite[Prop.\ 8.2.2(ii)]{Kas94},
\cite[\S3.1]{Kas02}).
\begin{enumerate}
\item[$\bullet$] If $\Lambda\in (\fh^*)_{\mathrm{int}}^+$ then 
$L(\Lambda)$ is the simple $\mathbf{U}$-module of highest weight $\Lambda$ (see \cite[(10.4.6)]{Kac}). 
\item[$\bullet$] If $\Lambda\not\in (\fh^*)_{\mathrm{int}}^+$ then 
$L(\Lambda)$ is not a highest weight module.
\item[$\bullet$] If $\Lambda\in (\fh^*)_{\mathrm{int}}^-$ then 
$L(\Lambda)$ is the simple $\mathbf{U}$-module of lowest weight $\Lambda$.
\end{enumerate}
The finite dimensional simple modules, denoted
$L^{\mathrm{fin}}(a(u)))$, are integrable weight modules which are not extremal weight modules.
The connection between the $L^{\mathrm{fin}}(a(u))$ and the $L(\lambda)$ for
$\lambda\in (\fh^*)^0_{\mathrm{int}}$ is given by Theorem \ref{findimmods} below.

The module $L(\Lambda)$ is universal (see \cite[\S2.6]{Kas05}, \cite[\S2.1]{B02}, \cite[\S2.5]{Nak02}).  One way to formulate this universality
is to let $\mathbf{U}_0$ be the subalgebra generated by $K^{\pm1}_1,\ldots, K^{\pm1}_n, 
C^{\pm\frac12}, D^{\pm1}$, 
let $\mathrm{intInd}$ be an induction functor in the category of integrable $\mathbf{U}$-modules
and write
$$L(\Lambda) 
= \mathrm{int}\Ind_{\mathbf{U}_0 \rtimes \cB^{\mathrm{ad}}}^{\mathbf{U}}(S(\Lambda)),
\qquad\hbox{where}\quad
S(\Lambda) = \hbox{span}\{ u_{w\Lambda}\ |\ w\in W^{\mathrm{ad}} \}
$$
is the $\mathbf{U}_0\rtimes \cB^{\mathrm{ad}}$-module
with action given by
$T_i u_{w\Lambda} = (-q)^{\langle w\Lambda_i, \alpha_i^\vee\rangle} u_{s_iw\Lambda}$
and
$K_i u_{w\Lambda} = q^{\langle w\Lambda, \alpha_i^\vee \rangle} u_{w\Lambda},$
for $i\in \{0,1,\ldots, n\}$ and $w\in W^{\mathrm{ad}}$.

\subsection{Demazure submodules $L(\Lambda)_{\le w}$ }\label{Demazurechars}

Let $w\in W^{\mathrm{ad}}$.  
The \emph{Demazure module} $L(\Lambda)_{\le w}$ is the 
$\mathbf{U}^+$-submodule of $L(\Lambda)$ given by
$$L(\Lambda)_{\le w} = \mathbf{U}^+ u_{w\Lambda}
\qquad\hbox{and}\qquad
\mathrm{char}(L(\Lambda)_{\le w}) 
= \sum_{p\in B(\Lambda)_{\le w}} e^{\mathrm{wt}(p)},$$
since $L(\Lambda)_{\le w}$ has a crystal $B(\Lambda)_{\le w}$. 
The \emph{BGG-Demazure operator} on 
$\CC[\fh^*_\ZZ] 
= \hbox{$\CC$-span}\{ X^\lambda\ |\ \lambda\in \fh^*_\ZZ\}$ is given by
$$D_i = (1+s_i)\frac{1}{1-X^{-\alpha_i}}, \qquad
\hbox{for $i\in \{0, 1,\ldots, n\}$}.$$
Let $\Lambda\in (\fh^*)_{\mathrm{int}}$, $w\in W$ and $i\in \{0, 1, \ldots, n\}$.
\begin{equation*}
\hbox{If $\Lambda\in (\fh^*)_{\mathrm{int}}^+$\quad then \qquad}
D_i \mathrm{char}(L(\Lambda)_{\le w})) = \begin{cases}
\mathrm{char}(L(\Lambda)_{\le s_iw}), &\hbox{if $s_i w \posgeq \ w$}, \\
\mathrm{char}(L(\Lambda)_{\le w}), &\hbox{if $s_iw\posleq \ w$;}
\end{cases}
\end{equation*}
$$\hbox{if $\lambda\in (\fh^*)_{\mathrm{int}}^0$\quad then \qquad}
D_i \mathrm{char}(L(\lambda)_{\le w}) = \begin{cases}
\mathrm{char}(L(\lambda)_{\le s_iw}), &\hbox{if $s_i w \zerogeq\  w$}, \\
\mathrm{char}(L(\lambda)_{\le w}), &\hbox{if $s_iw\zeroleq\ w$;}
\end{cases}
$$
$$\hbox{if $\Lambda\in (\fh^*)_{\mathrm{int}}^-$\quad then \qquad}
D_i \mathrm{char}(L(\Lambda)_{\le w}) = \begin{cases}
\mathrm{char}(L(\Lambda)_{\le s_iw}), &\hbox{if $s_i w \neggeq\ w$}, \\
\mathrm{char}(L(\Lambda)_{\le w}), &\hbox{if $s_iw\negleq\ w$;}
\end{cases}
$$
(see \cite[Theorem 8.2.9]{Kum}, \cite{Kas93}, \cite[\S2.8]{Kas05} and 
\cite[Theorems 4.7 and 4.11]{Kt16}).

\subsection{An alternate presentation for level 0 extremal weight modules}
\label{levelzeroloopsection}

For $\lambda = m_1\omega_1+\cdots +m_n\omega_n\in (\fh^*)_{\mathrm{int}}$ and let
$x_{1,1}, \ldots, x_{m_1,1}$,\ 
$x_{1,2}, \ldots, x_{m_2,2}, \ \ldots\ $
$x_{1,n}, \ldots, x_{m_n,n}$
be $n$ sets of formal variables.
Letting $e_i^{(j)} = e_i(x_{1,j}, \ldots, x_{m_j,j})$ denote the elementary symmetric function
in the variables $x_{1,j}, \ldots, x_{m_j,j}$ define
\begin{align*}
RG_\lambda 
&= \CC[x^{\pm1}_{1,1}, \ldots, x^{\pm1}_{m_1,1}]^{S_{m_1}}\otimes
\cdots \otimes \CC[x^{\pm1}_{1,n}, \ldots, x^{\pm1}_{m_n,n}]^{S_{m_n}} \\
&= 
\CC[e_1^{(1)}, \ldots, e_{m_1-1}^{(1)}, (e_{m_1}^{(1)})^{\pm1}]
\otimes \cdots \otimes
\CC[e_1^{(n)}, \ldots, e_{m_n-1}^{(n)}, (e_{m_n}^{(n)})^{\pm1}], \\
RG^+_\lambda &= 
\CC[x_{1,1}, \ldots, x_{m_1,1}]^{S_{m_1}}
\otimes \cdots \otimes
\CC[x_{1,n}, \ldots, x_{m_n,n}]^{S_{m_n}}, \quad\hbox{and} \\
RG^-_\lambda &= 
\CC[x^{-1}_{1,1}, \ldots, x^{-1}_{m_1,1}]^{S_{m_1}}
\otimes \cdots \otimes
\CC[x^{-1}_{1,n}, \ldots, x^{-1}_{m_n,n}]^{S_{m_n}}.
\end{align*}

Let 
\begin{align*}
e^{(i)}_+(u) &= (1-x_{1,i}u)(1-x_{2,i}u)\cdots (1-x_{m_i,i}u)
\quad\hbox{and} \\
e^{(i)}_-(u^{-1}) &= (1-x_{1,i}^{-1}u^{-1})(1-x_{2,i}^{-1}u^{-1})\cdots (1-x_{m_i,i}^{-1}u^{-1}).
\end{align*}
Let $\mathbf{U}'$ be the subalgebra of $\mathbf{U}$ without the generator $D$.

\begin{thm}\label{exwtequalsDVerma}
\emph{(see \cite[\S3.4]{Nak02})}
The extremal weight module $L(\lambda)$ is the
$(\mathbf{U}' \otimes_\ZZ RG_\lambda)$-module
generated by a single vector $m_\lambda$ with relations 
$$\mathbf{x}_{i,r}^+ m_\lambda = 0, \qquad
K_i m_\lambda = q^{m_i} m_\lambda,
\quad
Cm_\lambda = m_\lambda,
$$
$$\mathbf{q}^{(i)}_+(u) m_\lambda = K_i \frac{e^{(i)}_+(q^{-1}u) }{e^{(i)}_+ (qu) }m_\lambda
\qquad\hbox{and}\qquad
\mathbf{q}^{(i)}_-(u^{-1}) m_\lambda = K_i \frac{e^{(i)}_-(qu^{-1}) }{e^{(i)}_- (q^{-1}u^{-1}) }m_\lambda,
$$
where $\mathbf{q}^{(i)}_+(u)$ and $\mathbf{q}^{(i)}_-(u^{-1})$ 
are as defined in \eqref{qplusseries} and \eqref{qminusseries}.
\end{thm}

\noindent
In this form $L(\lambda)$ has been termed the \emph{universal standard module} 
\cite[\S3.4]{Nak02}) or the \emph{global Weyl module} \cite[\S2]{CP01}.  See
\cite[Theorem 2]{Nak02} and \cite[Remark 2.15]{Nak02} for discussion of how to see that 
the extremal weight module, the universal standard module and the global Weyl module coincide.

\begin{remark}\label{RGchar}
Let
$$
0_q = \frac{1}{1-q} + \frac{q^{-1}}{1-q^{-1}} = \cdots + q^{-3}+q^{-2} + q^{-1} +1+q+q^2+\cdots,$$
(although
$\frac{q^{-1}}{1-q^{-1}} = \frac{1}{q-1} = \frac{-1}{1-q}$, 
it is important to note that $0_q$ is \emph{not} equal to $0$, it is a doubly infinite formal series 
in $q$ and $q^{-1}$).
Since $\deg(e_j^{(i)}) = j$,
$$
\mathrm{gchar}(RG_\lambda^+)
=\left(\prod_{i=1}^n \prod_{k=1}^{m_i} \frac{1}{1-q^k}\right)
\qquad
 \mathrm{gchar}(RG_\lambda^-)
= 
\left(\prod_{i=1}^n \prod_{k=1}^{m_i} \frac{1}{1-q^{-k} }\right)
\quad\hbox{and}
$$
$$
\mathrm{gchar}(RG_\lambda)
=
\Big(0_{q^{m_1}}\prod_{k=1}^{m_1-1} \frac{1}{1-q^k}\Big)
\Big(0_{q^{m_2}}\prod_{k=1}^{m_2-1} \frac{1}{1-q^k}\Big)
\cdots \Big(0_{q^{m_n}}\prod_{k=1}^{m_n-1} \frac{1}{1-q^k}\Big).$$
\end{remark}

\subsection{Level 0 $L(\lambda)$ and finite dimensional simple 
$\mathbf{U}$-modules $L^{\mathrm{fin}}(a(u))$}\label{findimmodules}

The loop presentation provides a triangular decomposition of $\mathbf{U}$ (different
form the usual triangular decomposition coming from the Kac-Moody presentation).
The extremal weight module $L(\lambda)$ is
the standard (Verma type) module for the loop triangular decomposition
(see \cite[Theorem 2.3(b)]{CP95},
\cite[Lemma 2.14]{Nak02} and \cite[outline of proof of Theorem 12.2.6]{CP94}).

A \emph{Drinfeld polynomial} is an $n$-tuple of polynomials 
$a(u) = (a^{(1)}(u), \ldots, a^{(n)}(u))$ with
$a^{(i)}(u)\in \CC[u]$, represented as
$$a(u) = a^{(1)}(u)\omega_1+\cdots +a^{(n)}(u)\omega_n,
\qquad\hbox{with}\quad
a^{(i)}(u) = (u-a_{1,i})\cdots (u-a_{m_i,i})
$$
so that the coefficient of $u^j$ in $a^{(i)}(u)$ is $e^{(i)}_{m_i-j}(a_{1,i},\ldots, a_{m_i,i})$,
the $(m_i-j)$th elementary symmetric function evaluated at the values $a_{1,i},\ldots, a_{m_i,i}$.
The \emph{local Weyl module} (a finite dimensional standard module) is defined by
$$M^{\mathrm{fin}}(a(u)) = L(\lambda) \otimes_{RG_\lambda} m_{a(u)},
\quad\hbox{where}\quad
e_k^{(i)}(x_{1,i},x_{2,i}, \ldots) m_{a(u)} = e_k^{(i)}(a_{1,i}, \ldots, a_{m_i,i}) m_{a(u)}
$$
specifies the $RG_\lambda$-action on $m_{a(u)}$.  In other words, the module
$M^{\mathrm{fin}}(a(u))$ is $L(\lambda)$ except that variables $x_{j,i}$ from \S \ref{levelzeroloopsection} 
have been specialised to the values $a_{j,i}$.
As in Theorem \ref{exwtequalsDVerma}, 
let $\mathbf{U}'$ be the subalgebra of $\mathbf{U}$ without the generator $D$.

\begin{thm}\label{findimmods}\emph{ (see \cite[Theorem 2]{Dr88} and \cite[Theorem 3.3]{CP95})}
The standard module
$M^{\mathrm{fin}}(a(u))$ has a unique simple quotient
$L^{\mathrm{fin}}(a(u))$ and
$$\begin{matrix}
\{ \hbox{Drinfeld polynomials}\}
&\longrightarrow 
&\{ \hbox{finite dimensional simple $\mathbf{U'}$-modules}\} \\
a(u) = a^{(1)}(u)\omega_1+\cdots+a^{(n)}(u)\omega_n
&\longmapsto
&L^{\mathrm{fin}}(a(u))
\end{matrix}
$$
is a bijection.
\end{thm}


\subsection{Path models for the crystals $B(\Lambda)$}

The work of Littelmann \cite{Li94}, \cite{Li95}
provided a particularly convenient model for the crystals $B(\Lambda)$
when $\Lambda$ is positive or negative level.  This model realizes the crystal as a set of 
paths $p\colon \RR_{[0,1]}\to \fh^*$ with combinatorially defined 
Kashiwara operators $\tilde e_0, \ldots, \tilde e_n,
\tilde f_0, \ldots, \tilde f_n$.  In the LS (Lakshmibai-Seshadri) model the generator of the crystal
$B(\Lambda)$ is the straight line path to $\Lambda$.

When $\lambda\in (\fh^*)_{\mathrm{int}}^0$ using 
$$
\begin{matrix}
p_\lambda\colon &\RR_{[0,1]} &\to &\fh^* \\
&t &\to &t\lambda\, ,
\end{matrix}
\qquad\hbox{the straight line path from $0$ to $\lambda$,}
$$
as a generator for $B(\lambda)$ may not be the optimal choice.  Remarkably,
Naito and Sagaki 
(see \cite[Definition 3.1.4 and Theorem 3.2.1]{INS16} and \cite[Theorem 4.6.1(b)]{NNS15}), 
have shown that  $B(\lambda)$
can be constructed with sequences of Weyl group elements and rational numbers
as in \cite[\S1.2, 1.3 and 2.2]{Li94} but with the positive level length $\ell^+$ 
and Bruhat order $\posleq\ $ replaced
by the level zero length $\ell^0$ and Bruhat order $\zeroleq\ $.
However, when working with the Naito-Sagaki construction one must be very careful 
not to identify the Naito-Sagaki sequences with actual paths 
(piecewise linear maps from $\RR_{[0,1]}$ to $\fh^*$)
because the natural map from
Naito-Sagaki sequences to paths  is \emph{not always injective}
(an example is provided by \cite[Remark 5.10]{Kas02}).

If $\Lambda\in (\fh^*)_{\mathrm{int}}^+$ then
$-\Lambda\in (\fh^*)_{\mathrm{int}}^-$ and
$$B(\Lambda) = \{ \tilde f_{i_1}\cdots \tilde f_{i_k} p_\Lambda\ |\ 
\hbox{$k\in \ZZ_{\ge 0}$ and ${i_1}, \ldots, i_k\in \{0,1,\ldots, n\}$} \}$$
$$B(-\Lambda) = \{ \tilde e_{i_1}\cdots \tilde e_{i_k} p_{-\Lambda}\ |\ 
\hbox{$k\in \ZZ_{\ge 0}$ and $i_1, \ldots, i_k\in \{0,1,\ldots, n\}$} \}.
$$
are each a single connected component and their characters are determined by the Weyl-Kac character formula
\cite[Theorem 11.13.3]{Kac}.

%
%
%

\subsection{Crystals for level 0 extremal weight modules $L(\lambda)$}

For general $\lambda\in (\fh^*)_{\mathrm{int}}^0$ the crystal $B(\lambda)$ is not
connected (as a graph with edges determined by the Kashiwara
operators $\tilde e_0, \ldots, \tilde e_n,
\tilde f_0, \ldots \tilde f_n$).  
Let $\lambda = m_1\omega_1+\cdots+m_n\omega_n$, with $m_1, \ldots, m_n\in \ZZ_{\ge 0}$.
By \cite[Corollary 4.15]{BN02}, the map
\begin{equation}
\begin{matrix}
\Phi_\lambda\colon &L(\lambda) &\longrightarrow 
&L(\omega_1)^{\otimes m_1}\otimes \cdots \otimes L(\omega_n)^{\otimes m_n} \\
&u_\lambda &\longmapsto &u_{\omega_1}^{\otimes m_1}\otimes\cdots\otimes u_{\omega_n}^{\otimes m_n}
\end{matrix}
\qquad\hbox{is injective}
\label{tensorinjection}
\end{equation}
and gives rise to an injection of crystals
\begin{equation*}
B(\lambda) \hookrightarrow B(\omega_1)^{\otimes m_1}\otimes \cdots \otimes B(\omega_n)^{\otimes m_n}
\end{equation*}
which takes the connected component of $B(\lambda)$ containing $b_\lambda$ to
the connected component of $B(\omega_1)^{\otimes m_1}\otimes \cdots \otimes B(\omega_n)^{\otimes m_n}$ containing $b_{\omega_1}^{\otimes m_1}\otimes\cdots\otimes b_{\omega_n}^{\otimes m_n}$.
Kashiwara \cite[Theorem 5.15]{Kas02} fully described the structure of 
$L(\omega_i)$ (see \cite[Theorem 2.16]{BN02}).
Beck-Nakajima analyzed the PBW basis of \eqref{PBWbasis}  
and use \eqref{tensorinjection} to show that
the connected components of $B(\lambda)$ are labeled by $n$-tuples of partitions 
$\kappa = (\kappa^{(1)}, \ldots, \kappa^{(n)})$ such that $\ell(\kappa^{(i)})<m_i$.
Together with a result of Fourier-Littelmann 
which shows that the
crystal of the level 0 module $M^{\mathrm{fin}}(a(u))$ is isomorphic to a level one Demazure
crystal $B(\nu+\Lambda_0)_{\le w}$, 
the full result is as detailed in Theorem \ref{levelzerostructure} below.

A labeling set for a basis of $RG_\lambda/\langle e^{(i)}_{m_i}=1\rangle$ is
\begin{align*}
S^\lambda
&=\{
\vec \kappa = (\kappa^{(1)}, \ldots, \kappa^{(n)}) \mid \text{$\kappa^{(i)}$ is a partition with $\ell(\kappa^{(i)}) < m_i$ for $i\in \{1, \ldots, n\}$}
\}.
\end{align*}
The connected component of $b_\lambda$ in $B(\lambda)$ is
$$B(\lambda)_0
=  \Big\{ \tilde r_{i_1}\cdots \tilde r_{i_k} b_\lambda \ |\ 
\hbox{$k\in \ZZ_{\ge 0}$ and  $\tilde r_{i_1}, \ldots, \tilde r_{i_k}\in 
\{\tilde e_0, \ldots, \tilde e_n, \tilde f_0, \ldots, \tilde f_n\}$} \Big\}.
$$
Define $B^{\mathrm{fin}}(\lambda)$ to be the ``crystal'' which has a crystal graph which is
the ``quotient'' of the crystal graph of $L(\lambda)$ obtained by identifying
the vertices $b$ and $b'$ if there is an element $\mathbf{s}\in RG_\lambda$ such that
$\mathbf{s}G(b) = G(b')$, where $G(b)$ denotes the canonical basis element of $L(\lambda)$ corresponding
to $b$.
$$\hbox{$B^{\mathrm{fin}}(\lambda)$ is the crystal of the finite dimensional standard module
$M^{\mathrm{fin}}(a(u))$.}$$

\begin{thm} \label{levelzerostructure}  \emph{(see \cite[Theorem 4.16]{BN02}, 
\cite[\S3.4]{Nak02}) and \cite[Proposition 3]{FL05})}
Let $\lambda = m_1\omega_1+\cdots + m_n\omega_n\in 
(\fh^*)_{\mathrm{int}}^0$.  As in Theorem \ref{MacisDem},
let $\nu\in \fa_\ZZ^*$, $j\in \ZZ_{\ge 0}$ and $w\in W^{\mathrm{ad}}$ such
that $w(\nu+\Lambda_0) = -j\delta + \lambda+\Lambda_0$ and $w$ is minimal length.
Then 
$$B(\lambda) \simeq B(\lambda)_0\times S^\lambda
\quad\hbox{and}\quad
B(\lambda)_0 \simeq \ZZ^k \times B^{\mathrm{fin}}(\lambda)
\quad\hbox{and}\quad
B^{\mathrm{fin}}(\lambda) \simeq B(\nu+\Lambda_0)_{\le w},
$$
where $k$ is the number of elements of $m_1, \ldots, m_n$ which are nonzero.
\end{thm}

\noindent
Additional useful references for Theorem \ref{levelzerostructure}
are \cite[Theorem 1]{Nak02} and \cite[Theorem 1]{B02}.
The first two statements in Theorem \ref{levelzerostructure}
are reflections of the very important fact that $L(\lambda)$ is free as an $RG_\lambda$-module.
This fact that $L(\lambda)$ is free as an $RG_\lambda$-module was deduced
geometrically, via the quiver variety, in \cite[Theorem 7.3.5]{Nak99} 
(note that property $(T_{G_\mathbf{w}\times\CC^\times})$ 
there includes the freeness, see the definition of property $(T_G)$ after \cite[(7.1.1)]{Nak99}).
This freeness was understood more algebraically in the work of Fourier-Littelmann \cite{FL05}
and Chari-Ion \cite[Cor.\ 2.10]{CI13}.  Further understanding of the 
$RG_\lambda$-action in terms of the geometry of the semi-infinite flag variety 
is in \cite[\S5.1]{BF13}.

The last statement of Theorem \ref{levelzerostructure}
is proved  by considering the map
$$
\begin{matrix}
B(\lambda) &\longrightarrow &B(\Lambda_0)\otimes B(\lambda) \\
\cup\vert &&\cup\vert  \\
B^{\mathrm{fin}}(\lambda) &\stackrel{\sim}{\longrightarrow} &B(\nu+\Lambda_0)_{\le w} 
\end{matrix}
\qquad\hbox{given by}\qquad
b \mapsto b_{\Lambda_0}\otimes b.
$$
where $b_{\Lambda_0}$ is the highest weight of the crystal $B(\Lambda_0)$.
Combining the isomorphism $B^{\mathrm{fin}}(\lambda) \simeq B(\nu+\Lambda_0)_{\le w}$ with
Theorem \ref{MacisDem} and the positive level formula in \S \ref{Demazurechars} gives
$$
\mathrm{char}(M(a(u))) = \mathrm{char}(L(\nu+\Lambda_0)_{\le w}) \\
= q^{-j} X^{-\Lambda_0}\tilde E_\lambda(q,0)
\qquad\hbox{and}$$
\begin{equation}
\mathrm{char}(L(\lambda)_{\le w_0})
= \mathrm{gchar}(RG_\lambda^+)\tilde E_{w_0\lambda}(q,0),
\qquad
\mathrm{char}(L(\lambda)) = \mathrm{gchar}(RG_\lambda) \tilde E_{w_0\lambda}(q,0).
\label{levzerochar}
\end{equation}

\section{Examples for $\fg = \widehat{\fsl}_2$}


Let $\fg = \widehat{\mathfrak{sl}}_2$ with
$\fh^*=\mathbb{C}\omega_1\oplus \mathbb{C} \Lambda_0\oplus \mathbb{C} \delta$
with affine Cartan matrix
$$\begin{pmatrix}
\alpha_0(h_0) &\alpha_0(h_1) \\
\alpha_1(h_0) &\alpha_1(h_1)
\end{pmatrix}
= \begin{pmatrix} 2 &-2 \\ -2 &2\end{pmatrix}
\qquad\hbox{and}\qquad
\begin{array}{ll}
\theta = \alpha_1 = 2\omega_1,\qquad
&\theta^\vee = \alpha_1^\vee = h_1, \\
\Lambda_1 = \omega_1+\Lambda_0,
&\alpha_0 = -\alpha_1 + \delta,
\end{array}
$$
Using that $\langle \alpha_1, \alpha_1\rangle =2$,
if $k\in \ZZ$ and $\mu^\vee = k\alpha_1^\vee = k\alpha_1 = k 2\omega_1 = 2k\omega_1$
so that $\mu_1^\vee = 2k$ and
$-\hbox{$\frac12$}\langle \mu^\vee, \mu^\vee\rangle
=-\hbox{$\frac12$} (k2k) = -k^2$.
Thus, following \eqref{tmumatrix}, \eqref{simatrix} and \eqref{szeromatrix},
in the basis $\{ \delta, \omega_1, \Lambda_0\}$ of $\fh^*$, 
$$t_{k\alpha^\vee} = \begin{pmatrix} 1 &-k &-k^2 \\
0 &1 &2k \\
0 &0 &1 \end{pmatrix},
\qquad
s_1 = \begin{pmatrix} 1 &0 &0\\
0 &-1 &0 \\
0 &0 &1 \end{pmatrix},
\qquad
s_0 = \begin{pmatrix} 1 &1 &-1 \\
0 &-1 &2 \\
0 &0 &1 \end{pmatrix}.
$$
These matrices are used to compute the $W$-orbits pictured in Section
\ref{Worbits}.

\subsection{Macdonald polynomials}

The Demazure operators are given by
$$D_1f = \frac{f - X^{-2\omega_1}(s_1f)}{1-X^{-2\omega_1}}
\qquad\hbox{and}\qquad
D_0f =  \frac{f - X^{-\alpha_0}(s_0f)}{1-X^{-\alpha_0}}
=
\frac{f - q^{-1}X^{2\omega_1}(s_0f)}{1-q^{-1}X^{2\omega_1}}.
$$
The normalized Macdonald polynomials $\tilde E_{\omega_1}(q,t)$ 
and $\tilde E_{-\omega_1}(q,t)$
are
$$\tilde E_{\omega_1}(q,t) = X^{\omega_1}
\qquad\hbox{and}\qquad
\begin{array}{cc}
\beginpicture
\setcoordinatesystem units <.7cm,.7cm>         
\setplotarea x from -2.5 to 2.5, y from -0.5 to 0.5  
    \plot -2.2 0 2.2 0 /
    \plot  -2 -0.2  -2 0.5 /
    \plot  -1 -0.2  -1 0.5 /
    \plot  0 -0.2  0 0.5 /
    \plot  1 -0.2  1 0.5 /
    \plot 2 -0.2 2 0.5 /
    \put{$\bullet$} at 0 0
    \arrow <5pt> [.2,.67] from 0.5 0.2 to -0.5 0.2
    \put{$\displaystyle{\tilde E_{-\omega_1}(q,t) = X^{-\omega_1}}$} at 0 -1.5
    \put{$\displaystyle{\phantom{\tilde E_{-omega_1}(q,t) }= X^{-\omega_1}}$} at 0 -3
\endpicture
&
\beginpicture
\setcoordinatesystem units <.7cm,.7cm>         
\setplotarea x from -2.5 to 2.5, y from -0.5 to 0.5  
    \plot -2.2 0 2.2 0 /
    \plot  -2 -0.2  -2 0.5 /
    \plot  -1 -0.2  -1 0.5 /
    \plot  0 -0.2  0 0.5 /
    \plot  1 -0.2  1 0.5 /
    \plot 2 -0.2 2 0.5 /
    \put{$\bullet$} at 0 0
    \plot 0.5 0.2 0 0.2 /
    \plot 0.95 0.2 0.95 0.35 /
    \arrow <5pt> [.2,.67] from 0 0.35 to 0.5 0.35
    \put{$\displaystyle{+\frac{(1-t)}{1-qt} X^{\omega_1}}$} at 0 -1.5
    \put{$\displaystyle{+\frac{(1-t^{-1})q^{-1}}{1-q^{-1}t^{-1}} X^{\omega_1} }$} at 0 -3
\endpicture
\end{array}
$$
giving $\tilde E_{\omega_1}(0,t) = \tilde E_{\omega_1}(\infty,t)=\tilde E_{\omega_1}(q,0)
=\tilde E_{\omega_1}(q,\infty) = X^{\omega_1}$ and
$$\begin{array}{ll}
\displaystyle{ \tilde E_{-\omega_1}(0,t) = X^{-\omega_1}+(1-t)X^{\omega_1}, }\qquad
&\displaystyle{ \tilde E_{-\omega_1}(\infty,t) = X^{-\omega_1}  } \\
\displaystyle{ \tilde E_{-\omega_1}(q,0) = X^{-\omega_1}+X^{\omega_1}, }
&\displaystyle{ \tilde E_{-\omega_1}(q,\infty) = X^{-\omega_1}+q^{-1}X^{\omega_1}.}
\end{array}
$$
The normalized Macdonald polynomials 
$\tilde E_{2\omega_1}(q,t)$ and $\tilde E_{-2\omega_1}(q,t)$ are
$$
\begin{array}{cc}
\beginpicture
\setcoordinatesystem units <.7cm,.7cm>         
\setplotarea x from -2.5 to 2.5, y from -0.5 to 0.5  
    \plot -2.2 0 2.2 0 /
    \plot  -2 -0.2  -2 0.5 /
    \plot  -1 -0.2  -1 0.5 /
    \plot  0 -0.2  0 0.5 /
    \plot  1 -0.2  1 0.5 /
    \plot 2 -0.2 2 0.5 /
    \put{$\bullet$} at 0 0
    \arrow <5pt> [.2,.67] from 0.5 0.2 to 1.5 0.2
    \put{$\displaystyle{\tilde E_{2\omega_1}(q,t) = X^{2\omega_1}}$} at 0 -1.5
    \put{$\displaystyle{\phantom{\tilde E_{2\omega_1}(q,t)} = X^{2\omega_1}}$} at 0 -3
\endpicture
&
\beginpicture
\setcoordinatesystem units <.7cm,.7cm>         
\setplotarea x from -2.5 to 2.5, y from -0.5 to 0.5  
    \plot -2.2 0 2.2 0 /
    \plot  -2 -0.2  -2 0.5 /
    \plot  -1 -0.2  -1 0.5 /
    \plot  0 -0.2  0 0.5 /
    \plot  1 -0.2  1 0.5 /
    \plot 2 -0.2 2 0.5 /
    \put{$\bullet$} at 0 0
    \plot 0.5 0.2 1 0.2 /
    \plot 0.95 0.2 0.95 0.35 /
    \arrow <5pt> [.2,.67] from 1 0.35 to 0.5 0.35
    \put{$\displaystyle{+\frac{(1-t)q}{1-qt}}$} at 0 -1.5
    \put{$\displaystyle{+\frac{1-t^{-1}}{1-q^{-1}t^{-1}} }$} at 0 -3
\endpicture
\end{array}
\qquad\quad\hbox{and}
$$
$$
\begin{array}{cccc}
\beginpicture
\setcoordinatesystem units <.7cm,.7cm>         
\setplotarea x from -2.5 to 2.5, y from -0.5 to 0.5  
    \plot -2.2 0 2.2 0 /
    \plot  -2 -0.2  -2 0.5 /
    \plot  -1 -0.2  -1 0.5 /
    \plot  0 -0.2  0 0.5 /
    \plot  1 -0.2  1 0.5 /
    \plot 2 -0.2 2 0.5 /
    \put{$\bullet$} at 0 0
    \arrow <5pt> [.2,.67] from 0.5 0.2 to -0.5 0.2
    \arrow <5pt> [.2,.67] from -0.5 0.2 to -1.5 0.2
    \put{$\displaystyle{\tilde E_{-2\omega_1}(q,t) = X^{-2\omega_1}}$} at 0 -1.5
    \put{$\displaystyle{\phantom{\tilde E_{-2\omega_1}(q,t)} =X^{-2\omega_1}}$} at 0 -3.5
\endpicture
&
\beginpicture
\setcoordinatesystem units <.7cm,.7cm>         
\setplotarea x from -2.5 to 2.5, y from -0.5 to 0.5  
    \plot -2.2 0 2.2 0 /
    \plot  -2 -0.2  -2 0.5 /
    \plot  -1 -0.2  -1 0.5 /
    \plot  0 -0.2  0 0.5 /
    \plot  1 -0.2  1 0.5 /
    \plot 2 -0.2 2 0.5 /
    \put{$\bullet$} at 0 0
    \arrow <5pt> [.2,.67] from 0.5 0.2 to -0.5 0.2
    \plot -0.5 0.2 -1 0.2 /
    \plot -0.95 0.2 -0.95 0.35 /
    \arrow <5pt> [.2,.67] from -1 0.35 to -0.5 0.35
    \put{$\displaystyle{+\frac{1-t}{1-qt}}$} at 0 -1.5
    \put{$\displaystyle{+\frac{(1-t^{-1})q^{-1}}{1-q^{-1}t^{-1}}}$} at 0 -3.5
\endpicture
&
\beginpicture
\setcoordinatesystem units <.7cm,.7cm>         
\setplotarea x from -2.5 to 2.5, y from -0.5 to 0.5  
    \plot -2.2 0 2.2 0 /
    \plot  -2 -0.2  -2 0.5 /
    \plot  -1 -0.2  -1 0.5 /
    \plot  0 -0.2  0 0.5 /
    \plot  1 -0.2  1 0.5 /
    \plot 2 -0.2 2 0.5 /
    \put{$\bullet$} at 0 0
    \plot 0.5 0.2 0 0.2 /
    \plot 0.05 0.2 0.05 0.35 /
    \arrow <5pt> [.2,.67] from 0 0.35 to 0.5 0.35
    \arrow <5pt> [.2,.67] from 0.5 0.35 to 1.5 0.35
    \put{$\displaystyle{+\frac{1-t}{1-q^2t}X^{2\omega_1}}$} at 0 -1.5
    \put{$\displaystyle{+\frac{(1-t^{-1})q^{-2}}{1-q^{-2}t^{-1}}X^{2\omega_1}}$} at 0 -3.5
\endpicture
&
\beginpicture
\setcoordinatesystem units <.7cm,.7cm>         
\setplotarea x from -2.5 to 2.5, y from -0.5 to 0.5  
    \plot -2.2 0 2.2 0 /
    \plot  -2 -0.2  -2 0.5 /
    \plot  -1 -0.2  -1 0.5 /
    \plot  0 -0.2  0 0.5 /
    \plot  1 -0.2  1 0.5 /
    \plot 2 -0.2 2 0.5 /
    \put{$\bullet$} at 0 0
    \plot 0.5 0.2 0 0.2 /
    \plot 0.15 0.2 0.15 0.35 /
    \arrow <5pt> [.2,.67] from 0 0.35 to 0.5 0.35
    \plot 0.5 0.35 1 0.35 /
    \plot 0.95 0.35 0.95 0.5 /
    \arrow <5pt> [.2,.67] from 1 0.5 to 0.5 0.5
    \put{$\displaystyle{+\frac{(1-t)}{(1-q^2t)}\frac{(1-t)q}{(1-qt)}}$} at 0 -1.5
    \put{$\displaystyle{+\frac{(1-t^{-1})q^{-2}}{(1-q^{-2}t^{-1})}\frac{(1-t^{-1})}{(1-q^{-1}t^{-1})}}$} at 0 -3.5
\endpicture
\end{array}
$$
giving
$$\begin{array}{ll}
\displaystyle{ \tilde E_{2\omega_1}(0,t) = X^{2\omega_1}, }
&\displaystyle{ \tilde E_{2\omega_1}(\infty,t) = X^{2\omega_1} + (1-t^{-1}) } \\
\displaystyle{ \tilde E_{2\omega_1}(q,0) = X^{2\omega_1}+q, }
&\displaystyle{ \tilde E_{2\omega_1}(q,\infty) = X^{2\omega_1}+1}
\end{array}
$$
and
$$\begin{array}{ll}
\tilde E_{-2\omega_1}(0,t) = X^{-2\omega_1} + (1-t)X^{2\omega_1} +(1-t), \qquad
&\tilde E_{-2\omega_1}(\infty, t) = X^{-2\omega_1}, \\
\tilde E_{-2\omega_1}(q,0) = X^{-2\omega_1} + X^{2\omega_1} +1+q, \qquad
&\tilde E_{-2\omega_1}(q,\infty) = X^{-2\omega_1} + q^{-2}X^{2\omega_1}
+q^{-1}+q^{-2}.
\end{array}
$$

\subsection{The crystal $B(\Lambda_0)$ and $B(\omega_1+\Lambda_0)$}
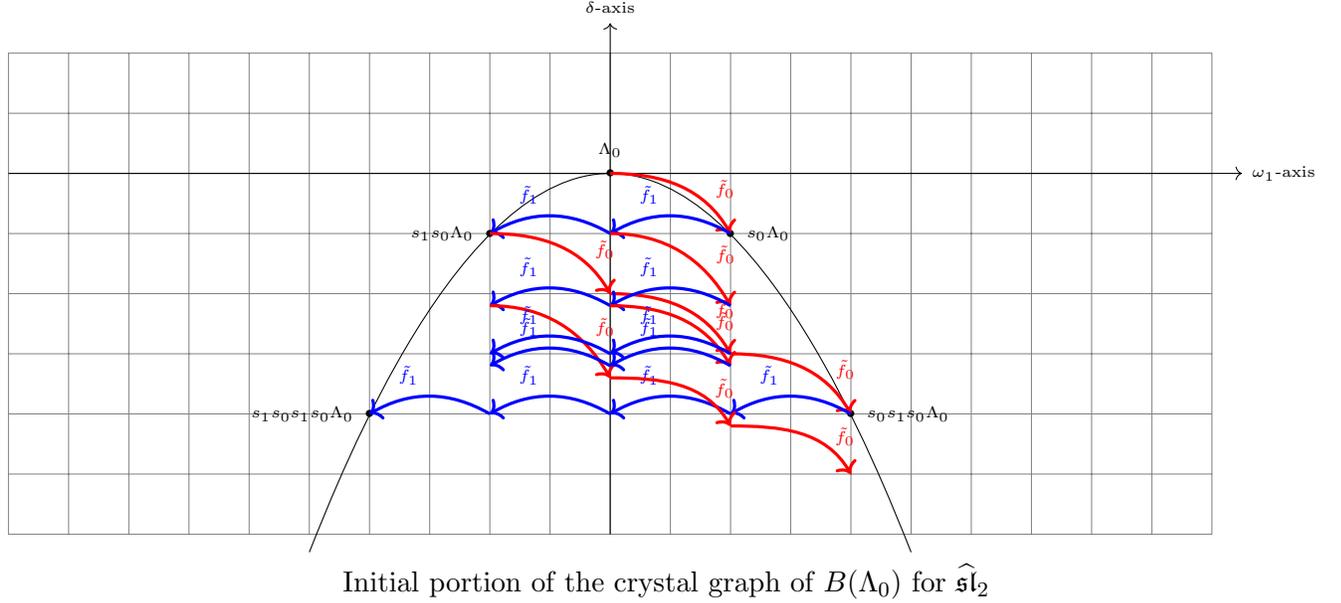
\label{Ionthmexamples}

$$
\begin{matrix}
\begin{tikzpicture}[scale=0.8]
\tikzstyle{every node}=[font=\tiny]    
    \draw[style=help lines,step=1cm] (-10,-6) grid (10,2);   
    \draw[->] (-10,0) -- (10.5,0) node[anchor=west]{$\omega_1$-axis};
    \draw[->] (0,-6) -- (0, 2.5) node[anchor=south]{$\delta$-axis};

\node at (0,0) {$\bullet$};
        \filldraw (0,0.1) node[anchor=south] {$\Lambda_0$};
        \filldraw (2.1,-1) node[anchor=west] {$s_0\Lambda_0$};
        \filldraw (-2.1,-1) node[anchor=east] {$s_1s_0\Lambda_0$};
        \filldraw (4.1,-4) node[anchor=west] {$s_0s_1s_0\Lambda_0$};
        \filldraw (-4.1,-4) node[anchor=east] {$s_1s_0s_1s_0\Lambda_0$};
\node at (2,-1) {$\bullet$};
\node at (-2, -1) {$\bullet$};  
\node at (4, -4) {$\bullet$};  
\node at (-4, -4) {$\bullet$};  
		\draw (0, 0)parabola (-5, -6.3) ;
		\draw (0, 0)parabola (5, -6.3) ;
\draw [<-,red, in=0, out=120,  very thick] (2,-1) to node[right, xshift=9, yshift=-1]{$\tilde{f}_0$} (0,0);
\draw [<-,blue, out=30, in=150, very thick] (0,-1) to node[left, yshift=8]{$\tilde{f}_1$} (2,-1);
\draw [<-,blue, out=30, in=150, very thick] (-2,-1) to node[left, yshift=8]{$\tilde{f}_1$} (0,-1);
\draw [<-,red, in=0, out=120,  very thick] (2,-2.2) to node[right, xshift=9, yshift=-1]{$\tilde{f}_0$} (0,-1);
\draw [<-,red, in=0, out=120,  very thick] (0,-2) to node[right, xshift=9, yshift=-1]{$\tilde{f}_0$} (-2,-1);
\draw [<-,red, in=0, out=120,  very thick] (2,-3) to node[right, xshift=9, yshift=-1]{$\tilde{f}_0$} (0,-2);
\draw [<-,red, in=0, out=120,  very thick] (4,-4) to node[right, xshift=9, yshift=-1]{$\tilde{f}_0$} (2,-3);
\draw [<-,blue, out=30, in=150, very thick] (0,-2.2) to node[left, yshift=8]{$\tilde{f}_1$} (2,-2.2);
\draw [<-,blue, out=30, in=150, very thick] (-2,-2.2) to node[left, yshift=8]{$\tilde{f}_1$} (0,-2.2);
\draw [<-,blue, out=30, in=150, very thick] (2,-4) to node[left, yshift=8]{$\tilde{f}_1$} (4,-4);
\draw [<-,blue, out=30, in=150, very thick] (0,-4) to node[left, yshift=8]{$\tilde{f}_1$} (2,-4);
\draw [<-,blue, out=30, in=150, very thick] (-2,-4) to node[left, yshift=8]{$\tilde{f}_1$} (0,-4);
\draw [<-,blue, out=30, in=150, very thick] (-4,-4) to node[left, yshift=8]{$\tilde{f}_1$} (-2,-4);
\draw [<-,red, in=0, out=120,  very thick] (2,-3.2) to node[right, xshift=9, yshift=-1]{$\tilde{f}_0$} (0,-2.2);
\draw [<-,red, in=0, out=120,  very thick] (0,-3.4) to node[right, xshift=9, yshift=-1]{$\tilde{f}_0$} (-2,-2.2);
\draw [<-,red, in=0, out=120,  very thick] (2,-4.2) to node[right, xshift=9, yshift=-1]{$\tilde{f}_0$} (0,-3.4);
\draw [<-,red, in=0, out=120,  very thick] (4,-5) to node[right, xshift=9, yshift=-1]{$\tilde{f}_0$} (2,-4.2);
\draw [<-,blue, out=30, in=150, very thick] (0,-3) to node[left, yshift=8]{$\tilde{f}_1$} (2,-3);
\draw [<-,blue, out=30, in=150, very thick] (-2,-3) to node[left, yshift=8]{$\tilde{f}_1$} (0,-3);
\draw [<-,blue, out=30, in=150, very thick] (0,-3.2) to node[left, yshift=8]{$\tilde{f}_1$} (2,-3.2);
\draw [<-,blue, out=30, in=150, very thick] (-2,-3.2) to node[left, yshift=8]{$\tilde{f}_1$} (0,-3.2);
%

\end{tikzpicture}
\\
\hbox{Initial portion of the crystal graph of $B(\Lambda_0)$ for $\widehat{\fsl}_2$}
\end{matrix}
$$
The characters of the first few Demazure modules in $L(\Lambda_0)$ are
\begin{align*}
\mathrm{char}(L(\Lambda_0)_{\le 1}) &= X^{\Lambda_0} 
= X^{\Lambda_0} \tilde E_0(q^{-1},0),
\\
\mathrm{char}(L(\Lambda_0)_{\le s_0}) 
&= D_0 X^{\Lambda_0}
= X^{\Lambda_0}(1+qX^{2\omega_1}) 
=qX^{\Lambda_0}(X^{2\omega_1}+q^{-1})
= qX^{\Lambda_0} \tilde E_{2\omega_1}(q^{-1},0),
\\
\mathrm{char}(L(\Lambda_0)_{\le s_1s_0}) &= D_1D_0 X^{\Lambda_0}
= X^{\Lambda_0}(1+qX^{2\omega_1}+q+qX^{-2\omega_1}) \\
&= qX^{\Lambda_0}(X^{-2\omega_1}+X^{2\omega_1}+1+q^{-1})
= qX^{\Lambda_0}\tilde E_{-2\omega_1}(q^{-1},0).
\end{align*}

\noindent
The crystal graph of $B(\omega_1+\Lambda_0)$ is pictured in Plate B and
\begin{align*}
\mathrm{char}(L(\omega_1+\Lambda_0)_{\le 1}) 
&= X^{\Lambda_0}X^{\omega_1} = X^{\Lambda_0} \tilde E_{\omega_1}(q^{-1},0),
\\
\mathrm{char}(L(\omega_1+\Lambda_0)_{\le s_1}) 
&= D_1 X^{\omega_1+\Lambda_0}
= X^{\Lambda_0}(X^{\omega_1}+X^{-\omega_1}) = X^{\Lambda_0} \tilde E_{-\omega_1}(q^{-1},0),
\\
\mathrm{char}(L(\omega_1+\Lambda_0)_{\le s_0s_1}) 
&= D_0D_1 X^{\omega_1+\Lambda_0}
= X^{\Lambda_0}(X^{\omega_1}+X^{-\omega_1}+qX^{\omega_1}+q^2X^{3\omega_1}) \\
&= q^2X^{\Lambda_0} \tilde E_{3\omega_1}(q^{-1},0),
\end{align*}

\subsection{The crystal $B(\omega_1)$}

The crystal
$B(\omega_1)= \{ p_{-\omega_1+k\delta}, p_{\omega_1+k\delta}\ |\ k\in \ZZ\}$
is a single connected component.  The crystal graph of $B(\omega_1)$ is pictured in Plate B.
Following \cite[(12.1.9)]{Kac} and putting $q=X^{-\delta}$
noting that $\tilde E_{-\omega_1}(q^{-1},0) = X^{\omega_1}+X^{-\omega_1}$
and $\tilde E_{-\omega_1}(q^{-1},\infty) = X^{-\omega_1}+q X^{\omega_1}$, then
\begin{align*}
\mathrm{char}(L(\omega_1)_{\le s_1}) 
&= \frac{1}{1-q^{-1}} \tilde E_{-\omega_1}(q^{-1},0)
\qquad  \hbox{and} \\ 
\mathrm{char}(L(\omega_1) )
&= 0_q \tilde E_{-\omega_1}(q^{-1},0)
= 0_q \tilde E_{-\omega_1}(q^{-1},\infty),
\end{align*}
where
$0_q = \cdots + q^{-3}+q^{-2}+q^{-1}+1+q+q^2+\cdots$ as in Remark \ref{RGchar}.

\subsection{The crystal $B(2\omega_1)$}


On crystals, the injective $\mathbf{U}$-module homomorphism
$L(2\omega_1) \hookrightarrow L(\omega_1)\otimes L(\omega_1)$
given in \eqref{tensorinjection} is the inclusion
$$
\begin{array}{cl}
B(\omega_1)\otimes B(\omega_1) 
&= \{ p_{w_1\omega_1}\otimes p_{w_2\omega_1}\ |\  w_1, w_2\in W^{\mathrm{ad}}\} \\
\cup\vert \\
B(2\omega_1) 
&= \{ p_{w_1\omega_1}\otimes p_{w_2\omega_1}\ |\  
\hbox{$w_1, w_2\in W^{\mathrm{ad}}$ with $w_1\zerogeq\ w_2$} \}
\end{array}
$$
The connected components of $B(2\omega_1)$ are determined by
$$B(2\omega_1) = B(2\omega_1)_0\times S^{2\omega_1},
\qquad\hbox{where}\quad
S^{2\omega_1} = \{ \hbox{partitions $\kappa$ with $\ell(\kappa)<2$} \}
=\ZZ_{\ge 0}.$$ 
For $\kappa\in \ZZ_{\ge 0}$, the connected component corresponding to $\kappa$,
as a subset of $B(\omega_1)\otimes B(\omega_1)$, is
$$B(2\omega_1)_\kappa = \{ p_{w_1\omega_1}\otimes p_{w_2\omega_1}\ |\ 
\hbox{$w_1\zerogeq\ w_2$ and $\ell^0(w_1)-\ell^0(w_2)\in \{\kappa, \kappa+1\}$\ } \}
$$
Representatives of the components and
the crystal graph of the connected component $B(2\omega_1)_0$
are pictured in Plate C.
Inspection of the crystal graphs gives
\begin{align*}
\mathrm{char}(L(2\omega_1)_{\le s_1} )
&=  \frac{1}{1-q^{-2}} \frac{1}{1-q^{-1}}  \tilde E_{-2\omega_1}(q^{-1},0), 
\quad\hbox{and}\\
\mathrm{char}(L(2\omega_1)) 
&= 0_q\frac{1}{1-q} \tilde E_{-2\omega_1}(q^{-1},0)
= 0_q \frac{1}{1-q} \tilde E_{-2\omega_1}(q^{-1},\infty).
\end{align*}

\section{Examples for $\fg = \widehat{\fsl}_3$}

Let $\fg = \widehat{\fsl}_3$.  The affine Dynkin diagram is 
$
\begin{tikzpicture}
	\draw (0,0)--(1,0);
	\draw (1,0)--(0.5, 0.5);
	\draw (0.5,0.5)--(0,0);
	\draw[fill=white] (0,0) circle (2.5pt) ;
	\draw[fill=white] (1,0) circle (2.5pt) ; 
         \draw[fill=white] (0.5,0.5) circle (2.5pt) ; 
\end{tikzpicture}
$ 
and the affine Cartan matrix is
$$\begin{pmatrix}
\alpha_0(h_0) &\alpha_0(h_1) &\alpha_0(h_2) \\
\alpha_1(h_0) &\alpha_1(h_1) &\alpha_1(h_2) \\
\alpha_2(h_0) &\alpha_2(h_1) &\alpha_2(h_2) 
\end{pmatrix}
= \begin{pmatrix} 2 &-1 &-1 \\ -1 &2 &-1 \\ -1 &-1 &2\end{pmatrix}
$$
and
$
h_0 = -(h_1+h_2)+K,\quad
\alpha_0 = -(\alpha_1+\alpha_2) + \delta, \quad
\Lambda_1 = \omega_1+\Lambda_0,\quad
\Lambda_2 = \omega_2+\Lambda_0.
$

In the basis $\{ \delta, \omega_1, \omega_2, \Lambda_0\}$ of $\fh^*$
the action of the affine Weyl group $W^{\mathrm{ad}}$ is given by
$$
s_1 = \begin{pmatrix} 1 &0 &0 &0 \\
0 &-1 &0 &0 \\
0 &1 &1 &0 \\
0 &0 &0 &1
\end{pmatrix}
\qquad
s_2 = \begin{pmatrix} 1 &0 &0 &0 \\
0 &1 &1 &0 \\
0 &0 &-1 &0 \\
0 &0 &0 &1
\end{pmatrix}
\qquad
s_0 = \begin{pmatrix} 1 &1 &1 &-1 \\
0 &0 &-1 &1 \\
0 &-1 &0 &1 \\
0 &0 &0 &1 \end{pmatrix},
\quad\hbox{and}
$$
$$t_{k_1h_1+k_2h_2} 
= \begin{pmatrix} 
1 &-k_1 &-k_2 &-k_1^2-k_2^2+k_1k_2 \\
0 &1 &0 &2k_1-k_2 \\
0 &0 &1 &2k_2-k_1 \\
0 &0 &0 &1 \end{pmatrix}
\qquad\hbox{for $k_1, k_2\in \ZZ$.}
$$
These matrices are used to compute the orbits pictured at the end of Section \ref{Worbits}.

\subsection{The extremal weight modules $L(\omega_1)$ and $L(\omega_2)$}

Letting $\CC^3 = \hbox{$\CC$-span}\{v_1, v_2, v_3\}$ be the standard representation of 
$\mathring{\fg} = \fsl_3$, the extremal weight modules $L(\omega_1)$ and $L(\omega_2)$ 
for $\fg = \widehat{\fsl}_3$ are
\begin{align*}
L(\omega_1) = \CC^3\otimes_\CC \CC[\epsilon, \epsilon^{-1}]
\qquad\hbox{and}\qquad
L(\omega_2) = (\Lambda^2\CC^3)\otimes_\CC \CC[\epsilon, \epsilon^{-1}],
\end{align*}
with $u_{\omega_1} = v_1$ and $u_{\omega_2} = v_1\wedge v_2$, respectively.
The crystals $B(\omega_1)$ and $B(\omega_2)$ have realizations as sets of straight line paths:
\begin{align*}
B(\omega_1) 
&= \{ p_{\omega_1+k\delta}, p_{s_1\omega_1+k\delta}, p_{-\omega_2+k\delta}\ |\  k\in \ZZ\} ,
\qquad
B(\omega_2) 
&= \{ p_{\omega_2+k\delta}, p_{s_2\omega_2+k\delta}, p_{-\omega_1+k\delta}\ |\  k\in \ZZ\}.
\end{align*}
$$
\begin{matrix}
\resizebox{5cm}{!}{
\begin{tikzpicture}[scale=1.1,every node/.style={minimum size=1cm}]
		\begin{scope}[
		yshift=0,
		every node/.append style={yslant=\yslant,xslant=\xslant},
		yslant=\yslant,xslant=\xslant
	]
		\draw [fill=black]
			(1.5,-1.25) circle (.1) node[anchor=north west] {$\omega_1+\delta$}
			(-1.5,-1.25) circle (.1) node[anchor=north west] {$-\omega_2+\delta$}
			(0,1.25) circle (.1) node[anchor=south east] {$s_1\omega_1+\delta$};
			\draw[-latex, out=90, in=0, thick, blue]
			(1.5,-1.25) to node[right]{$\tilde{f}_1$}  (0.2,1.25);
			\draw[-latex, out=-150, in=180, thick, blue]
			(0,1.25) to node[left]{$\tilde{f}_2$}  (-1.7,-1.25);
			\coordinate (-omega2+delta) at (-1.5,-1.25);
			\coordinate (omega1+delta) at (1.5,-1.25);
			\coordinate (s1omega1+delta) at (0,1.25);
	\end{scope}

		\begin{scope}[
		yshift=-100,
		every node/.append style={yslant=\yslant,xslant=\xslant},
		yslant=\yslant,xslant=\xslant
	]
		\draw [fill=black]
			(1.5,-1.25) circle (.1) node[anchor=north west] {$\omega_1$}
			(-1.5,-1.25) circle (.1) node[anchor=north west] {$-\omega_2$}
			(0,1.25) circle (.1) node[anchor=south east] {$s_1\omega_1$};
			\draw[-latex, out=90, in=0, thick, blue]
			(1.5,-1.25) to node[right]{$\tilde{f}_1$}  (0.2,1.25);
			\draw[-latex, out=-150, in=180, thick, blue]
			(0,1.25) to node[left]{$\tilde{f}_2$}  (-1.7,-1.25);
			\coordinate (-omega2) at (-1.5,-1.25);
			\coordinate (omega1) at (1.5,-1.25);
			\coordinate (s1omega1) at (0,1.25);
	\end{scope} 
	\begin{scope}[
		yshift=-200,
		every node/.append style={yslant=\yslant,xslant=\xslant},
		yslant=\yslant,xslant=\xslant
	] 
		\draw [fill=black]
			(1.5,-1.25) circle (.1) node[anchor=north west] {$\omega_1-\delta$}
			(-1.5,-1.25) circle (.1) node[anchor=north west] {$-\omega_2-\delta$}
			(0,1.25) circle (.1) node[anchor=south east] {$s_1\omega_1-\delta$};
			\draw[-latex, out=90, in=0, thick, blue]
			(1.5,-1.25) to node[right]{$\tilde{f}_1$}  (0.2,1.25);
			\draw[-latex, out=-150, in=180, thick, blue]
			(0,1.25) to node[left]{$\tilde{f}_2$}  (-1.7,-1.25);
			\coordinate (-omega2-delta) at (-1.5,-1.25);
			\coordinate (omega1-delta) at (1.5,-1.25);
			\coordinate (s1omega1-delta) at (0,1.25);
		\end{scope}
		\draw[-latex, out=-30, in=-150, thick, red] (-omega2+delta) to node[below]{$\tilde{f}_0$} (omega1);
		\draw[-latex, out=-30, in=-150, thick, red] (-omega2) to node[below]{$\tilde{f}_0$} (omega1-delta);
			
		\draw[thick, black] (-omega2+delta) to (-omega2-delta);
		\draw[thick, black] (s1omega1+delta) to (s1omega1-delta);
		\draw[thick, black] (omega1+delta) to (omega1-delta);
		
		\draw[thick,black] (-omega2+delta) to (omega1+delta);
		\draw[thick,black] (-omega2+delta) to (s1omega1+delta);
		\draw[thick,black] (omega1+delta) to (s1omega1+delta);
		
		\draw[thick,black] (-omega2) to (omega1);
		\draw[thick,black] (-omega2) to (s1omega1);
		\draw[thick,black,dashed] (omega1) to (s1omega1);
		
		\draw[thick,black] (-omega2-delta) to (omega1-delta);
		\draw[thick,black] (-omega2-delta) to (s1omega1-delta);
		\draw[thick,black,dashed] (omega1-delta) to (s1omega1-delta);
			
\end{tikzpicture}
}
\qquad\qquad
&\resizebox{5cm}{!}{
\begin{tikzpicture}[scale=1.1,every node/.style={minimum size=1cm}]
		\begin{scope}[
		yshift=0,
		every node/.append style={yslant=\yslant,xslant=\xslant},
		yslant=\yslant,xslant=\xslant
	]
		\draw [fill=black]
			(-1.5,1.25) circle (.1) node[anchor=south east] {$-\omega_1+\delta$}
			(0,-1.25) circle (.1) node[anchor=north west] {$s_2\omega_2+\delta$}
			(1.5,1.25) circle (.1) node[anchor=north west] {$\omega_2+\delta$};
			\draw[-latex, out=-90, in=30, thick, blue]
			(1.5,1.25) to node[right]{$\tilde{f}_2$}  (0,-1.25);
			\draw[-latex, out=-150, in=180, thick, blue]
			(0,-1.25) to node[left]{$\tilde{f}_1$}  (-1.5,1.25);
			\coordinate (-omega1+delta) at (-1.5,1.25);
			\coordinate (omega2+delta) at (1.5,1.25);
			\coordinate (s2omega2+delta) at (0,-1.25);
	\end{scope}

		\begin{scope}[
		yshift=-100,
		every node/.append style={yslant=\yslant,xslant=\xslant},
		yslant=\yslant,xslant=\xslant
	]
		\draw [fill=black]
			(-1.5,1.25) circle (.1) node[anchor=south east] {$-\omega_1$}
			(0,-1.25) circle (.1) node[anchor=north west] {$s_2\omega_2$}
			(1.5,1.25) circle (.1) node[anchor=north west] {$\omega_2$};
			\draw[-latex, out=-90, in=30, thick, blue]
			(1.5,1.25) to node[right]{$\tilde{f}_2$}  (0,-1.25);
			\draw[-latex, out=-150, in=180, thick, blue]
			(0,-1.25) to node[left]{$\tilde{f}_1$}  (-1.5,1.25);
			\coordinate (omega2) at (1.5,1.25);
			\coordinate (-omega1) at (-1.5,1.25);
			\coordinate (s2omega2) at (0,-1.25);
	\end{scope} 
	\begin{scope}[
		yshift=-200,
		every node/.append style={yslant=\yslant,xslant=\xslant},
		yslant=\yslant,xslant=\xslant
	] 
		\draw [fill=black]
			(-1.5,1.25) circle (.1) node[anchor=south east] {$-\omega_1-\delta$}
			(0,-1.25) circle (.1) node[anchor=north west] {$s_2\omega_2-\delta$}
			(1.5,1.25) circle (.1) node[anchor=north west] {$\omega_2-\delta$};
			\draw[-latex, out=-90, in=30, thick, blue]
			(1.5,1.25) to node[right]{$\tilde{f}_2$}  (0,-1.25);
			\draw[-latex, out=-150, in=180, thick, blue]
			(0,-1.25) to node[left]{$\tilde{f}_1$}  (-1.5,1.25);
			\coordinate (omega2-delta) at (1.5,1.25);
			\coordinate (s2omega2-delta) at (0,-1.25);
			\coordinate (-omega1-delta) at (-1.5,1.25);
		\end{scope}
		
		\draw[-latex, out=-30, in=130, thick, red] (-omega1+delta) to node[below]{$\tilde{f}_0$} (omega2);
		\draw[-latex, out=-30, in=130, thick, red] (-omega1) to node[below]{$\tilde{f}_0$} (omega2-delta);
		
		\draw[thick, black] (-omega1+delta) to (-omega1-delta);
		\draw[thick, black] (s2omega2+delta) to (s2omega2-delta);
		\draw[thick, black, dashed] (omega2+delta) to (omega2-delta);
		
		\draw[thick,black] (-omega1+delta) to (omega2+delta);
		\draw[thick,black] (-omega1+delta) to (s2omega2+delta);
		\draw[thick,black] (omega2+delta) to (s2omega2+delta);
		
		\draw[thick,black, dashed] (-omega1) to (omega2);
		\draw[thick,black] (-omega1) to (s2omega2);
		\draw[thick,black] (omega2) to (s2omega2);
		
		\draw[thick,black,dashed] (-omega1-delta) to (omega2-delta);
		\draw[thick,black] (-omega1-delta) to (s2omega2-delta);
		\draw[thick,black,dashed] (omega2-delta) to (s2omega2-delta);
\end{tikzpicture}
}
\\
\hbox{$B(\omega_1)$ crystal graph}
&\hbox{$B(\omega_2)$ crystal graph}
\end{matrix}
$$
Each of the crystals $B(\omega_1)$ and $B(\omega_2)$ has
a single connected component, all weight spaces are one dimensional and
\begin{align*}
\mathrm{char}(L(\omega_1)_{\le s_2s_1} ) &=
\frac{1}{1-q^{-1}}(X^{\omega_1} + X^{s_1\omega_1} + X^{-\omega_2}  )
=\frac{1}{1-q^{-1}} \tilde E_{-\omega_2}(q^{-1},0), \quad\hbox{and} \\
\mathrm{char}(L(\omega_1)) 
&= 0_q (X^{\omega_1} + X^{-\omega_2} + X^{s_1\omega_1})
=0_q \tilde E_{-\omega_2}(q^{-1},0)  \\
&= 0_q (q^{-1}X^{\omega_1} + q^{-1}X^{\omega_2-\omega_1} + X^{-\omega_2}  )
= 0_q \tilde E_{-\omega_2}(q^{-1},\infty).
\end{align*}
where
$0_q = \cdots + q^{-3}+q^{-2}+q^{-1}+1+q+q^2+\cdots$ as in Remark \ref{RGchar}.

For $a\in \CC$, the crystals of $M^{\mathrm{fin}}((u-a)\omega_1) \cong \CC^3$ 
and $M^{\mathrm{fin}}((u-a)\omega_2) \cong \Lambda^2(\CC^3)$ have crystal graphs $B^{\mathrm{fin}}(\omega_1)$
and $B^{\mathrm{fin}}(\omega_2)$.
$$
\begin{matrix}
\resizebox{4cm}{!}{
\begin{tikzpicture}[scale=0.95,every node/.style={minimum size=1cm}]
\tikzstyle{every node}=[font=\small]
	 \coordinate (omega2) at (-0.5,0.86603);
	 \coordinate (2omega2) at  (-1,2*0.86603);
	 \coordinate (2omega2+omega1) at (-0.5, 3*0.86603);
	 \coordinate (omega2+2omega1) at (0.5,3*0.86603);
	 \coordinate (2omega1) at (1,2*0.86603);
	 \coordinate (omega1) at (0.5,0.86603);
	 \coordinate (-omega1) at (-0.5,-0.86603);
	 \coordinate (-omega1+omega2) at (-1,0);
	 \coordinate (-omega2+omega1) at (1,0);
	 \coordinate (-omega2) at  (0.5,-0.86603);
	 \coordinate (-2omega1-omega2) at (-0.5,-3*0.86603);
	 \coordinate (-2omega1) at (-1, -2*0.86603);
	 \coordinate (-2omega2) at (1,-2*0.86603);
	 \coordinate (-2omega2-omega1) at (0.5,-3*0.86603);
	 
	\draw[thin] ($1.5*(-omega1)$) to  ($1.5*(omega1)$);
	\draw[thin] ($1.5*(-omega2)$) to ($1.5*(omega2)$);
	\draw[thin] ($1.5*(-omega1+omega2)$) to ($1.5*(-omega2+omega1)$);
	
	\draw[fill=black]
	(omega1) circle (.04) node[anchor=south west, xshift=4, yshift=-7] {$\omega_1$}
	(-omega1+omega2) circle (.04) node[anchor=north east, yshift=0] {$s_1\omega_1$}
	(-omega2) circle (.04) node[anchor=north west, yshift=7] {$-\omega_2$};

	\draw[-latex, out=160, in=60, thin, blue]  (omega1) to node[above, xshift=3]{$\tilde{f}_1$} (-omega1+omega2);
	\draw[-latex, out=-60, in=180, thin, blue]  (-omega1+omega2) to node[below, xshift=5, yshift=-2]{$\tilde{f}_2$} (-omega2);
	\draw[-latex, out=40, in=-40, thin, red]  (-omega2) to node[right, yshift=8]{$\tilde{f}_0$} (omega1);
\end{tikzpicture}
}
\qquad\qquad
&
\resizebox{4cm}{!}{
\begin{tikzpicture}[scale=0.95,every node/.style={minimum size=1cm}]
\tikzstyle{every node}=[font=\small]
	 \coordinate (omega2) at (-0.5,0.86603);
	 \coordinate (2omega2) at  (-1,2*0.86603);
	 \coordinate (2omega2+omega1) at (-0.5, 3*0.86603);
	 \coordinate (omega2+2omega1) at (0.5,3*0.86603);
	 \coordinate (2omega1) at (1,2*0.86603);
	 \coordinate (omega1) at (0.5,0.86603);
	 \coordinate (-omega1) at (-0.5,-0.86603);
	 \coordinate (-omega1+omega2) at (-1,0);
	 \coordinate (-omega2+omega1) at (1,0);
	 \coordinate (-omega2) at  (0.5,-0.86603);
	 \coordinate (-2omega1-omega2) at (-0.5,-3*0.86603);
	 \coordinate (-2omega1) at (-1, -2*0.86603);
	 \coordinate (-2omega2) at (1,-2*0.86603);
	 \coordinate (-2omega2-omega1) at (0.5,-3*0.86603);
	 
	\draw[thin] ($1.5*(-omega1)$) to  ($1.5*(omega1)$);
	\draw[thin] ($1.5*(-omega2)$) to ($1.5*(omega2)$);
	\draw[thin] ($1.5*(-omega1+omega2)$) to ($1.5*(-omega2+omega1)$);
	
	\draw[fill=black]
	(omega2) circle (.04) node[anchor=south east, xshift=0, yshift=-7] {$\omega_2$}
	(-omega2+omega1) circle (.04) node[anchor=south west, yshift=3] {$s_2\omega_2$}
	(-omega1) circle (.04) node[anchor=north west,xshift=-5 ,yshift=0] {$-\omega_1$};

	\draw[-latex, out=0, in=120, thin, blue]  (omega2) to node[above, xshift=-3, yshift=2]{$\tilde{f}_2$} (-omega2+omega1);
	\draw[-latex, out=-90, in=0, thin, blue]  (-omega2+omega1) to node[below, xshift=10, yshift=7]{$\tilde{f}_1$} (-omega1);
	\draw[-latex, out=120, in=-120, thin, red]  (-omega1) to node[left, yshift=-8]{$\tilde{f}_0$} (omega2);
\end{tikzpicture}
}
\\
\hbox{$B^{\mathrm{fin}}(\omega_1)$}
&\hbox{$B^{\mathrm{fin}}(\omega_2)$} 
\end{matrix}
$$

\subsection{The extremal weight module $L(\omega_1+\omega_2)$}

To construct the crystal $B(\omega_1+\omega_2)$ use 
$$B(\omega_1)\otimes B(\omega_2) 
= \{ p_{v\omega_1+k\delta}\otimes p_{w\omega_2+\ell\delta}
\ |\ v\in \{1, s_1, s_2s_1\}, w\in \{1, s_2, s_1s_2\}, k\in \ZZ, \ell\in \ZZ\}$$
with the tensor product action for crystals given by (see, for example, \cite[Prop.\ 5.7]{Ra06})
$$\tilde f_i(p_1\otimes p_2) = \begin{cases}
\tilde f_i p_1 \otimes p_2, &\hbox{if $d_i^+(p_1)>d_i^-(p_2)$,} \\
p_1 \otimes \tilde f_i p_2, &\hbox{if $d_i^+(p_1)\le d_i^-(p_2)$,} 
\end{cases}
\quad
\tilde e_i(p_1\otimes p_2) = \begin{cases}
\tilde e_i p_1 \otimes p_2, &\hbox{if $d_i^+(p_1)\ge d_i^-(p_2)$,} \\
p_1 \otimes \tilde e_i p_2, &\hbox{if $d_i^+(p_1)< d_i^-(p_2)$.} 
\end{cases}
$$
where $d_i^{\pm}(p)$ are determined by
$$\tilde f_i^{d_i^+(p)} p \ne 0 \ \hbox{and}\ 
\tilde f_i^{d_i^+(p)+1} p = 0 ,
\qquad\hbox{and}\qquad
\tilde e_i^{d_i^-(p)} p \ne 0 \ \hbox{and}\ 
\tilde e_i^{d_i^-(p)+1} p = 0 .
$$
The crystal $B(\omega_1+\omega_2)$ is realized as a subset of 
$B(\omega_1)\otimes B(\omega_2)$ via the crystal embedding
$$
\begin{matrix}
B(\omega_1+\omega_2) &\hookrightarrow &B(\omega_1)\otimes B(\omega_2) \\
b_{\omega_1+\omega_2} &\longmapsto &p_{\omega_1}\otimes p_{\omega_2}
\end{matrix}
$$
By Theorem \ref{levelzerostructure}, $B(\omega_1+\omega_2)$ is
connected and is generated by $b_{\omega_1+\omega_2}$.  The crystal
$B(\omega_1+\omega_2)$ is pictured and its character is computed in Plate D.

\section{Alcove walks for affine flag varieties}\label{affineflags}

\subsection{The affine Kac-Moody group $G$}

The most visible form of the loop group is  $G = \mathring{G}(\CC((\epsilon)))$, 
where $\CC((\epsilon))$ is the field of formal power series in a variable $\epsilon$
and $\mathring{G}$ is a reductive algebraic group.  The favourite example is 
when
$$\hbox{when}\quad\mathring{G} = GL_n
\qquad\hbox{and the loop group is}\quad
G=GL_n(\CC((\epsilon))),
$$
the group of $n\times n$ invertible matrices 
with entries in $\CC((\epsilon))$.
A slightly more extended, and extremely powerful, point of view is to let 
$G$ be the Kac-Moody group whose Lie algebra is the affine Lie algebra
$\fg$ of Section \ref{affLiealgebra}, so that 
$G$ is a central extension of a semidirect product (where the semidirect
product comes from the action of $\CC^\times$ on $\mathring{G}(\CC((\epsilon))$ 
by ``loop rotations")
\begin{equation}
\{1\} \to \CC^\times \to G \to \mathring{G}(\CC((\epsilon)))\rtimes \CC^\times \to \{1\}
\qquad\hbox{so that}\qquad
G=\exp(\fg).
\label{affKMgroup}
\end{equation}

In this section we wish to work with $G$ via generators and relations.
Fortunately, presentations of $G$ are well established in the work of Steinberg
\cite{St67} and Tits \cite{Ti87} and others (see \cite[\S3]{PRS} for a survey).
More specifically, up to an extra (commutative) torus $T\cong (\CC^\times)^k$ for some $k$,
the group $G$ is generated by subgroups isomorphic to $SL_2(\CC)$,
the images of homomorphisms
\begin{equation}
\begin{matrix}
\varphi_i \colon &SL_2(\CC) &\longrightarrow &G \\
&\begin{pmatrix} 1 &c \\ 0 &1\end{pmatrix} &\longmapsto &x_{\alpha_i}(c) \\
&\begin{pmatrix} 1 &0 \\ c &1\end{pmatrix} &\longmapsto &x_{-\alpha_i}(c) \\
&\begin{pmatrix} d &0 \\ 0 &d^{-1}\end{pmatrix} &\longmapsto &h_{\alpha_i^\vee}(d) \\
&\begin{pmatrix} c &1 \\ -1 &0\end{pmatrix} &\longmapsto &y_i(c)
\end{matrix}
\qquad
\begin{array}{l}
\hbox{for each vertex}\quad
\begin{tikzpicture}
	\draw[fill=white] (0,0) circle (2.5pt) node[above=1pt] {\small $\phantom{{}_j}i\phantom{{}_j}$};
\end{tikzpicture}
\\
\hbox{of the affine Dynkin diagram.}
\end{array}
\label{yidefn}
\end{equation}
For each $\alpha\in \mathring{R}^+$ there is a homomorphism
$\varphi_{\alpha}\colon SL_2(\CC((\epsilon)))\to G$ and we let
\begin{equation}
x_\alpha(f) = \varphi_\alpha\begin{pmatrix} 1 &f \\ 0 &1\end{pmatrix} 
\qquad\hbox{and}\qquad
x_{-\alpha}(f) = \varphi_\alpha\begin{pmatrix} 1 &0 \\ f &1\end{pmatrix}
\quad\hbox{for $f\in \CC((\epsilon))$.}
\label{xalphadefn}
\end{equation}
For each $z\in W^{\mathrm{ad}}$ fix a reduced word $z={s_{j_1}}\cdots s_{j_r}$,
define
\begin{equation}
n_z = y_{j_1}(0)\cdots y_{j_r}(0)
\quad\hbox{and define}\quad
x_{\pm\alpha+r\delta}(c) = x_{\pm\alpha}(c\epsilon^r)
\quad\hbox{for $\alpha\in \mathring{R}^+$ and $r\in \ZZ$.}
\label{rootgpelts}
\end{equation}

In the following sections, for simplicity, we will work only with the loop group 
$G = \mathring{G}(\CC((\epsilon)))$ rather than the full affine Kac-Moody group as in 
\eqref{affKMgroup}.
We shall also use a slightly more general setting where $\CC$ is replaced by
an arbitrary field $\kk$ so that, in Proposition \ref{Heckeactions}, we may let $\kk = \FF_q$ be the finite field
with $q$ elements.

\subsection{The affine flag varieties $G/I^+$,
$G/I^0$ and $G/I^-$}

Let $\kk$ be a field.
Define subgroups of $\mathring{G}(\kk)$ by
\begin{align*}
\hbox{``unipotent upper triangular matrices''}\qquad\qquad
U^+(\kk) &= \langle x_\alpha(c)\ |\ \alpha\in \mathring{R}^+, c\in \kk\rangle, \\
\hbox{``diagonal matrices''}\qquad\qquad
H(\kk^\times) &= \langle h_{\alpha_i^\vee}(c)\ |\ i\in \{1, \ldots, n\}, \ c\in \kk^\times \rangle, \\
\hbox{``unipotent lower triangular matrices''}\qquad\qquad
U^-(\kk) &= \langle x_{-\alpha}(c)\ |\ \alpha\in \mathring{R}^+, c\in \kk\}\rangle.
\end{align*}
Let
$$
\begin{array}{llll}
&\hbox{$\mathbb{F}=\kk((\epsilon))$,\quad which has}\quad 
&\fo=\kk[[\epsilon]] \quad\hbox{and}\quad
&\fo^\times = \{ p\in \kk[[\epsilon]]\ |\ p(0)\ne 0\}, \\
\hbox{or let} &\hbox{$\mathbb{F} = \kk[\epsilon, \epsilon^{-1}]$,\quad which has}\quad
&\fo = \kk[\epsilon]\quad\hbox{and}\quad
&\fo^\times = \kk^\times.
\end{array}
$$
Define subgroups of $G=\mathring{G}(\FF)$ by
$$U^-(\FF) = \langle x_{-\alpha+k\delta}(c)\ |\ c\in \kk, \alpha\in \mathring{R}^+, k\in \ZZ
\rangle
\qquad\hbox{and}\qquad
H(\fo^\times) = \langle h_{\alpha_i^\vee}(c)\ |\ c\in \fo^\times\rangle.$$
Let $g(0)$ denote $g$ evaluated at $\epsilon=0$ and let
$g(\infty)$ denote $g$ evaluated at $\epsilon^{-1}=0$.
Following, for example, \cite[Lect.\ 11 Theorem (B)]{Pet97} and \cite[\S 11]{LuICM}, 
define subgroups of $G(\FF)$ by
\begin{align*}
\hbox{(positive Iwahori)} \qquad\qquad 
I^+ &= \{ g\in G \ |\ \hbox{$g(0)$ exists and $g(0)\in U^+(\kk)H(\kk^\times)$} \}, \\
\hbox{(level 0 Iwahori)} \qquad\qquad 
I^0 &= U^-(\FF) H(\fo^\times), \\
\hbox{(negative Iwahori)} \qquad\qquad 
I^- &= \{ g\in G\ |\ \hbox{$g(\infty)$ exists and $g(\infty) \in U^-(\kk)H(\kk^\times)$} \},
\end{align*}
Then 
\begin{align*}
G/I^+ 
&\text{ is the positive level (thin) affine flag variety,} \\
G/I^0
&\text{ is the level 0 (semi-infinite) affine flag variety,} \\
G/I^-
&\text{ is the negative level (thick) affine flag variety.} 
\end{align*}
These are studied with the aid of the decompositions
(relax notation and write $zI^+$ for $n_zI^+$)
$$
G=\bigsqcup_{x\in W^{\mathrm{ad}}} I^+xI^+, \qquad
G=\bigsqcup_{y\in W^{\mathrm{ad}}} I^0yI^+, \qquad
G=\bigsqcup_{z\in W^{\mathrm{ad}}} I^-zI^+,
$$

\subsection{Labeling points of $I^+wI^+$, $I^+wI^+\cap I^0vI^+$ and $I^+wI^+\cap I^-zI^+$}

Recall the bijection between elements of $W^{\mathrm{ad}}$ and alcoves given in
\eqref{alcovedefn} and note that if $v\in W^{\mathrm{ad}}$ and $i\in \{0, \ldots, n\}$
then the hyperplane separating the alcoves corresponding
to $v$ and $s_iv$ is 
$$\fh^{v\alpha^\vee_i} = \{ x+\Lambda_0 \ |\ 
\hbox{$x\in \fa_\RR^*$ and $\langle x+\Lambda_0, v\alpha^\vee_i\rangle = 0$}\}.$$
A \emph{blue labeled step of type $i$} is
$$
\begin{matrix}
\begin{tikzpicture}
\node[label=below:$c$] at (0.3,0.6){};
\node[label=above:$\fh^{v\alpha^\vee_i}$] at (0, 0.8){};
\node[label=above:$v$] at (-0.7,0.35){};
\node[label=above:$vs_i$] at (0.7,0.35){};
\draw (0,0) to (0,1);
\draw[blue][->](-0.7,0.5) -- (0.7, 0.5);
\end{tikzpicture}
\end{matrix}
\qquad\hbox{with\quad $vs_i \posgeq\ v$.}
\qquad\hbox{Let}\quad
\Phi^+
\begin{pmatrix}
\begin{tikzpicture}
\node[label=below:$c$] at (0.3,0.6){};
\node[label=above:$\fh^{v\alpha^\vee_i}$] at (0, 0.8){};
\node[label=above:$v$] at (-0.7,0.35){};
\node[label=above:$vs_i$] at (0.7,0.35){};
\draw (0,0) to (0,1);
\draw[blue][->](-0.7,0.5) -- (0.7, 0.5);
\end{tikzpicture}
\end{pmatrix}
= y_i(c),
$$
where $y_i(c)$ is as in \eqref{yidefn}. 
A \textit{blue labeled walk of type $\vec w = (i_1,\ldots,i_{\ell})$} is a sequence
$(p_1, \ldots, p_\ell)$ 
where $p_k$ is a blue labeled step of type $i_k$ and 
which begins at the alcove $1$.

\smallskip\noindent
A \emph{red labeled step of type $i$} is
$$
\begin{matrix}
\begin{tikzpicture}
\node[label=below:$c$] at (0.3,0.6){};
\node[label=above:$\fh^{v\alpha^\vee_i}$] at (0, 0.8){};
\node[label=above:$v$] at (-0.7,0.35){};
\node[label=above:$vs_i$] at (0.7,0.35){};
\draw (0,0) to (0,1);
\draw[red][->](-0.7,0.5) -- (0.7, 0.5);
\end{tikzpicture} \\
c\in \kk
\end{matrix} 
\quad\hbox{or}\quad
\begin{matrix}
\begin{tikzpicture}
\node[label=below:$0$] at (-0.4,0.6){};
\node[label=above:$\fh^{-v\alpha^\vee_i}$] at (0, 0.8){};
\node[label=above:$v$] at (-0.7,0.35){};
\node[label=above:$vs_i$] at (0.7,0.35){};
\draw (0,0) to (0,1);
\draw[red][->](0.7, 0.5) -- (-0.7,0.5);
\end{tikzpicture} \\
\phantom{\CC}
\end{matrix} 
\quad\hbox{or}\quad
\begin{matrix}
\begin{tikzpicture}
\node[label=above:$v$] at (-0.7,0.35){};
\node[label=above:$vs_i$] at (0.7,0.35){};
\node[label=below:$c^{-1}$] at (0.35,0.6){};
\node[label=above:$\fh^{-v\alpha^\vee_i}$] at (0, 0.8){};
\draw (0,0) to (0,1);
\draw[red][->] (0.7, 0.4) -- (0.1,0.4) -- (0.1,0.6) --(0.7, 0.6);
\end{tikzpicture} \\
c\in \kk^\times
\end{matrix} 
\qquad\hbox{where\quad $vs_i\zerogeq\ v$.}
$$
With notations as in \eqref{rootgpelts}, let
$$
\Phi^0
\begin{pmatrix}
\begin{tikzpicture}
\node[label=above:$v$] at (-0.5,0.35){};
\node[label=above:$vs_i$] at (0.5,0.35){};
\node[label=below:$c$] at (0.3,0.6){};
\node[label=above:$\fh^{v\alpha^\vee_i}$] at (0, 0.8){};
\draw (0,0) to (0,1);
\draw[red][->](-0.7,0.5) -- (0.7, 0.5);
\end{tikzpicture}
\end{pmatrix} 
= x_{v\alpha_i}(c),
\quad
\Phi^0
\begin{pmatrix}
\begin{tikzpicture}
\node[label=above:$v$] at (-0.5,0.35){};
\node[label=above:$vs_i$] at (0.5,0.35){};
\node[label=below:$0$] at (-0.4,0.6){};
\node[label=above:$\fh^{-v\alpha^\vee_i}$] at (0, 0.8){};
\draw (0,0) to (0,1);
\draw[red][->](0.7, 0.5) -- (-0.7,0.5);
\end{tikzpicture}
\end{pmatrix} 
= x_{-v\alpha_i}(0),
\quad
\Phi^0
\begin{pmatrix}
\begin{tikzpicture}
\node[label=above:$v$] at (-0.5,0.35){};
\node[label=above:$vs_i$] at (0.5,0.35){};
\node[label=below:$c^{-1}$] at (0.35,0.6){};
\node[label=above:$\fh^{-v\alpha^\vee_i}$] at (0, 0.8){};
\draw (0,0) to (0,1);
\draw[red][->] (0.7, 0.4) -- (0.1,0.4) -- (0.1,0.6) --(0.7, 0.6);
\end{tikzpicture}
\end{pmatrix} 
= x_{-v\alpha_i}(c^{-1}).
$$
A \textit{red labeled walk of type $\vec w = (i_1,\ldots,i_{\ell})$} is a sequence
$(p_1, \ldots, p_\ell)$ 
where $p_k$ is a red labeled step of type $i_k$ and 
which begins at the alcove $1$.

\smallskip\noindent
A \emph{green labeled step of type $i$} is
$$
\begin{matrix}
\begin{tikzpicture}
\node[label=below:$c$] at (0.4,0.65){};
\node[label=above:$\fh^{v\alpha^\vee_i}$] at (0, 0.8){};
\draw (0,0) to (0,1);
\draw[teal][->](-0.7,0.5) -- (0.7, 0.5);
\node[label=above:$v$] at (-0.5,0.5){};
\node[label=above:$vs_i$] at (0.5,0.5){};
\end{tikzpicture} \\
c\in \kk
\end{matrix} 
\qquad\hbox{or}\qquad
\begin{matrix}
\begin{tikzpicture}
\node[label=below:$0$] at (0.4,0.65){};
\node[label=above:$\fh^{-v\alpha^\vee_i}$] at (0, 0.8){};
\draw (0,0) to (0,1);
\draw[teal][->](0.7,0.5) -- (-0.7, 0.5);
\node[label=above:$v$] at (-0.5,0.5){};
\node[label=above:$vs_i$] at (0.5,0.5){};
\end{tikzpicture} \\
\phantom{\kk}
\end{matrix} 
\qquad\hbox{or}\qquad
\begin{matrix}
\begin{tikzpicture}
\node[label=below:$c^{-1}$] at (0.5,0.6){};
\node[label=above:$\fh^{-v\alpha^\vee_i}$] at (0, 0.8){};
\node[label=above:$v$] at (-0.5,0.5){};
\node[label=above:$vs_i$] at (0.5,0.5){};
\draw (0,0) to (0,1);
\draw[teal][->] (0.7, 0.4) -- (0.1,0.4) -- (0.1,0.6) --(0.7, 0.6);
\end{tikzpicture} \\
c\in \kk^\times
\end{matrix} 
\qquad\hbox{where\quad $vs_i \neggeq\ v$.}
$$
Let
$$
\Phi^-
\begin{pmatrix}
\begin{tikzpicture}
\node[label=below:$c$] at (0.4,0.65){};
\node[label=above:$\fh^{v\alpha^\vee_i}$] at (0, 0.8){};
\draw (0,0) to (0,1);
\draw[teal][->](-0.7,0.5) -- (0.7, 0.5);
\node[label=above:$v$] at (0.5,0.5){};
\node[label=above:$vs_i$] at (-0.5,0.5){};
\end{tikzpicture}
\end{pmatrix} 
= x_{v\alpha_i}(c),
\quad
\Phi^-
\begin{pmatrix}
\begin{tikzpicture}
\node[label=below:$0$] at (-0.4,0.65){};
\node[label=above:$\fh^{-v\alpha^\vee_i}$] at (0, 0.8){};
\draw (0,0) to (0,1);
\draw[teal][->](0.7,0.5) -- (-0.7, 0.5);
\node[label=above:$v$] at (-0.5,0.5){};
\node[label=above:$vs_i$] at (0.5,0.5){};
\end{tikzpicture}
\end{pmatrix} 
= x_{-v\alpha_i}(0),
\quad
\Phi^-
\begin{pmatrix}
\begin{tikzpicture}
\node[label=below:$c^{-1}$] at (0.5,0.6){};
\node[label=above:$\fh^{-v\alpha^\vee_i}$] at (0, 0.8){};
\node[label=above:$v$] at (-0.5,0.5){};
\node[label=above:$vs_i$] at (0.5,0.5){};
\draw (0,0) to (0,1);
\draw[teal][->] (0.7, 0.4) -- (0.1,0.4) -- (0.1,0.6) --(0.7, 0.6);
\end{tikzpicture}
\end{pmatrix} 
= x_{-v\alpha_i}(c^{-1}).
$$
A \textit{green labeled walk of type $\vec w = (i_1,\ldots,i_{\ell})$} is a sequence
$(p_1, \ldots, p_\ell)$ 
where $p_k$ is a green labeled step of type $i_k$ and 
which begins at the alcove $1$.

\begin{thm} \label{alcovewalks}
Let $v,w\in W^{\mathrm{ad}}$ and fix a reduced expression $w=s_{i_1}\ldots s_{i_{\ell}}$ for $w$.
The maps $\Psi^+_{\vec w}$, $\Psi^0_{\vec w, v}$, $\Psi^-_{\vec w, v}$ are bijections.
\begin{equation*}
\begin{matrix}
\Phi^+_{\vec w}\colon 
&\left\{ \begin{array}{l}
\hbox{blue labeled paths of type} \\
\hbox{$\vec w = (i_1, \ldots, i_\ell)$}
\end{array} \right\}
&\stackrel{\sim}{\longrightarrow}
&(I^+ w I^+)/I^+ \\
&(p_1,\ldots, p_\ell) &\longmapsto 
&\Phi^+(p_1)\ldots \Phi^+(p_{\ell}) I^+  \\
\\
\Phi^0_{\vec w, v}\colon 
&\left\{ \begin{array}{l}
\hbox{red labeled paths of type} \\
\hbox{$\vec w = (i_1, \ldots, i_\ell)$ ending in $v$}
\end{array} \right\}
&\stackrel{\sim}{\longrightarrow}
&(I^0 v I^+ \cap I^+ w I^+)/I^+ \\
&(p_1,\ldots, p_\ell) &\longmapsto 
&\Phi^0(p_1)\ldots \Phi^0(p_{\ell})n_vI^+ \\
\\
\Phi^-_{\vec w, v}\colon 
&\left\{ \begin{array}{l}
\hbox{green labeled paths of type} \\
\hbox{$\vec w = (i_1, \ldots, i_\ell)$ ending in $v$}
\end{array} \right\}
&\stackrel{\sim}{\longrightarrow}
&(I^- v I^+ \cap I^+ w I^+)/I^+ \\
&(p_1,\ldots, p_\ell) &\longmapsto 
&\Phi^-(p_1)\ldots \Phi^-(p_{\ell})n_vI^+ 
\end{matrix}
\end{equation*}
\end{thm}

The proof of Theorem \ref{alcovewalks} 
is by induction on the length of $w$ following \cite[Theorem 4.1 and \S7]{PRS} where
(a) and (b) are proved.  The induction step 
can be formulated as the following proposition.

\begin{prop}\label{multsummary}
Let $v,w\in W^{\mathrm{ad}}$ and fix a reduced expression $\vec w=s_{i_1}\ldots s_{i_{\ell}}$ for $w$.
Let $j\in \{0, \ldots, n\}$ and $c\in \CC$.  
If $ws_j\posl\ w$ then assume that the reduced word for $w$ is chosen with $i_\ell = j$.
Let
$$\hbox{$\tilde c\in \CC$ and $\tilde b_1\in I^+$
\quad be the unique elements such that}\quad
b_1 y_j(c) = y_j(\tilde c) \tilde b_1.
$$
\item[(a)] 
Let $p=(p_1,\ldots, p_\ell)$ be a blue labeled path of type $(i_1, \ldots, i_\ell)$
and let $\Phi^+_{\vec w}(p) = y_{i_1}(c_1)\ldots y_{i_\ell}(c_\ell)$.
Then
\begin{align*}
&(y_{i_1}(c_1)\cdots y_{i_{\ell}}(c_\ell)b_1) (y_j(c)b_2) \\
&=\begin{cases}
y_{i_1}(c_1)\cdots y_{i_{\ell}}(c_\ell) y_j(\tilde c)\tilde b_1b_2, 
&\hbox{if $w\posl\  ws_j$,} \\
y_{i_1}(c_1)\cdots y_{i_{\ell-1}}(c_{\ell-1}) y_{i_\ell}(c_\ell-\tilde c^{-1})
h_{\alpha_1^\vee}(\tilde c)x_{\alpha_{i_\ell}}(-\tilde c^{-1})\tilde b_1b_2, 
&\hbox{if $ws_j\posl\  w$ and $\tilde c\ne 0$,} \\
y_{i_1}(c_1)\cdots y_{i_{\ell-1}}(c_{\ell-1})h_{\alpha_{i_\ell}^\vee}(-1)x_{\alpha_{i_\ell}}(c)
\tilde b_1 b_2, 
&\hbox{if $ws_j\posl\  w$ and $\tilde c=0$,}
\end{cases}
\end{align*}
\item[(b)] 
Let $p=(p_1,\ldots, p_\ell)$ be a red labeled path of type $(i_1, \ldots, i_\ell)$ ending in $v$
and let \hfil\break
$\Phi^0_{\vec w,v}(p) = x_{\gamma_1}(c_1)\cdots x_{\gamma_\ell}(c_\ell) n_vI^+$ where the notation is as in \eqref{rootgpelts}.
Then
\begin{align*}
&(x_{\gamma_1}(c_1)\cdots x_{\gamma_{\ell}}(c_{\ell})n_vb_1)( y_j(c) b_2) \\
&=\begin{cases}
x_{\gamma_1}(c_1)\cdots x_{\gamma_{\ell}}(c_{\ell})x_{v\al_j}(\pm \tilde{c})n_{vs_j} \tilde b_1 b_2.
&\hbox{if $w\zerol\  ws_j$,} \\
x_{\gamma_1}(c_1)\cdots x_{\gamma_{\ell}}(c_{\ell})x_{-v\al_j}(\tilde{c}^{-1})n_v
x_{\alpha_j}(-\tilde{c})h_{\alpha_j^{\vee}}(\tilde{c}) \tilde b_1 b_2,
&\hbox{if $ws_j\zerol\  w$ and $\tilde c\ne 0$,} \\
x_{\gamma_1}(c_1)\cdots x_{\gamma_{\ell}}(c_{\ell}) 
x_{-v\alpha_j}(0)n_{vs_j} \tilde b_1 b_2,
&\hbox{if $ws_j\zerol\  w$ and $\tilde c=0$,}
\end{cases}
\end{align*}
\item[(c)] 
Let $p=(p_1,\ldots, p_\ell)$ be a green labeled path of type $(i_1, \ldots, i_\ell)$ ending in $v$
and let \hfil\break
$\Phi^-_{\vec w,v}(p) = x_{\gamma_1}(c_1)\cdots x_{\gamma_\ell}(c_\ell) n_vI^+$ where the notation is as in \eqref{rootgpelts}.
Then
\begin{align*}
&(x_{\gamma_1}(c_1)\cdots x_{\gamma_{\ell}}(c_{\ell})n_vb_1)( y_j(c) b_2) \\
&=\begin{cases}
x_{\gamma_1}(c_1)\cdots x_{\gamma_{\ell}}(c_{\ell})x_{v\al_j}(\pm \tilde{c})n_{vs_j} \tilde b_1 b_2.
&\hbox{if $w\negl\  ws_j$,} \\
x_{\gamma_1}(c_1)\cdots x_{\gamma_{\ell}}(c_{\ell})x_{-v\al_j}(\tilde{c}^{-1})n_v
x_{\alpha_j}(-\tilde{c})h_{\alpha_j^{\vee}}(\tilde{c}) \tilde b_1 b_2,
&\hbox{if $ws_j\negl\  w$ and $\tilde c\ne 0$,} \\
x_{\gamma_1}(c_1)\cdots x_{\gamma_{\ell}}(c_{\ell}) 
x_{-v\alpha_j}(0)n_{vs_j} \tilde b_1 b_2,
&\hbox{if $ws_j\negl\  w$ and $\tilde c=0$,}
\end{cases}
\end{align*}
\end{prop}

\subsection{Actions of the affine Hecke algebra}

The \emph{affine flag representation} is \quad
$$C(G/I^+)=\CC \linspan \left\{y_{\vec w}(\vec c)I^+ \ |\ w\in W^{\mathrm{ad}}, \vec c = (c_1, \ldots, c_\ell)\in \CC^{\ell(w)} \right\},$$
where, for a fixed reduced word $w=s_{i_1}\cdots s_{i_\ell}$
$$y_{\vec w}(\vec c) = y_{i_1}(c_1)\cdots y_{i_\ell}(c_\ell)I^+
=\sum_{g\in y_{i_1}(c_1)\cdots y_{i_\ell}(c_\ell)I^+} g
\qquad (\hbox{a formal sum}).
$$
Let
$$
\begin{array}{lll}
C(I^+\backslash G/I^+) = \CC\linspan \{ T_w \ |\ w\in W^{\mathrm{ad}}\},
\qquad&\hbox{where}\quad
&T_w = I^+ w I^+, \\
C(I^0\backslash G/I^+) = \CC\linspan \{X^w \ |\ w\in W^{\mathrm{ad}}\},
\qquad&\hbox{where}\quad
&X^w=I^0 w I^+, \\
C(I^-\backslash G/I^+) = \CC\linspan \{L^w \ |\ w\in W^{\mathrm{ad}}\},
\qquad&\hbox{where}\quad
&L^w=I^- w I^+, 
\end{array}
$$

\medskip\noindent
The \emph{affine Hecke algebra} is \quad
$h_{I^+}(G/I^+)=\hbox{$\CC$-span} \{ T_w \ |\ w\in W \}$,
$$\hbox{where}\quad
T_w = I^+ n_w I^+ = \sum_{x\in I^+ n_w I^+} x
\qquad\qquad\hbox{(a formal sum).}
$$
The formal sums allow us to view all of these elements as elements of the group algebra
of $G$, where infinite formal sums are allowed (to do this precisely one should use Haar
measure and a convolution product).
Proposition \ref{flagrep} computes the (right) action of $h_{I^+}(G/I^+)$ on $C(G/I^+)$,
and Proposition \ref{Heckeactions} computes 
the (right) action of $h_{I^+}(G/I^+)$ on $C(I^+\backslash G/I^+)$
on $C(I^0\backslash G/I^+)$, and
on $C(I^-\backslash G/I^+)$.  
Proposition \ref{flagrep} follows from  Proposition \ref{multsummary}(a) by summing over $c$
and Proposition \ref{Heckeactions} follows from Proposition \ref{multsummary} by summing over the appropriate
double cosets.
We use the convention that the normalization (Haar measure) is
such that $I^+\cdot I^+ = I^+$.

\begin{prop} \label{flagrep}  Let $w\in W^{\mathrm{ad}}$, let $\vec w= s_{i_1}\cdots s_{i_\ell}$
be a reduced word for $w$ and let $j\in \{0, \ldots, n\}$.  Assume that if $ws_j\posl\ w$ then $i_\ell = j$  and let
$$y_{\overrightarrow{ws_j}}(\vec c) = 
y_{i_1}(c_1)\cdots y_{i_{\ell-1}}(c_{\ell-1})
\quad\hbox{and}\quad
y_{\vec w}(\vec c) 
= y_{i_1}(c_1)\cdots y_{i_{\ell}}(c_{\ell}) = y_{\overrightarrow{ws_j}}(\vec c)y_j(c_\ell).
$$
Then
$$
y_{\vec w}(\vec c)I^+\cdot T_{s_j}
=\begin{cases}
\displaystyle{
y_{i_1}(c_1) \cdots y_{i_{\ell-1}}(c_{\ell-1})I^+ 
+ \sum_{\tilde{c}\in \kk^{\times}}y_{i_1}(c_1) \cdots y_{i_{\ell}}(c_{\ell}-\tilde{c}^{-1})I^+, }
&\hbox{if $ws_j \posl\  w$,} \\
\displaystyle{
\sum_{\tilde{c}\in \kk}y_{i_1}(c_1)\cdots y_{i_{\ell}}(c_{\ell})y_j(\tilde{c})I^+, }
&\hbox{if $w \posl\ ws_j$.} 
\end{cases}
$$
\end{prop}

\begin{prop} \label{Heckeactions} Let $\kk=\FF_q$ the finite field with $q$ elements.
Let $w\in W^{\mathrm{ad}}$ and $j\in \{0, \ldots, n\}$.   Then
$$
T_wT_{s_j}=\begin{cases}
T_{ws_j}, &\hbox{if $w\posl\  ws_j$,} \\
(q-1)T_w+qT_{ws_j}. &\hbox{if $ws_j\posl\ w$,}
\end{cases}
\qquad
L^wT_{s_j} = \begin{cases}
L^{ws_j}, &\hbox{if $w \negl\  ws_j$,} \\
qL^{ws_j}+(q-1)L^w, &\hbox{if $ws_j \negl\ w$,}
\end{cases}
$$
$$
\hbox{and}\qquad X^wT_{s_j} = \begin{cases}
X^{ws_j}, &\hbox{if $w \zerol\  ws_j$,} \\
qX^{ws_j}+(q-1)X^w, &\hbox{if $ws_j \zerol\  w$.}
\end{cases}
$$
\end{prop}

\end{document}